\documentclass[10pt]{article}
\usepackage[english]{babel}
\usepackage{hyperref}
\usepackage{amssymb, amsfonts}
\usepackage{stmaryrd}
\usepackage{mathtools}

\usepackage{amsmath}
\usepackage{latexsym}
\usepackage{amssymb}

\usepackage{amsthm}

\font\tengoth=eufm10 at 10pt
\font\sevengoth=eufm7 at 6pt
\newfam\gothfam
\textfont\gothfam=\tengoth
\scriptfont\gothfam=\sevengoth

\newcommand{\mlabel}[1]{\marginpar{#1}\label{#1}}

\newcommand{\fB}{{\mathfrak B}}

\newcommand{\fS}{{\mathfrak S}}

\newcommand{\s}{{\mathfrak s}}

\renewcommand{\:}{\colon}
\newcommand{\1}{\mathbf{1}}

\newcommand{\cA}{\mathcal{A}}
\newcommand{\cB}{\mathcal{B}}
\newcommand{\cC}{\mathcal{C}}

\newcommand{\cF}{\mathcal{F}}

\newcommand{\cH}{\mathcal{H}}
\newcommand{\cK}{\mathcal{K}}
\newcommand{\cL}{\mathcal{L}}
\newcommand{\cM}{\mathcal{M}}
\newcommand{\cN}{\mathcal{N}}

\newcommand{\cS}{\mathcal{S}}
\newcommand{\cT}{\mathcal{T}}

\newcommand{\eset}{\emptyset}

\renewcommand{\phi}{\varphi}
\newcommand{\dd}{{\tt d}}

\newcommand{\trile}{\trianglelefteq}
\newcommand{\subeq}{\subseteq}
\newcommand{\supeq}{\supseteq}

\newcommand{\into}{\hookrightarrow}
\newcommand{\eps}{\varepsilon}

\newcommand{\shalf}{{\textstyle{\frac{1}{2}}}}

\newcommand{\N}{{\mathbb N}}
\newcommand{\Z}{{\mathbb Z}}
\newcommand{\R}{{\mathbb R}}
\newcommand{\C}{{\mathbb C}}

\newcommand{\T}{{\mathbb T}}

\newcommand{\bS}{{\mathbb S}}

\renewcommand{\hat}{\widehat}

\renewcommand{\tilde}{\widetilde}

\newcommand{\SO}{\mathop{{\rm SO}}\nolimits}

\newcommand{\U}{\mathop{\rm U{}}\nolimits}

\newcommand{\Sp}{\mathop{{\rm Sp}}\nolimits}

\newcommand{\Ad}{\mathop{{\rm Ad}}\nolimits}

\newcommand{\Hom}{\mathop{{\rm Hom}}\nolimits}

\newcommand{\Rep}{\mathop{{\rm Rep}}\nolimits}

\newcommand{\Heis}{\mathop{{\rm Heis}}\nolimits}

\newcommand{\cI}{{\mathcal I}}

\newcommand{\Aut}{\mathop{{\rm Aut}}\nolimits}

\newcommand{\id}{\mathop{{\rm id}}\nolimits}

\renewcommand{\dim}{\mathop{{\rm dim}}\nolimits}

\newcommand{\Spann}{\mathop{{\rm span}}\nolimits}

\newcommand{\Rarrow}{\Rightarrow}

\newcommand{\oline}{\overline}

\newcommand{\la}{\langle}
\newcommand{\ra}{\rangle}

\newcommand{\res}{\vert}

\newcommand{\Spec}{{\rm Spec}}

\newcommand{\ssssarr}{\hbox to 15pt{\rightarrowfill}}
\newcommand{\sssarr}{\hbox to 20pt{\rightarrowfill}}
\newcommand{\ssarr}{\hbox to 30pt{\rightarrowfill}}
\newcommand{\sarr}{\hbox to 40pt{\rightarrowfill}}
\newcommand{\arr}{\hbox to 60pt{\rightarrowfill}}
\newcommand{\larr}{\hbox to 60pt{\leftarrowfill}}
\newcommand{\Arr}{\hbox to 80pt{\rightarrowfill}}

\def\theoremname{Theorem}
\def\propositionname{Proposition}
\def\corollaryname{Corollary}
\def\lemmaname{Lemma}
\def\remarkname{Remark}
\def\conjecturename{Conjecture}

\def\definitionname{Definition}
\def\exercisename{Exercise}
\def\examplename{Example}
\def\examplesname{Examples}
\def\problemname{Problem}
\def\problemsname{Problems}
\def\proofname{Proof}

\def\satzname{Satz}
\def\koroname{Korollar}
\def\folgname{Folgerung}
\def\bemerkname{Bemerkung}
\def\aufgname{Aufgabe}

\def\beisname{Beispiel}
\def\beissname{Beispiele}
\def\bewname{Beweis}

\def\@thmcounter#1{\noexpand\arabic{#1}}
\def\@thmcountersep{}
\def\@begintheorem#1#2{\it \trivlist \item[\hskip
\labelsep{\bf #1\ #2.\quad}]}
\def\@opargbegintheorem#1#2#3{\it \trivlist
      \item[\hskip \labelsep{\bf #1\ #2.\quad{\rm #3}}]}
\makeatother
\newtheorem{theor}{\theoremname}[section]
\newtheorem{propo}[theor]{\propositionname}
\newtheorem{coro}[theor]{\corollaryname}
\newtheorem{lemm}[theor]{\lemmaname}

\newenvironment{thm}{\begin{theor}\it}{\end{theor}}

\newenvironment{prop}{\begin{propo}\it}{\end{propo}}

\newenvironment{cor}{\begin{coro}\it}{\end{coro}}

\newenvironment{lem}{\begin{lemm}\it}{\end{lemm}}

\newtheorem{rema}[theor]{\remarkname}

\newenvironment{rem}{\begin{rema}\rm}{\end{rema}}

\newtheorem{stepnow}[theor]{}

\newtheorem{defin}[theor]{\definitionname} 

\newenvironment{definition}{\begin{defin}\rm}{\end{defin}}
\newenvironment{defn}{\begin{defin}\rm}{\end{defin}}

\newtheorem{exerc}[theor]{\exercisename}

\newtheorem{exa}[theor]{\examplename}

\newenvironment{ex}{\begin{exa}\rm}{\end{exa}}

\newtheorem{exas}[theor]{\examplesname}

\newenvironment{exs}{\begin{exas}\rm}{\end{exas}}

\newtheorem{conj}[theor]{\conjecturename}

\newtheorem{pro}[theor]{\problemname}

\newenvironment{prob}{\begin{pro}\rm}{\end{pro}}

\newtheorem{prs}[theor]{\problemsname}

\newenvironment{Proof*}{\begin{trivlist}\item[\hskip%
\labelsep{\bf\proofname.\quad}]}%
{\end{trivlist}}

\newenvironment{prf}{\begin{proof}}{\end{proof}}

{\hfill\qed\end{trivlist}}

\newenvironment{beweis*}{\begin{trivlist}\item[\hskip%
\labelsep{\bf\bewname.\quad}]}%
{\end{trivlist}}

\newtheorem{satzn}[theor]{\satzname}

\newtheorem{koro}[theor]{\koroname}

\newtheorem{folg}[theor]{\folgname}

\newtheorem{bem}[theor]{\bemerkname}

\newtheorem{aufg}[theor]{\aufgname}

\newtheorem{aufgn}[theor]{\aufgname}

\newtheorem{beis}[theor]{\beisname}

\newtheorem{beiss}[theor]{\beissname}

\renewcommand\mlabel{\label}

\DeclareMathAlphabet{\Ma}{U}{msa}{m}{n}
\DeclareMathAlphabet{\Mb}{U}{msb}{m}{n}
\DeclareMathAlphabet{\Meuf}{U}{euf}{m}{n}

\def\got#1{\Meuf{#1}}
\def\br#1.{\llbracket #1 \rrbracket}
\def\lbr{\llbracket}
\def\rbr{\rrbracket}

\renewcommand{\baselinestretch}{1.2}
\normalsize
 \def\ot #1.{{\got{#1}}}

\newcommand{\Ind}{\mathop{{\rm Ind}}\nolimits}

\def\fB{\mathfrak{B}}
\def\al#1.{{\cal #1}}
\def\wt{\widetilde}
\def\slim{\mathop{\hbox{\rm s-lim}}}
\def\s #1.{_{\smash{\lower2pt\hbox{\mathsurround=0pt $\scriptstyle #1$}}\mathsurround=5pt}}
\def\XP#1!{\renewcommand{\baselinestretch}{.7}\marginpar{{\footnotesize
$\leftarrow$#1}\hfil}\renewcommand{\baselinestretch}{1.2}}
\def\ccr #1,#2.{\overline{\Delta(#1,\,#2)}}
\def\rsl{\mathord{\al R.(X,\sigma)}}
\def\b#1.{{\bf #1}}

\begin{document}


\title{Crossed products of $C^*$-algebras for singular actions}

\author{
 {\sc Hendrik Grundling}                                            \\[1mm]
 {\footnotesize Department of Mathematics,}                         \\
 {\footnotesize University of New South Wales,}                      \\
 {\footnotesize Sydney, NSW 2052, Australia.}                             \\
 {\footnotesize hendrik@maths.unsw.edu.au}                  \\
 {\footnotesize FAX: +61-2-93857123}   \\
\and
  {\sc Karl--Hermann Neeb}              \\[1mm]
 {\footnotesize Department of Mathematics,}                           \\
 {\footnotesize
Friedrich-Alexander-Universit\"at Erlangen-N\"urnberg,}  \\
{\footnotesize  Cauerstr. 11,}
{\footnotesize  91058 Erlangen, Germany}                        \\
 {\footnotesize  karl-hermann.neeb@math.uni-erlangen.de}}

\maketitle

\begin{abstract}
\noindent
We consider group actions $\alpha \: G \to \Aut(\cA)$ of topological groups $G$
on $C^*$-algebras $\cA$ of the type which occur in many physics models. These are singular actions
in the sense that they need not be strongly continuous, or the group need not be locally compact.
We develop a ``crossed product host'' in analogy to the usual crossed product for strongly continuous
actions of locally compact groups, in the sense that its representation theory is in
a natural bijection
with the covariant representation theory of the action $\alpha \: G \to \Aut(\cA)$.
We prove a uniqueness theorem for crossed product hosts, and analyze existence conditions.
We also present a number of examples where a crossed product host exists, but the usual crossed product does not.
For actions where a crossed product host does not exist, we obtain a ``maximal" invariant subalgebra
for which a crossed product host exists. We further
study the case of a discontinuous action
$\alpha \: G \to \Aut(\cA)$  of a locally compact group in detail.
\end{abstract}

\section{Introduction}

A large part
of the  $C^*$-algebra literature
is concerned with
crossed products for strongly continuous group actions $\alpha \: G \to \Aut(\cA)$ of locally compact groups
on $C^*$-algebras. A particularly useful property of a crossed product is that
there is a natural bijection of its representation theory with the
covariant representation theory of the triple $(\cA, G, \alpha)$.
This allows
one to analyze the covariant representation theory of an action $\alpha$
by analyzing the representation theory of the crossed product,
using the usual $C^*$-algebraic tools and structure theory.
In addition, crossed products  allow one to  extend tools of group representation
theory to covariant representations, e.g.\ Mackey's theory of induced representations.
According to \cite{Ta67} this was one of the main motivations to
introduce crossed products. We refer to \cite{Wil07} for references to the basic literature.

The paper of Doplicher, Kastler and Robinson~\cite{DKR66}, where crossed products were first defined
in full generality,  was however motivated by physics. The main idea here was that the field $C^*$-algebra of a
(relativistic) quantum system should encode in bounded form all the algebraic relations of the concrete system
of operators describing the system, but in addition, its
representation theory should comprise only the physically relevant representations. In the quantum field theory context
this meant that the field algebra should have as representations all those covariant representations w.r.t.\ the action of the Poincar\'e group,
which satisfy a certain spectral condition for the translation subgroup. In~\cite{DKR66}, the crossed product for a strongly continuous
action $\alpha \: G \to \Aut(\cA)$ of a locally compact group $G$ acting on a $C^*$-algebra $\cA$ was constructed
to carry the covariant representation theory. In a later section of the paper, it was shown that one could factor out
by an ideal of the crossed product in order to obtain a $C^*$-algebra whose representations correspond to covariant representations
which satisfy the desired spectral condition.

In contrast to the mathematical success of crossed products, the application of crossed products to physics has
been limited. As pointed out
on many occasions by Borchers~\cite{Bo83, Bo96}, in quantum field theory the natural actions $\alpha \: G \to \Aut(\cA)$
are usually not strongly continuous, and the groups $G$ need not be locally compact, e.g.\ they may be   infinite--dimensional Lie groups,
such as the group of gauge transformations.
Here we want to address this particular problem, by
developing a ``crossed product'' for such singular actions,
 i.e. we want to construct a $C^*$-algebra $\cal C$ whose  representations
correspond to covariant representations of $\alpha \: G \to \Aut(\cA)$ in the same manner than for
the usual crossed product. Such an algebra $\cal C$  will be called a crossed product host, and our
definition for it below borrows ideas
from Raeburn (cf.~\cite{Rae88}), and sidesteps the usual construction via twisted convolution algebras.
We will include spectral conditions into the framework from the
start, because it is sometimes easier to construct the crossed product host whose representations satisfy
a spectral condition, than one who does not. We will obtain existence and uniqueness results
for crossed product hosts,
and provide a number of examples where a crossed product host exists, but the usual crossed product does not exist.
We will study the case of a discontinuous action
$\alpha \: G \to \Aut(\cA)$  of a locally compact group $G$ in detail.

A general crossed product is not particularly useful for physics, because, as a rule, not
all representations of the field algebra are physically acceptable;- there is usually a specified subclass of representations,
 and only covariant representations from this subclass are allowed. This means we are sometimes forced to
 deal with a situation where we know that a representation is not the representation  of a crossed product host
 for the given action. Our strategy will be to reduce $\cA$ to a subalgebra $\cA_\cL$ preserved by the action $\alpha_G$,
for which the given representation becomes that of a crossed product host, and it is maximal in a suitable sense.
The justification is as follows.

The construction of a $C^*$-algebra whose representation theory satisfies
some physical condition, can be thought of as an application of quantum constraint theory. In quantum constraint theory
one starts with some convenient initial $C^*$-algebra which encodes algebraically
the system we want to describe. We then impose physically motivated restrictions
on its representations, and seek the optimum $C^*$-algebra which will have as representations only
those which satisfy the restrictions. These restrictions can be where one requires some operator identity to hold
(``supplementary conditions'', cf.~\cite{Gr06}), or requiring that some subgroup of a symmetry group
acts trivially (cf.~\cite{La95}), or that some one-parameter unitary groups must be strong operator continuous
(regularity condition for the Weyl algebra~\cite{GrN09}), or that some symmetry group is covariantly represented and satisfies some
spectral condition (cf.~\cite{Bo96}). Typically there will be some elements of the original algebra which
are incompatible with the constraint conditions, and these need to be excluded before one constructs the
final physical algebra, and we will find a similar situation below.

The structure of this paper is as follows. We start with a short list of essential notation and concepts in Section~2, and
in the next section we recall the concept of a host algebra, which generalizes the concept of a group algebra
beyond locally compact groups, and also for subtheories of the continuous representation theory rather than
the full continuous representation theory. In Section~\ref{CPH-unique} we introduce our definition of a
crossed product host relative to a choice of host algebra $\cL$,
 explore its basic properties, give examples, and obtain a uniqueness theorem.
Then in In Section~\ref{CPH-exist} we analyze existence criteria for crossed product hosts, and these
are combined in the concept of a cross-representation, which always allow for the construction of a
crossed product host. We study the properties of these representations, and construct a
universal covariant $\cL$-representation. This is a cross representation if and only if there is a
crossed product host for the full theory of covariant
$\cL$-representations. We give examples where a crossed product host exists, but the usual crossed product does not exist.
In Section~\ref{XrepSC} we consider special situations where we can prove existence of crossed product hosts,
i.e. that a given covariant representation ${(\pi,U)}$ of $\alpha \: G \to \Aut(\cA)$ is a
cross representation.
These include where the host $\cL$ is represented by compact operators, where $G$ is compact
and $\Spec(U)$ is finite or $U$ is of finite multiplicity, and finally, if we have a semidirect product
$G=N \rtimes_\varphi  H$ of locally compact groups, and ${(\pi,U\!\restriction\! N)}$ and ${(\pi,U\!\restriction\! H)}$
 are cross representations. For each special situation we develop a number of further interesting examples.
In Section~\ref{NonX} we consider  non-cross representations, i.e. covariant representations which do not allow
construction of a crossed product host. The motivation comes from physics, as explained above.
We isolate a natural subalgebra $\cA_\cL$ preserved by the action $\alpha_G$,
on which the given representation restricts to a cross representation.
In Section~\ref{CPH-discont} we study the special case of a discontinuous action
$\alpha \: G \to \Aut(\cA)$  of a locally compact group $G$, and give very concrete criteria for
the existence of a crossed product host.
Finally, in Section~\ref{SectExmp} we develop three examples which are useful for
bosonic quantum systems, and where crossed products do not exist in the usual sense.
We conclude with a number of appendices, in which we give technical lemmas required
in the proofs as well as further related material which do not easily fit into the main text.

\tableofcontents

\section{Basic concepts and notation}

Below we will need the following.
\begin{defn}
\mlabel{def:1.1c}
\begin{itemize}
\item[(i)]
For a $C^*$-algebra $\cA$, we write $M({\cal A})$ for the
multiplier algebra of ${\cal A}$. If
${\cal A}$ has a unit, $\U({\cal A})$ denotes  its unitary group.
There is an injective morphism of $C^*$-algebras
$\iota_{\cal A} \: {\cal A} \to M({\cal A})$ and we will just write
${\cal A}$ for its
image in $M({\cal A})$. Then ${\cal A}$ is dense
in  $M({\cal A})$ with respect to
the {\it strict topology},
which is the locally convex topology defined by the seminorms
$$ p_A(M) := \|M \cdot A\| + \|A \cdot M\|,
\qquad A\in {\cal A},\; M\in M({\cal A})$$
(cf.\ \cite[Prop.~3.5]{Bu68} and \cite[Prop.~2.2]{Wo95}).
\item[(ii)]
Let $\cA$ and $\cL$ be $C^*$-algebras and
$\phi \: \cA \to M(\cL)$ be a morphism of $C^*$-algebras. We
say that $\phi$ is {\it non-degenerate} if ${\rm span}(\phi(\cA)\cL)$
 is dense in $\cL$ (cf.\ \cite{Rae88}).
A representation $\pi:\cA\to{\cal B}({\cal H})$ is called {\it nondegenerate}
if $\pi(\cA){\cal H}$ is dense in the Hilbert space ${\cal H}$.
\end{itemize}
\end{defn}
\noindent Nondegeneracy is an important concept, and in Appendix~A we list
useful properties for it.
If $\phi \: \cA \to M(\cB)$ is a morphism of $C^*$-algebras
which is non-degenerate, then we  write $\tilde\phi \: M(\cA) \to M(\cB)$ for its
uniquely determined extension to the multiplier algebras
(cf.\ \cite[Prop.~10.3]{Ne08}).
\medskip

For a complex Hilbert space ${\cal H}$, we write $\Rep({\cal A},{\cal H})$ for the
set of non-degenerate representations of ${\cal A}$ on ${\cal H}$,
and $\ot S.(\cA)$ for the set of states of $\cA.$
To avoid set--theoretic subtleties, we will express our results below
concretely, i.e., in terms of $\Rep({\cal A},{\cal H})$ for given Hilbert spaces
$\al H..$
We have an injection
$$ \Rep({\cal A}, {\cal H}) \into \Rep(M({\cal A}),{\cal H}), \quad \pi \mapsto \tilde\pi
\quad \hbox{ with } \quad \tilde\pi \circ \iota_{\cal A} = \pi, $$
which identifies a non-degenerate representation $\pi$ of
${\cal A}$ with the representation $\tilde\pi$ of its multiplier algebra
which extends $\pi$.  The representations of $M(\cA)$ on $\cH$ arising
from this extension process are characterized as those representations which are
continuous with respect to the strict topology on $M(\cA)$
and the strong operator topology on $\cB(\cH)$, or equivalently
by non-degeneracy of their restriction to
$\cA$ (cf.\ \cite[Prop.~10.4]{Ne08}).
We will refer to $\tilde\pi$ as the {\it multiplier extension}
of $\pi$.
It can be obtained by
\[
\wt{\pi}(M)=\slim_{\lambda\to\infty}\pi(M E_\lambda)\quad\forall\,M\in M(\al A.)
\]
where $(E_\lambda)_{\lambda\in\Lambda}$ is any approximate
identity  of $\al A..$
\medskip

If $\cB$ is a $C^*$-algebra, and  $X$ is a left Banach $\cB\hbox{--module,}$
then the closed span of $\cB X$ satisfies $\overline{\rm span}(\cB X)=\cB X=
 \{ Bx \mid B \in \cB, x \in X\}$ (cf.~\cite[Th.~II.5.3.7]{Bla06} or \cite[Th.~5.2.2]{Pa94}).
 In particular it implies that if $\phi \: \cA \to M(\cL)$ is non-degenerate, then
 $\cL=\phi(\cA)\cL$, and if $\pi:\cA\to\cB(\cH)$ is a non-degenerate representation, then
 $\overline{\rm span}(\pi(\cA)\cH)=\pi(\cA)\cH$.
 We will use the notation $\br S.:=\overline{\rm span}(S)$ where $S\subset Y$
 and $Y$ is a Banach space.
 \medskip

For  topological groups $G$ and $H$ we write
$\Hom(G,H)$ for the set of continuous group homomorphisms $G \to H$.
We also write
$\Rep(G,{\cal H})$ for the set of all (strong operator)
continuous unitary representations of
$G$ on ${\cal H}\,.$ 
Endowing $\U({\cal H})$ with the strong operator topology
turns it into a topological group, denoted $\U({\cal H})_s$,
so that $\Rep(G,{\cal H}) = \Hom(G,\U({\cal H})_s)$.

\begin{defn}\mlabel{def:1.1}
(i) We write  $(\cA, G, \alpha)$ for a triple, where
$\cA$ is a $C^*$-algebra, $G$ a topological group and
$\alpha \: G \to \Aut(\cA)$ is a homomorphism.
We call $\alpha$ {\it strongly continuous} if for every $A \in \cA$, the
orbit map $\alpha^A \: G \to \cA,$ $g \mapsto \alpha_g(A)$ is continuous.
If $\alpha$ is strongly continuous,
we call $(\cA, G, \alpha)$ a {\it $C^*$-dynamical system} (cf.~\cite{Pe89},
\cite[Def.~2.7.1]{BR02}), or say that the action is $C_0$.
Unless otherwise stated, we will not assume that $\alpha$ has this property
and simply speak of the triple $(\cA, G, \alpha)$ as an {\it
$C^*$-action}. The {\it usual case} will mean that the action is $C_0$, and the group
$G$ is locally compact.

(ii) A {\it covariant representation of $(\cA, G, \alpha)$}
is a pair $(\pi,U)$, where
$\pi \: \cA \to \cB(\cH)$ is a nondegenerate representation of $\cA$
on the Hilbert space $\cH$ and $U \: G \to \U(\cH)$ is a
continuous unitary representation satisfying
\begin{equation}
  \label{eq:covar}
U(g)\pi(A)U(g)^* = \pi(\alpha_g(A)) \quad
\mbox{ for } \quad g \in G, a \in \cA.
\end{equation}
For a fixed Hilbert space $\cH$, we write
${\rm Rep}(\alpha,\cH)$ for the set of covariant representations
$(\pi,U)$ of $(\cA, G, \alpha)$ on $\cH$.

\end{defn}

\section{Generalizing group algebras -- host algebras}

In the usual case for
$\alpha \: G \to \Aut(\cA)$, the
crossed product $\cA \rtimes_\alpha G$  is constructed from
the group algebra $C^*(G)$ and $\cA$ (cf.\ \cite{Rae88}).
Here we want to follow a similar strategy.

For our case, as we do not necessarily make the usual assumptions for  $(\cA, G, \alpha)$,
we need to generalize
the concept of a group algebra (cf.~\cite{Gr05, Ne08}). This is to allow for
topological groups which are not locally compact, or for representation theories restricted by
some given condition (e.g.\ spectral conditions - see \cite{GrN12}).
The main task of the usual group algebra
is to carry the continuous unitary representation
theory of the group. This is the property which we will generalize,
by seeking a $C^*$-algebra which will carry the desired subclass of
representations of $G$ in a natural way.
\begin{definition} \mlabel{def:2.1a} Let $G$ be a topological group.
A {\it host algebra for $G$} is a pair
$({\cal L}, \eta)$, where  ${\cal L}$ is a $C^*$-algebra and
$\eta \: G \to \U(M({\cal L}))$ is a group homomorphism
such that:
\begin{description}
\item[\rm(H1)] For each non-degenerate representation $(\pi, {\cal H})$
of $\cL$, the representation $\tilde\pi \circ \eta$ of $G$ is continuous.
\item[\rm(H2)] For each complex Hilbert space
${\cal H}$, the corresponding map
$$ \eta^* \: \Rep({\cal L},{\cal H}) \to \Rep(G, {\cal H}), \quad
\pi \mapsto \tilde\pi \circ \eta $$
is injective.
\end{description}
We write $\Rep(G,{\cal H})_\eta$ for the range of $\eta^*$,
and its elements are called {\it $\cL$-representations of $G$}.
Note that $\eta^*$ depends on $\cal H$.\\[3mm]
We call $({\cal L}, \eta)$  a {\it full host algebra} if, in addition, we have:
\begin{description}
\item[\rm(H3)] $\Rep(G,{\cal H})_\eta = \Rep(G,{\cal H})$
for each Hilbert space~${\cal H}$.
\end{description}

Given $U \in \Rep(G,\cH)_\eta$, we write
$U_\cL :=(\eta^*)^{-1}(U)\in \Rep(\cL,\cH)$ for the unique representation
of $\cL$ such that $\widetilde{U_\cL} \circ \eta = U$.
\end{definition}

Thus by (H2) and (H3), a full host algebra, when it exists, carries precisely the continuous
unitary representation theory of $G$, and if it is not full, it carries some subtheory of the
continuous unitary representations of $G$. In particular, if we want to impose additional restrictions,
e.g.\ a spectral condition, then we will specify a host algebra which is not full.
Since specification of a non-full host $\cL$ restricts the representation theory of
$G$, it can be thought of as a quantum constraint. This will make more sense below in
the context where $G$ acts on $\cA$.
In general, host algebras need not exist, as there are topological groups with continuous
unitary representations, but without
irreducible ones (cf.~\cite{GN01}), and $\eta^*$ preserves irreducibility
for host algebras. The existence of a host algebra for the fixed subclass of representations
of $G$ means that this class of representations is ``isomorphic'' to the representation theory
of a $C^*$-algebra. The standard example is when $G$ is locally compact, and we take
$\cL=C^*(G)$ with the canonical map
$\eta \: G \to \U(M(C^*(G)))$. At the end of this section there is a list
of other examples.

We consider the strict continuity of the map  ${\eta \: G\to \U(M(\cL))}$
in the definition (cf.~\cite[Rem.~4.14]{Ne08}):

\begin{prop} \mlabel{prop:3.3}
Let $G$ be a topological group, ${\cal L}$ be a $C^*$-algebra
and $\eta \: G \to \U(M({\cal L}))$ be a group homomorphism.
\begin{description}
\item[\rm(i)] If $\eta \: G\to \U(M(\cL))$ is  strictly continuous,  then
$\tilde\pi \circ \eta$  is continuous
for each $\pi\in \Rep({\cal L},{\cal H})$, i.e. (H1) is satisfied.
\item[\rm(ii)] $\eta(G)$ spans a strictly dense subalgebra of $M(\cL)$
if and only if  $\eta^*$ is injective for all $\cH$.
\item[\rm(iii)] Let $\cL_c^L:=\{L \in \cL\,\mid\,
G \to \cL,  g\mapsto \eta(g)L\quad\hbox{is continuous}\}$
 and $\cL_c := \cL_c^L \cap (\cL_c^L)^*$. Then $\cL_c^L$
 is a closed right ideal, biinvariant under $G$,
 and $\cL_c$ is a $C^*$-subalgebra of $\cL$ which is
$G$-biinvariant, and the corresponding homomorphism
$\eta_c \: G \to \U(M(\cL_c))$ is strictly continuous.
\item[\rm(iv)] Let $\eta \: G \to \U(M(\cL))$ be such that
$\eta^*$ is injective for all $\cH$.
Then $\cL_c$ is a closed two-sided ideal of $M(\cL)$ and the corresponding
morphism $\gamma \: M(\cL) \to M(\cL_c)$ is
strictly continuous and satisfies $\gamma \circ \eta = \eta_c$.
Moreover $(\cL_c, \eta_c)$ is a host algebra for $G$.
\end{description}
\end{prop}

\begin{prf} (i) This follows from the
 strict continuity of the multiplier extension $\tilde\pi$.

(ii) If $\eta(G)$ spans a strictly dense subalgebra of $M(\cL)$, then
as  the extension $\tilde\pi$ of a non-degenerate representation
$\pi$ of $\cL$ is strictly continuous, $\eta^*$ is injective for all $\cH$.
If, conversely, $\eta^*$ is injective for all $\cH$, then
$\Spann(\eta(G))$ is strictly dense in $M(\cL)$ by
\cite[Prop.~2.2]{Wo95}.

(iii) Most of the claims are direct verifications.  Note that $\cL_c^L$
 is closed because $G$ acts by isometries on $\cL$.

 (iv) Since $\eta^*$ is injective for all $\cH$, (ii) implies
that $\eta(G)$ spans a strictly dense subalgebra of $M(\cL)$.
Then the $G$-biinvariant closed subalgebra
$\cL_c$ of $\cL$ is a two-sided ideal of $M(\cL)$, so $\gamma$
has the claimed properties. To verify the host properties of
$\cL_c$, note that (H1) follows from the strict
continuity of $\eta_c$. Moreover $\gamma(M(\cL))$ contains $\cL_c$, so that it is strictly
dense in $M(\cL_c)$ (\cite[Prop.~2.2]{Wo95}). Therefore
$\eta_c(G) =  \gamma(\eta(G))$ spans a strictly dense subalgebra of
$M(\cL_c)$, which implies that $\eta^*$ is injective for all $\cH$, as we have seen in (ii).
\end{prf}
\begin{rem} \mlabel{vNalgHost}
If $({\cal L}, \eta)$  a   host algebra for $G$ and
 $(\pi, {\cal H})$ is a non-degenerate representation of $\cL$, then
 $\pi(\cL)''=(\tilde\pi \circ \eta)(G)''$. To see this,
note that by Proposition~\ref{prop:3.3}(ii),  $\Spann\eta(G)\subset M(\cL)$ as well as
$\cL\subset M(\cL)$ are both strictly dense algebras. Thus using the strictly continuous
extension of $\pi$ from $\cL$ to $M(\cL)$, we get that the strong operator closures of
$\pi(\cL)$ and  $\Spann (\tilde \pi \circ \eta)(G)$ are the same. It follows for the  commutants that
$\pi(\cL)'=(\tilde \pi \circ \eta)(G)'$, hence that
$\pi(\cL)''=(\tilde \pi \circ \eta)(G)''$.
\end{rem}

\begin{defn} \mlabel{def:morhost}
Let $G$ and $H$ be topological groups and
$(\cL_G,\eta^G)$, resp. $(\cL_H,\eta^H)$, be host algebras for $G$,
resp.\ $H$. Then a {\it morphism of host algebras} is a pair
$(\phi, \psi)$, where
$\phi \:  H \to G$ is a continuous group homomorphism and
$\psi \: \cL_H \to \cL_G$ is a non-degenerate
morphism of $C^*$-algebras whose extension $\tilde \psi \: M(\cL_H) \to
M(\cL_G)$ satisfies
\[ \tilde\psi \circ \eta^H = \eta^G \circ \phi.\]
\end{defn}

\begin{ex} \mlabel{ex:3.5} 
Typical examples of morphisms of host algebras
in the case of locally compact groups with group algebras as hosts, arise from
inclusions  $\phi \: H \into G$ of open subgroups.
Inclusions of closed subgroups
$H \into G$ in general do not not induce inclusions
$C^*(H) \to C^*(G)$. For example
the inclusion $\phi \: H = \T \into G = \T^2, t \mapsto (t,1)$
cannot produce an inclusion $C^*(H) \to C^*(G)$
because the regular representation of $G$ on $L^2(G)$ represents
$C^*(G)$ by compact operators but the operators coming from
$C^*(H)$ have infinite-dimensional eigenspaces, hence
are not compact if they are non-zero.

\end{ex}

\begin{exs} \mlabel{Exmp3.2}
Presently we know the following sources
of host algebras:
\begin{itemize}
\item[{\rm (1)}] If $G$ is locally compact,
we take
$\cL=C^*(G)$ with the canonical map
$\eta \: G \to \U(M(C^*(G)))$ which is strictly continuous,
i.e. continuous with respect to the strict topology on \break $\U(M(C^*(G)))$
and defines on $C^*(G)$ the structure of a full host algebra for
the class of continuous representations of $G$ (cf.\ \cite[Sect.~13.9]{Dix77}).
\item[{\rm (2)}] If $\cL$ is a host algebra of $G$ and
$\cI \subeq \cL$ is a closed ideal, then
$\cL/\cI$ also is a $C^*$-algebra, the quotient map $q \: \cL \to \cL/\cI$
induces a surjective homomorphism
$\tilde q \: M(\cA) \to M(\cL/\cI)$, and
$\tilde q \circ \eta_G \: G \to M(\cL/\cI)$ defines on
$\cL/\cI$ the structure of a host algebra for $G$.
By construction we then have a morphism from the host algebra
$\cL$ (w.r.t.~$G$) to the host algebra $\cL/\cI$ (w.r.t.~$G$).
\item[{\rm (3)}] Let $(V,\omega)$ be a countably dimensional symplectic space
and $\Heis(V,\omega) := \R \times V$ the corresponding Heisenberg group
with the multiplication
\[ (t,v)(t',v') := (t + t' + \shalf\omega(v,v'), v  + v')\]
endowed with the direct limit topology turning it into a
topological group.
We have shown in \cite{GrN09} that there exists a host algebra
$(\cL,\eta)$ of $\Heis(V,\omega)$ for which
$\Heis(V,\omega)_\eta$ is the class of all continuous unitary
representations $U$ satisfying
$U(t,0) = e^{it}\1$ for $t \in \R$.
If $\dim V < \infty$, then $\Heis(V,\omega)$ is locally compact, and so
the existence of such a host algebra
 already follows from Remark~\ref{FactorGpalg} below,
applied to $Z := \R \times \{0\}$
and $\chi(t,0) := e^{it}$.

\item[{\rm (4)}] Host algebras obtained from complex semigroups (cf.\ \cite[Sect.~III.2]{Ne00},
\cite{Ne08}, \cite{MN11}).
We shall discuss this class in more detail in the forthcoming paper (\cite{GrN12}). In this
case $G$ is a Lie group, possibly infinite-dimensional,
and we have a  complex involutive semigroup $S$ on which
$G$ acts by unitary multipliers via a
map $\eta \: G \to \U(M(S))$ satisfying certain
additional conditions.
The corresponding host algebra is
$\cL = C^*(S,\beta)$, the enveloping $C^*$-algebra of $S$
corresponding to a locally bounded absolute value $\beta$ on $S$ which is invariant under $G$.
This is closely related to spectral conditions.

\item[{\rm (5)}]
Host algebras  for subsets of norm continuous representations in
$\Rep(G,{\cal H})$ can  be obtained
for general topological groups.
If $(U,\cH)$ is a norm continuous unitary
representation of the topological group $G$ and
$\cL := C^*(U(G))$ is the $C^*$-subalgebra of $\cB(\cH)$ generated by $U(G)$, then
$\eta := U \: G \to \U(\cL)$ is continuous, so that every continuous representation
$\pi$ of $\cL$ leads to a norm continuous representation $\pi \circ \eta$ of $G$.
Hence $(\cL,\eta)$ is a host algebra for $G$  (see also  \cite[Prop.~III.2.18]{Ne00}).

\item[{\rm (6)}]
If $G$ is an abelian topological group,
there is a different method to obtain
hosts algebras  for its norm continuous representations.
Let
$\hat G := \Hom(G,\T)$ be its group of continuous characters, and
$\Sigma\subeq \hat G$ be an equicontinuous subset which is closed in the topology
of pointwise convergence, i.e., as a subset of the compact group
$\T^G$. Then the topology of pointwise
convergence turns $\Sigma$ into a
compact Hausdorff space, so that
$\cL := C(\Sigma)$ is a unital $C^*$-algebra.
We obtain a continuous homomorphism
\[ \eta \:  G \to \U(\cL), \quad \eta(g)(\chi) := \chi(g), \]
and it is easy to see that $(\cL,\eta)$ is a host algebra for $G$.
The $\cL$-representations of $G$ are precisely those representations given
by Borel spectral measures $P$ on the compact space $\Sigma$ via
$U_g = \int_{\Sigma} \chi(g)\, dP(\chi)$. All these representations
are norm continuous and satisfy $U = \pi \circ \eta$, where $\pi(f) = \int_\Sigma f(\chi)\, dP(\chi)$.
According to \cite[Prop.~7.6]{HoMo98}, typical examples are obtained
for locally compact abelian groups $G$ by compact subsets
$\Sigma$ of the dual group $\hat G$, endowed with the compact open
topology.

\end{itemize}
\end{exs}

\begin{rem} \mlabel{FactorGpalg}
By combining (1) with (2) in Examples~\ref{Exmp3.2},
note that when $G$ is locally compact, all quotients $C^*(G)/\cI$ are also host algebras
for $G$. Such quotients are necessary when we impose spectral conditions,
and we will develop this subject in the sequel \cite{GrN12}.
A particularly common occurrence of this  are quotients
defined by a character $\chi \: Z \to \T$ of a closed central
subgroup $Z \subeq G$. To define these quotients,
let $q \: G \to G/Z$ denote the quotient map
and observe that
\[ C_c(G;\chi) := \{ f \in C_c(G) \: (\forall g \in G,z  \in Z)\,
f(gz) = \chi(z)^{-1} f(g)\} \]
is closed under convolution and the involution $*$, so that we obtain
by completion with respect to the norm
$\|f\|_{1,\chi} := \int_{G/Z} |f(g)|\, d(gZ)$ a
Banach-$*$-algebra $L^1(G;\chi)$ whose enveloping
$C^*$-algebra is denoted $C^*(G;\chi)$. The map
\[ p_\chi \: C_c(G) \to C_c(G;\chi), \quad
p_\chi(f)(g) := \int_Z f(gz) \chi(z)\, dz \]
extends to a surjective morphism $p_\chi \: L^1(G) \to L^1(G;\chi)$
and hence to a surjective morphism $p_\chi \: C^*(G) \to C^*(G;\chi)$.
Since $C^*(G;\chi)$ is a quotient $C^*(G)\big/\ker(p_\chi)$, it
is a host algebra of $G$ and the corresponding
unitary representations $U$ of $G$ are those for which
$U(z) = \chi(z)\1$ for $z \in Z$
(see \cite[Lemma~A.VII.4]{Ne00} for details).
\end{rem}

\begin{rem} \mlabel{rem:5.11}
Given a continuous unitary representation $(U,\cH)$  of $G$,
Example~\ref{Exmp3.2}(5) raises the question  when
$\cL_U := C^*(U(G))\subset\cB(\cH)$ is a host algebra.
If $G$ is locally compact and abelian, then we now show that
 $\cL_U = C^*(U(G))$ is a host algebra if and only if $U$ is norm continuous.

Note first that by Lemma~\ref{lem:c.2} in the appendix, $U$ is norm continuous if and only if
$\Sigma := \Spec(U)  \subeq \hat G$ is compact. If this is the case, then the homomorphism
$U \: G \to \cL_U$ turns $\cL_U$ into a host algebra of $G$
as in Example~\ref{Exmp3.2}(5).

Conversely, assume that $U$ is not norm continuous, i.e., $\Sigma$ is not compact.
Considering $\Sigma$ as a Borel subset of
$\hat{G_d}$, the compact group of {\sl all} characters of
$G$ (i.e. the character group of the  discrete group $G_d$),
the spectrum of the representation $U$ of the discrete group $G_d$
is the closure $\Sigma_c$ of $\Sigma$ in $\hat{G_d}$. This implies that
$\cL_U \cong C(\Sigma_c)$, so that $\Sigma_c$ is the set of characters of $\cL_U$.
In particular, this $C^*$-algebra is not a host algebra for
$G$ if  $\Sigma_c \not=\Sigma$. That this is always
the case if $\Sigma$ is non-compact follows from
Glicksberg's Theorem \cite[Thm.~1.2]{Gl62}.
Thus $\cL_U = C^*(U(G))$ is a host algebra if and only if $U$ is norm continuous.

The preceding argument also shows that $\cL_U$ is isomorphic
to the group algebra $C^*(G_d) \cong C(\hat{G_d})$ if and only if
$\Sigma$ is dense in $\hat{G_d} = \Hom(G_d,\T)$.
\end{rem}

\section{Crossed product hosts -- uniqueness}
\label{CPH-unique}

One of the main features of  a crossed product  $\cA \rtimes_\alpha G$
in the usual case
is that its representations are in natural bijection with
the covariant representations of $(\cA, G, \alpha)$.
This is the central property which we want to generalize.

Inspired by Raeburn's approach to crossed
products (cf.~\cite{Rae88} and Example~\ref{ex:1.1} below),
we now define:

\begin{defn} \mlabel{def:a.5} Let $G$ be a topological group,
and let $(\cL,\eta)$ be a host algebra for $G$
and $(\cA,G, \alpha)$ be a $C^*$-action.
We call a triple $(\cC, \eta_\cA, \eta_\cL)$
a {\it crossed product host} for $(\alpha,\cL)$ if
\begin{description}
\item[\rm(CP1)] $\eta_\cA \: \cA \to M(\cC)$ and
$\eta_\cL \: \cL \to M(\cC)$ are morphisms of $C^*$-algebras.
\item[\rm(CP2)] $\eta_\cL$ is non-degenerate (cf. Definition~\ref{def:1.1c}).
\item[\rm(CP3)] The multiplier extension $\tilde\eta_\cL \: M(\cL)
\to M(\cC)$ satisfies in $M(\cC)$ the relations
\[ \tilde\eta_\cL(\eta(g))\eta_\cA(A)\tilde\eta_\cL(\eta(g))^*
= \eta_\cA(\alpha_g (A))\quad \mbox{for all} \quad A \in \cA,\;\hbox{and}\; g \in G.\]
\item[\rm(CP4)] $\eta_\cA(\cA) \eta_\cL(\cL) \subeq \cC$
and $\cC$ is generated by this set as a $C^*$-algebra.
\end{description}
A {\it full } crossed product host for $(\alpha,\cL)$ satisfies in addition:
\begin{description}
\item[\rm(CP5)] For every covariant representation
$(\pi, U)$ of $(\cA,\alpha)$ on $\cH$ for which $U$ is an $\cL$-representation
of $G$, there exists a unique
representation $\rho \: \cC \to \cB(\cH)$ with
\[ \rho(\eta_\cA(A)\eta_\cL(L)) = \pi(A) U_\cL(L)
\quad \mbox{ for } \quad A \in \cA, L \in \cL.\]
\end{description}
Two crossed product hosts $(\cC^{(i)}, \eta_\cA^{(i)}, \eta_\cL^{(i)})$, $i=1,2$,
are {\it isomorphic} if there is an isomorphism $\Phi:\cC^{(1)}\to\cC^{(2)}$ such that
  ${\big(\Phi(\cC^{(1)}),\tilde\Phi\circ\eta_\cA^{(1)},\tilde\Phi\circ\eta_\cL^{(1)}\big)}
=(\cC^{(2)}, \eta_\cA^{(2)}, \eta^{(2)}_\cL)$.
\end{defn}

\begin{rem} \mlabel{rem:2.5} Suppose that (CP1-5) are satisfied, then

(a) The action
$\eta_G:=\tilde\eta_\cL\circ\eta:G\to \U(M(\cC))$ is unitary since
$\eta(G) \subeq \U(M({\cal L}))$ and $\tilde\eta_\cL \: M(\cL)
\to M(\cC)$ is a unital homomorphism.
If $\eta$ is strictly continuous, it follows from Theorem~\ref{thm:a.2}(iii)
that $\eta_G$ is strictly continuous.

(b) That $\cC$ is generated by $\eta_\cA(\cA)\eta_\cL(\cL)$,
i.e., the second part of (CP4), is a consequence of (CP5)
when $\cC$ is full.
In fact, (CP5) implies that the $C^*$-subalgebra $\cC_0$ generated
by $\eta_\cA(\cA) \eta_\cL(\cL)$ has the property that all
representations $\pi \: \cC \to \cB(\cH)$ are uniquely determined
by their restrictions to $\cC_0$. Since every $C^*$-algebra embeds into
some $\cB(\cH)$, this means that the inclusion
$\cC_0 \to \cC$ is an epimorphism, and since epimorphisms
of $C^*$-algebras are surjective (\cite{HoNe95}), $\cC = \cC_0$ follows.

(c) The representations $\rho$ in (CP5) are automatically
non-degenerate. In fact, if $\rho(\cC)v = \{0\}$, then we also have
$\pi(\cA)U_\cL(\cL)v = \{0\}$, and since $\pi$ is non-degenerate,
we obtain $U_\cL(\cL)v = \{0\}$. As $U_\cL$ is also non-degenerate,
we further derive $v= 0$. Since $\rho$ is non-degenerate, it extends to a representation
$\tilde\rho \: M(\cC) \to \cB(\cH)$ and (CP5) immediately leads to
\begin{equation}
  \label{eq:mulrel}
\tilde\rho \circ \eta_\cA = \pi \quad \mbox{ and }
\tilde\rho \circ \eta_\cL = U_\cL.
\end{equation}

(d) If $\cA$ is unital, then $\eta_\cA(\1) = \1$ because  by
Proposition~\ref{prop:a.6}(ii) below, $\eta_\cA \: \cA \to M(\cL)$ is non-degenerate.
Thus $\cL\subset\cC$. If $\cL$ is unital, then $\eta_\cL(\1) = \1$ by (CP2), hence
 $\cA\subset\cC$. Though the usual group $C^*$-algebra $C^*(G)$
is unital if $G$ is discrete,
it is not true that all full host algebras of discrete groups are unital
(cf.~\cite[Example~4.16]{Ne08}).

(e) Definition~\ref{def:a.5}
generalizes crossed products in four directions:
\begin{itemize}
\item The group $G$ need not be locally compact,
\item the action $\alpha$ need not be strongly continuous,
\item the host algebra $\cL$, specifying the corresponding class of
$G$-representations does not have to coincide with $C^*(G)$
when $G$ is locally compact,
\item 
for a non-full crossed product host, we restrict to a subtheory of the
covariant $\cL\hbox{--representations}$ (see below).
\end{itemize}

For general $C^*$-actions
$\alpha \: G \to \Aut(\cA)$, a full crossed product host need not
exist (even if a host algebra $\cL$ exists, cf.\ Example~\ref{nonfullMax}),
which is why we also consider
non-full crossed product hosts.
We will analyze conditions for existence and provide more examples below.

\end{rem}
We first verify that Definition~\ref{def:a.5} coincides with the usual crossed product
in the case of $C^*$-dynamical systems of locally compact groups.

\begin{ex} \mlabel{ex:1.1}
Consider the usual case, i.e. we have a strongly continuous homomorphism
 $\alpha \: G \to \Aut(\cA)$ where
 $G$ is locally compact.  Then  $\cL=C^*(G)$ is a full host algebra for $G$
(Example~\ref{Exmp3.2}(1)).
 We claim that the crossed product algebra $\cA \rtimes_\alpha G$
is a full crossed product host for $(\alpha,\cL)$.

We use Raeburn's characterization (\cite[Prop.~3]{Rae88})
of the crossed product of
 $(\cA, G, \alpha)$ as a $C^*$-algebra
$\cC$ together with a homomorphism
$\eta_\cA \: \cA \to M(\cC)$ and a strictly continuous
homomorphism $\eta_G \: G \to \U(M(\cC))$ satisfying:
\begin{description}
\item[\rm(a)] $\eta_\cA(\alpha_g A) = \eta_G(g) \eta_\cA(A) \eta_G(g)^*$
for $g \in G$, $A \in \cA$.
\item[\rm(b)] For every covariant representation $(\pi, U)$ of
$(\cA, G, \alpha)$, there is a non-degenerate representation
$\pi \times U$ of $\cC$ with
$(\pi \times U)\,\tilde{} \circ \eta_\cA = \pi$ and
$(\pi \times U)\,\tilde{} \circ \eta_G = U$.
\item[\rm(c)] $\eta_\cA(\cA) \eta_G(C_c(G))$ is dense in
$\cC$, where $\eta_G(f)$, $f \in C_c(G)$, refers to the integrated
 representation $\eta_G \: C_c(G) \to M(\cC)$
of the convolution algebra $C_c(G)$ which is defined
because $\eta_G$ is strictly continuous.
\end{description}

Given this crossed product $\cC$,
we verify our conditions (CP1)-(CP5). Observe that the
strict continuity of $\eta_G \: G \to \U(M(\cC))$ leads by integration
to a morphism
$L^1(G) \to M(\cC)$ of Banach $*$-algebras and therefore extends
to a morphism $\eta_\cL \: C^*(G) \to M(\cC)$ of $C^*$-algebras
whose non-degeneracy follows immediately from (c)
(cf.\ Theorem~\ref{thm:a.2}), and now it is clear that
$\eta_G = \tilde\eta_\cL \circ \eta$. This proves (CP1)-(CP3).
The remaining two conditions (CP4), (CP5) follow from (b) and (c).
This proves that the crossed product $\cC=\cA \rtimes_\alpha G$ is a
full crossed product host.

Conversely, let
$(\cC, \eta_\cA, \eta_\cL)$ be a full crossed product host,
where $\cL= C^*(G)$. Then the non-degeneracy of $\eta_\cL$
implies the strict continuity of the multiplier
action $\eta_G:=\tilde\eta_\cL\circ\eta$ of $G$ on $\cC$. Hence (CP1)-(CP5) imply that
$(\cC,\eta_\cA)$ is a crossed product in the sense of Raeburn,
hence isomorphic to the crossed product $C^*$-algebra
$\cA \rtimes_\alpha G$.
\end{ex}

A property of central importance for a crossed product host, is that
it carries the covariant  $\cL$-representations of $G$:

\begin{defn}
Given $(\cA, G, \alpha)$, where the action
$\alpha \: G \to \Aut(\cA)$ need not be
strongly continuous, assume a
host algebra $(\cL,\eta)$ for $G$. Then
a covariant representation
$(\pi, U)$ of $(\cA, G, \alpha)$ on $\cH$ is called an
{\it $\cL$-representation} if $U$ is an $\cL$-representation
(cf.\ Definition~\ref{def:2.1a}). We write
$\Rep_\cL(\alpha, \cH)$ for the set of covariant
$\cL$-representations of $(\cA, G, \alpha)$ on $\cH$.
\end{defn}

\begin{thm}
\mlabel{Thm-HostCross}
Let $(\cC, \eta_\cA, \eta_\cL)$ be a crossed product host
for $(\alpha,\cL)$, and recall the homomorphism
$\eta_G:=\tilde\eta_\cL\circ\eta:G\to \U(M(\cC))$.
Then for each Hilbert space $\cH$ the map
$$ \eta^*_\times \: \Rep(\cC,{\cal H}) \to \Rep(\alpha, \cH),  \quad\hbox{given by}\quad
\eta^*_\times(\rho):=\big(\tilde\rho \circ \eta_\cA, \tilde\rho \circ \eta_G\big)$$
is injective, and its range
$\Rep(\alpha,{\cal H})_{\eta_\times}$ consists of $\cL$-representations
of $(\cA, G, \alpha)$.
If $\cC$ is full, then we also have
$\Rep(\alpha,{\cal H})_{\eta_\times} = \Rep_\cL(\alpha,{\cal H})$.
\end{thm}

\begin{prf}
Let $\rho\in\Rep(\cC,{\cal H})$. Then
$\tilde\rho \circ \eta_G =\tilde\rho \circ(\tilde\eta_\cL\circ\eta)$ is a unitary representation of $G$.
Since $\eta_\cA \: \cA \to M(\cC)$ is non-degenerate by
Proposition~\ref{prop:a.6} below, $\tilde\rho \circ \eta_\cA$ is
nondegenerate by Theorem~\ref{thm:a.2}(vii).
Then (CP3) implies that $\eta^*_\times(\rho)={\big(\tilde\rho \circ \eta_\cA, \tilde\rho \circ \eta_G\big)}
\in\Rep(\alpha, \cH)$. Since $\tilde\rho \circ \eta_G ={\tilde\rho \circ(\tilde\eta_\cL\circ\eta)}$,
it is also clear that $\tilde\rho \circ \eta_G$ is an  $\cL$-representation, hence
that $\Rep(\alpha,{\cal H})_{\eta_\times} \subseteq \Rep_\cL(\alpha,{\cal H})$.
If two representations $\rho_i\in\Rep(\cC,{\cal H})$,
$i =1,2$, produce the same covariant representation
${\big(\tilde\rho_i \circ \eta_\cA, \tilde\rho_i \circ \eta_G\big)}\in\Rep_\cL(\alpha,{\cal H})$,
by (CP4), they coincide on $\cC$, i.e. $\eta^*_\times$ is injective.
In the case that $\cC$ is full, then by
(CP5) we  obtain that $\Rep(\alpha,{\cal H})_{\eta_\times} = \Rep_\cL(\alpha,{\cal H})$.
\end{prf}

Thus   a full crossed product host $(\cC, \eta_\cA, \eta_\cL)$ is a host for the covariant
$\cL$-representations of $\alpha$ in the sense of \cite{Gr05}, and
a non--full crossed product host is a host for a subtheory of these.
It allows one to analyze these representations with the usual $C^*$-algebra
tools.
Since $\eta^*_\times$ preserves much of the structure of representations (e.g.\ direct sums,
subrepresentations, irreducibility, cf.\ \cite[Cor.~2.2]{Gr05}),
the existence of a full crossed product host
means that  the theory of
covariant $\cL$-representations of $\alpha$ is isomorphic
to the representation theory of a $C^*$-algebra. In general one does not
expect this to hold (see Example~\ref{nonfullMax}), and so one is also interested in subtheories
of the covariant $\cL$-representations for which a non--full
crossed product host exists.

Further useful properties of crossed product hosts are:
\begin{prop} \mlabel{prop:a.6} If $(\cC, \eta_\cA, \eta_\cL)$ is
a crossed product host
for $(\alpha,\cL)$, then
\begin{description}
\item[\rm(i)] $\Spann\!\big(\eta_\cA(\cA) \eta_\cL(\cL)\big)$ is dense in $\cC$ and
\item[\rm(ii)] $\eta_\cA \: \cA \to M(\cC)$ is also non-degenerate.
\item[\rm(iii)] If $G_d$ denotes the group
$G$ endowed with the discrete topology,
then there exists a  homomorphism $\eta_d:\cA\rtimes_{\alpha} G_d\to
M(\cC)$ given by $\eta_d(A\delta_g)=\eta_\cA(A)\eta_G(g)$ for $A\in\cA$, $g\in G$,
where $\delta_g\in\ell^1(G)$ is the characteristic function of $\{g\}$.
Moreover $\eta_d\big(\cA\rtimes_{\alpha} G_d\big)$ is strictly dense in $M(\cC)$.
\end{description}
\end{prop}

\begin{prf} (i)  As $(\cL,\eta)$ is a host algebra for $G$,
the subgroup $\eta(G)$ spans a strictly dense subspace of $M(\cL)$
(Proposition~\ref{prop:3.3}(ii)).
Since $\tilde\eta_\cL$ is strictly continuous, we have
for $A \in \cA$ and $L \in \cL$ in $M(\cC)$:
\[ \eta_\cL(L)\eta_\cA(A)
\in \oline{\tilde\eta_\cL(\Spann(\eta(G)))\eta_\cA(A)}
\subeq \br\eta_\cA(\cA)\tilde\eta_\cL(M(\cL))., \]
using (CP3) for the last inclusion, where $\br\cdot.$ denotes closed span.
This implies that
\[ \eta_\cL(L)\eta_\cA(A) \eta_\cL(\cL)
\subeq \br{\eta_\cA(\cA)\eta_\cL(\cL)}..\]
From the non-degeneracy of $\eta_\cL$, using Theorem~\ref{thm:a.2}(iv) and
$\eta_\cL(L)\eta_\cA(A)\in\cC$, it now follows that
\[ \eta_\cL(L)\eta_\cA(A)
\in  \oline{\eta_\cL(L)\eta_\cA(\cA)\eta_\cL(\cL)}
\subeq \br\eta_\cA(\cA)\eta_\cL(\cL)..\]
We conclude that the closed subalgebra generated by
$\eta_\cA(\cA)\eta_\cL(\cL)$ is $*$-invariant and that
\[ \eta_\cA(\cA)\eta_\cL(\cL)\eta_\cA(\cA)\eta_\cL(\cL)
\subeq \eta_\cA(\cA) \br\eta_\cA(\cA)\eta_\cL(\cL).
\subeq \br\eta_\cA(\cA)\eta_\cL(\cL).,\]
i.e., that $\Spann\big(\eta_\cA(\cA) \eta_\cL(\cL)\big)$ is dense in~$\cC$.

(ii) is an immediate consequence of (i) and
$\cA = \cA \cdot \cA$.

(iii) That $\eta_d$  defines a homomorphism, is due to the fact that
$\cA\rtimes_{\alpha} G_d$ is generated by the products $A\delta_g$ subject to the covariance condition.
Since $C^*\big(\eta_\cA(\cA)\eta_G(G)\big)$ separates all the representations of $\cC$ by the host property of
$\cL$, it follows that $C^*\big(\eta_\cA(\cA)\eta_G(G)\big)$  is strictly dense in $M(\cC)$
by \cite[Prop.~2.2]{Wo95}.
\end{prf}

\begin{rem} An important consequence of the density of
$\eta_\cA(\cA) \eta_\cL(\cL)$ in $\cC$ is that all states
$\phi \in \fS(\cC)$ are determined by their values on the products
$\eta_\cA(A)\eta_\cL(L)$, $A \in \cA$, $L \in \cL$, resp., by the
continuous bilinear map
\[ \tilde\phi \: \cA \times \cL \to \C, \quad
\tilde\phi(A,L) := \phi\big(\eta_\cA(A)\eta_\cL(L)\big).\]
\end{rem}

The following theorem generalizes Raeburn's characterization
of unital crossed products.

\begin{thm} {\rm(Uniqueness Theorem)}
Let $(\cL,\eta)$ be a host algebra for the topological group $G$,
and $(\cA, G, \alpha)$ be a $C^*$-action.
Let $(\cC, \eta_\cA, \eta_\cL)$ and
$(\cC^\sharp, \eta_\cA^\sharp, \eta_\cL^\sharp)$ be
crossed product hosts for $(\alpha,\cL)$,  such that
$\Rep(\alpha,{\cal H})_{\eta_\times} = \Rep(\alpha,{\cal H})_{\eta^\sharp_\times}$
for any Hilbert space $\cH$.
Then there exists a unique
isomorphism $\phi \: \cC \to \cC^\sharp$ with
$\tilde\phi \circ \eta_\cA = \eta_\cA^\sharp$ and
$\tilde\phi \circ \eta_\cL = \eta_\cL^\sharp$.
In particular, full crossed product hosts for $(\alpha,\cL)$ are isomorphic.
\end{thm}

\begin{prf}
Note first that an isomorphism $\phi \: \cC \to \cC^\sharp$ which satisfies
$\tilde\phi \circ \eta_\cA = \eta_\cA^\sharp$ and
$\tilde\phi \circ \eta_\cL = \eta_\cL^\sharp$ is uniquely determined by (CP4), so we only need to show existence.

Let $\rho^\sharp \: \cC^\sharp \to \cB(\cH)$ be a faithful non-degenerate
representation of $\cC^\sharp$.
Then  by the hypothesis $\Rep(\alpha,{\cal H})_{\eta_\times} = \Rep(\alpha,{\cal H})_{\eta^\sharp_\times}$
we conclude that
 $(\eta^\sharp)^*_\times(\rho)=(\pi_\cA^\sharp,U)\in
\Rep(\alpha,{\cal H})_{\eta_\times}$,  where \break
$U := \tilde\rho^\sharp \circ \tilde\eta_\cL^\sharp \circ \eta \: G \to \U(\cH)$
and $\pi_\cA^\sharp := \tilde\rho^\sharp \circ \eta_\cA^\sharp$.
By Theorem~\ref{Thm-HostCross}, the map $\eta^*_\times$ is injective, hence
there is a unique representation $\pi \: \cC \to \cB(\cH)$ with
\begin{equation}\label{eq:covrel2}
 \pi(\eta_\cA(A)\eta_\cL(L)) = \pi_\cA^\sharp(A) U_\cL(L) \quad
\mbox{ for } \quad A \in \cA, L \in \cL,
\end{equation}
and it satisfies (using (CP4)):
\begin{equation}\label{eq:covrel3}
 \tilde\pi \circ \eta_\cA = \pi_\cA^\sharp \quad \mbox{ and } \quad
 \tilde\pi \circ \eta_\cL = U_\cL.
\end{equation}
On the other hand, the relation
$U = \tilde\rho^\sharp \circ \tilde\eta_\cL^\sharp \circ \eta$ yields
$U_\cL = \tilde\rho^\sharp \circ \eta_\cL^\sharp$.
Thus (\ref{eq:covrel2}) becomes
\[
\pi(\eta_\cA(A)\eta_\cL(L)) = \tilde\rho^\sharp \big( \eta_\cA^\sharp(A)  \eta_\cL^\sharp(L)\big) \quad
\mbox{ for } \quad A \in \cA, L \in \cL,
\]
and so by (CP4) and \eqref{eq:covrel2}, we have
$\pi(\cC) = \rho^\sharp(\cC^\sharp)$. Then
$\phi := (\rho^\sharp)^{-1} \circ \pi \: \cC \to \cC^\sharp$
is a morphism of $C^*$-algebras whose multiplier extension
$\tilde \phi = (\tilde\rho^\sharp)^{-1} \circ \tilde\pi$
satisfies
\[ \tilde\phi \circ \eta_\cA = \eta_\cA^\sharp \quad \mbox{ and } \quad
 \tilde\phi \circ \eta_\cL = \eta_\cL^\sharp.\]
Changing the roles of $\cC$ and $\cC^\sharp$, we also find a  morphism
$\psi \:  \cC^\sharp \to \cC$
with
\[ \tilde\psi \circ \eta_\cA^\sharp = \eta_\cA \quad \mbox{ and } \quad
 \tilde\psi \circ \eta_\cL^\sharp = \eta_\cL.\]
Then $(\psi \circ \phi)\,\tilde{} \circ \eta_\cA = \eta_\cA$ and
$(\psi \circ \phi)\,\tilde{} \circ \eta_\cL = \eta_\cL$ lead to
$\psi \circ \phi = \id_\cC$. We likewise obtain
$\phi \circ \psi = \id_{\cC^\sharp}$.
This completes the proof.
\end{prf}

\section{Crossed product hosts -- existence}
\label{CPH-exist}

In this section we consider existence conditions and concrete constructions of
crossed product hosts. These are based on the concept of a cross representation
for the pair $(\alpha, \cL)$ defined below.

\begin{thm} \mlabel{thm:5.1}
Let $(\cL,\eta)$ be a host algebra
for the topological group $G$
and $(\cA,G, \alpha)$ be a  $C^*$-action.
\begin{description}
\item[\rm(a)] Let
$(\cC, \eta_\cA, \eta_\cL)$ be a  crossed product host for $(\alpha,\cL)$.
Then for the (faithful) universal representation $(\rho_u, \cH_u)$ of $\cC$,
the corresponding covariant $\cL$-representation
$(\pi, U)$ of $(\cA, G, \alpha)$ satisfies
\begin{eqnarray*}
 \rho_u(\eta_\cA(A)\eta_\cL(L)) &=& \pi(A) U_\cL(L)
\quad \mbox{ for } \quad A \in \cA, L \in \cL, \\[1mm]
\eta^*_\times(\rho_u)=\big(\tilde\rho_u \circ \eta_\cA, \tilde\rho_u \circ \eta_G\big)
&=& {(\pi,U)}\quad\hbox{ and  }\quad\rho_u(\al C.)=C^*\big(\pi(\cA) U_\cL(\cL)\big)\,.
\end{eqnarray*}
\item[\rm(b)] Conversely, let ${(\pi,U)}\in\Rep_\cL(\alpha,{\cal H})$ and put $\al C.:=C^*\big(\pi(\cA) U_\cL(\cL)\big)$.
Then $\pi(\cA)\cup U_\cL(\cL)\subset M(\cC)\subset \cB(\cH)$, and we
obtain morphisms ${\eta_\cA \: \cA \to M(\cC)}$ and
$\eta_\cL \: \cL \to M(\cC)$ determined by
$\eta_\cA(A)C:=\pi(A)C$ and $\eta_\cL(L)C:=U_\cL(L)C$ for
$A\in\cA$, $L\in\cL$ and $C\in\cC$.
Then the following are equivalent:
\begin{itemize}
\item[\rm(i)]$ (\cC, \eta_\cA, \eta_\cL)$ is  a  crossed product host.
\item[\rm(ii)] $\pi(\cA) U_{\cL}(\cL) \subeq
{U_{\cL}(\cL)\cB(\cH)}$.
\item[\rm(iii)] For every  approximate identity
$(E_j)_{j\in J}$ of $\cL$ we have
\[ \| U_{\cL}(E_j) \pi(A) U_{\cL}(L)
- \pi(A) U_{\cL}(L) \| \to 0 \quad \mbox{ for } \quad
A \in \cA, L \in \cL.\]
\item[\rm(iv)] There exists an approximate identity
$(E_j)_{j\in J}$ of $\cL$ such that
\[ \| U_{\cL}(E_j) \pi(A) U_{\cL}(L)
- \pi(A) U_{\cL}(L) \| \to 0 \quad \mbox{ for } \quad
A \in \cA, L \in \cL.\]
\end{itemize}
\item[\rm(c)]
Let $(\cC, \eta_\cA, \eta_\cL)$ be a  crossed product host.
Then any factor algebra of $\cC$ is also
a crossed product host. In particular, if $\Phi:\cC\to\cC/{\cal J}$ is a factor map where
${\cal J}$ is a closed two-sided ideal, then  ${\big(\cC/{\cal J},\tilde\Phi\circ\eta_\cA,\tilde\Phi\circ\eta_\cL\big)}$
is a crossed product host.
\end{description}
\end{thm}

\begin{prf}
(a) This is a direct consequence of Theorem~\ref{Thm-HostCross} if we start with the  universal representation
$\rho_u$ of $\cC$, which then produces the covariant pair
$\eta^*_\times(\rho_u):=\big(\tilde\rho_u \circ \eta_\cA, \tilde\rho_u \circ \eta_G\big)\in\Rep_\cL(\alpha,{\cal H}_u)$.
Take ${(\pi,U)}$ to be this pair, then this has the properties claimed.

(b) We first show that $\pi(\cA)\cup U_\cL(\cL)\subset M(\cC)\subset \cB(\cH)$.
Observe that the defining representation $\rho$ of $\cC$ on $\cH$ is non-degenerate.
In fact, since the $\cA$-representation $\pi$ is non-degenerate,
the relation $\pi(\cA)U_{\cL}(\cL)v = \{0\}$ implies
$U_{\cL}(\cL)v = \{0\}$, and since $U_{\cL}$ is also non-degenerate,
$v = 0$ follows.
Thus, since $\cC$ is faithfully and nondegenerately represented on $\cH$,
the idealizer of $\cC$ in $\cB(\cH)$ is isomorphic to its multiplier algebra
(cf.\ \cite[Prop.~VIII.1.20]{FD88}, or \cite[II.7.3.5]{Bla06}).
In other words, the map
\[ \varphi:\{ M \in \cB(\cH) \mid
M \cC \cup \cC M \subeq \cC \}\to M(\cC), \qquad
\varphi(M)C=MC\;\hbox{ for }\;C\in\cC, \]
is an isomorphism of $C^*$-algebras.
 Thus a $*$-invariant subset $\cM \subeq \cB(\cH)$ is contained in $M(\cC)$ if
and only if $\cM \cC \subeq \cC= C^*(\pi(\cA) U_{\cL}(\cL)).$
In fact it suffices to check this for the generating elements, i.e. that
$\cM \cdot\pi(\cA) U_{\cL}(\cL) \cup \cM\cdot U_{\cL}(\cL) \pi(\cA)
\subeq \cC$.

We first show that $\pi(\cA)\subeq M(\cC)$. Since  $\cA\cA=\cA$, we have
\[ \pi(\cA)\cdot \big(\pi(\cA)U_{\cL}(\cL)\big)=\pi(\cA)U_{\cL}(\cL)
\subset\cC,\]
and $\cL\cL = \cL$ implies that
\[ \pi(\cA)\cdot\big(U_{\cL}(\cL)\pi(\cA)\big)
\subseteq \pi(\cA)U_{\cL}(\cL)U_{\cL}(\cL) \pi(\cA)\subset \cC.\]
Thus $\pi(\cA)\cdot \cC\subseteq \cC$, hence
$\pi$ defines a homomorphism
$\eta_\cA \: \cA \to M(\cC)$ of $C^*$-algebras with
$\eta_\cA(A)L = \pi(A)L$ for $A \in \cA$ and $L \in \cC$.

Next we show that $U_{\cL}(\cL)\subeq M(\cC)$:
\[ U_{\cL}(\cL)\cdot\big(U_{\cL}(\cL)\pi(\cA)\big)\subseteq
U_{\cL}(\cL)\pi(\cA), \]
and $\cA = \cA\cA$ leads to
\[ U_{\cL}(\cL)\cdot\big(\pi(\cA)U_{\cL}(\cL)\big)
\subeq U_{\cL}(\cL)\pi(\cA)\pi(\cA) U_{\cL}(\cL)
\subeq \cC.\]
Thus $U_{\cL}(\cL){\cC}\subseteq\cC$, and
we obtain a homomorphism
$\eta_\cL \: \cL \to M({\cC})$ of $C^*$-algebras with
$\eta_\cL(L)C = U_{\cL}(L)C$ for $L \in \cL$ and $C \in \cC$.
Next, we verify equivalence of the conditions (i) to (iv).

(i) $\Rarrow$ (ii): Suppose that
$(\cC, \eta_\cA, \eta_\cL)$ is a crossed product host
of $(\cA, G, \alpha)$ and $(\cL, \eta)$.  In view of Proposition~\ref{prop:a.6},
the subset $\eta_\cL(\cL)\eta_\cA(\cA)$ spans a dense subspace of
$\cC$, which implies that
\[ \pi(\cA)U_{\cL}(\cL) =
\rho(\eta_\cA(\cA)\eta_\cL(\cL))
\subeq \rho(\cC) \subeq
\br\rho(\eta_\cL(\cL)\eta_\cA(\cA)).
= \br{U_{\cL}(\cL)\pi(\cA)}.\subseteq {U_{\cL}(\cL)\cB(\cH)}.\]
Notice that this deduction works for any covariant pair $(\pi,U)$ obtained
from a representation $\rho$ of $\cC$. We will need this observation in (c).
Notice that we have in fact shown that (i) implies
\begin{equation}
\label{equivCrossrep}
\pi(\cA)U_{\cL}(\cL) \subeq
\br\rho(\eta_\cL(\cL)\eta_\cA(\cA)).
=\br\rho(\eta_\cL(\cL)\eta_\cL(\cL)\eta_\cA(\cA)).
= \br{U_{\cL}(\cL)\rho(\cC)}.
={U_{\cL}(\cL)\rho(\cC)}
\end{equation}
which implies (ii), hence it is also an equivalent condition.

(ii) $\Rarrow$ (iii): Since $\|E_j\|\leq1$, the set of all
operators $B \in \cB(\cH)$ with
$\|U_{\cL}(E_j) B - B\| \to 0$ is a closed subspace
which obviously contains $U_{\cL}(\cL)\cB(\cH)$.
Now (ii) implies that it also contains $\pi(\cA) U_{\cL}(\cL)$.

(iii) $\Rarrow$ (iv) is trivial because every $C^*$-algebra possesses an
approximate identity (\cite[Thm.~3.1.1]{Mu90}).

(iv) $\Rarrow$ (i): Next, we want to verify conditions (CP1) to (CP4) for $(\cC, \eta_\cA, \eta_\cL)$.
By definition (CP1) and (CP4) hold. To verify (CP2), i.e. that $\eta_\cL$ is non-degenerate,
note that by condition (iv) and the trivially true
\[ \| U_{\cL}(E_j)  U_{\cL}(L)\pi(A)
-  U_{\cL}(L)\pi(A) \| \to 0 \quad \mbox{ for } \quad
A \in \cA, L \in \cL,\]
we obtain  $\| U_{\cL}(E_j) C
-  C \| \to 0$ for all $C\in\cC$, from which it is obvious
that  $\eta_\cL(\cL)\cC$ is dense in $\cC$.
Condition (CP3) is an immediate consequence of the
definitions and the non-degeneracy of $U_{\cL}$. This proves (i).

(c) It is clear that  ${\big(\cC/{\cal J},\tilde\Phi\circ\eta_\cA,\tilde\Phi\circ\eta_\cL\big)}$
satisfies (CP1), (CP3) and (CP4).
Let \break $\rho:\cC/{\cal J}\to \cB(\cH)$ be a faithful representation. Then
$\rho\circ\Phi:\cC\to \cB(\cH)$ is a nondegenerate representation, hence
$\big(\widetilde{(\rho\circ\Phi)} \circ \eta_\cA, \widetilde{(\rho\circ\Phi)} \circ \eta_G\big)\in\Rep_\cL(\alpha,{\cal H})$.
As $\cC$ is a crossed product host, we conclude from
(CP2) (see also (i) $\Rarrow$ (ii) above) that
\[
\big(\widetilde{(\rho\circ\Phi)} \circ \eta_\cA\big)(\cA)
\cdot \big(\widetilde{(\rho\circ\Phi)} \circ \eta_\cL\big)(\cL) \subeq
\big(\widetilde{(\rho\circ\Phi)} \circ \eta_\cL\big)(\cL)\cdot \cB(\cH),
\]
hence
\[ \big(\tilde\rho\circ \tilde\Phi\circ\eta_\cA\big)(\cA)
\cdot \big(\tilde\rho\circ \tilde\Phi\circ\eta_\cL\big)(\cL)\subeq
\big( \tilde\rho\circ\tilde\Phi\circ\eta_\cL\big)(\cL)\cdot \cB(\cH).\]
This implies via (b) that
\[ \rho(\cC/{\cal J})
=C^*\big((\tilde\rho\circ\tilde\Phi\circ\eta_\cA)(\cA)
\cdot(\tilde\rho\circ\tilde\Phi\circ\eta_\cL)(\cL)\big)\]
 is a crossed product host, and as $\rho$ is faithful,
a property which is inherited  by its extension $\tilde\rho$
to $M(\cC/{\cal J})$,
 it follows that ${\big(\cC/{\cal J},\tilde\Phi\circ\eta_\cA,\tilde\Phi\circ\eta_\cL\big)}$
 is a crossed product host.
\end{prf}

The preceding theorem shows how crossed product hosts can be constructed.
It  also isolates a distinguished class of representations:

\begin{defn}
\mlabel{crossrepDef}
Let $\alpha \: G \to \Aut(\cA)$ be a $C^*$-action and $(\cL,\eta)$ be a host algebra for $G$. Then
a covariant $\cL$-representation
$(\pi, U)\in\Rep_\cL(\alpha, \cH)$ is called a
{\it cross representation} for $(\alpha,\cL)$
if any of the equivalent conditions (b)(i)-(iv) of
Theorem~\ref{thm:5.1} hold.
We write $\Rep_\cL^\times(\alpha, \cH)$ for the set of cross representations for  $(\alpha,\cL)$
on $\cH$.
\end{defn}
These conditions are easily verified for covariant representations
in the usual case where $G$ is locally compact, $\cL = C^*(G)$ and $\alpha$ is strongly continuous (cf.\ Lemma~\ref{lem:a.1}(ii)(b)).
\begin{prop}{\rm(Permanence properties of cross representations)} \mlabel{lem:perm}
For cross representations for $(\cL, \alpha)$ of $(\cA,G,\alpha)$,
the following assertions hold:
\begin{description}
  \item[\rm(i)] Finite direct sums of cross representations are cross representations.
  \item[\rm(ii)] Subrepresentations of cross representations are cross representations.
  \item[\rm(iii)] Arbitrary multiples of cross representations are cross representations.
  \item[\rm(iv)] If $\cB\subset\cA$ is an $\alpha$-invariant $C^*$-subalgebra, then a
cross representation $(\pi,U)$ of $(\cL, \alpha)$ restricts on $\cB$ to a
cross representation for $(\cL,\alpha\!\restriction\!\cB)$.
\item[\rm(v)] Let $(\phi, \psi) \:
(\cL_H,\eta^H) \to (\cL_G,\eta^G)$ be a morphism of host algebras
(cf. Def.~\ref{def:morhost}) and define the  $C^*$-action
$\beta:=\alpha\circ\phi:H\to\Aut(\cA)$.
 Let $(\pi, U)$ be a covariant representation of $(\cA,G,\alpha)$.
If $(\pi, U\circ\varphi)$ is a cross representation for
$(\cL_H, \beta)$, then $(\pi, U)$ is a cross representation for $(\cL_G, \alpha)$.
%
\item[\rm(vi)] If $\cA$ is nonunital, let $\cA_1=\cA\oplus\C$ be the algebra $\cA$ with the identity adjoined.
Extend $\alpha_G$ to $\cA_1$ by setting $\alpha_g(\1)=\1$ for all $g\in G$ and extend
 covariant $\cL$-representations $(\pi, U)\in\Rep_\cL(\alpha,{\cal H})$ to $\cA_1$ by setting
 $\pi(\1)=\1$. Then $(\pi, U)$ is a cross representation for  $(\alpha,\cL)$ if and only if its
 extension to $(\cA_1,G,\alpha)$ is a cross representation.
\end{description}
\end{prop}

\begin{prf} (i)-(iv) are immediate consequences
of Theorem~\ref{thm:5.1}(b)(iii).

(v) Denote $U \circ \phi=:U^\phi$.
Since $(\pi, U^\phi)$ is a cross representation, we have
\[
\pi(\cA) U^\phi_{\cL_H}(\cL_H) \subeq U^\phi_{\cL_H}(\cL_H) \cB(\cH).
\]
The non-degeneracy of the morphism $\psi \: \cL_H \to \cL_G$ of
$C^*$-algebras further implies that
$\cL_G = \psi(\cL_H) \cL_G$.
Now
\[ U^\phi=U\circ\phi=\eta_G^*(\tilde U_{\cL_G})\circ\phi=\tilde U_{\cL_G}\circ\eta_G\circ\phi
=\tilde U_{\cL_G}\circ\tilde\psi\circ\eta_H=\big(\widetilde {U_{\cL_G}\circ\psi}\big)\circ\eta_H
=\eta^*_H(U_{\cL_G}\circ\psi) \] using
$\eta_G\circ\phi=\tilde\psi\circ\eta_H.$
An application of $(\eta^*_H)^{-1}$ produces the relation
$U_{\cL_G} \circ \psi = U^\phi_{\cL_H}$. Using these facts, we obtain:
\[ \pi(\cA) U_{\cL_G}(\cL_G)
= \pi(\cA) U^\phi_{\cL_H}(\cL_H) U_{\cL_G}(\cL_G)
\subeq U^\phi_{\cL_H}(\cL_H) \cB(\cH)
\subeq U_{\cL_G}(\cL_G) \cB(\cH).\]
Hence $(\pi, U)$ is a cross representation of $(\cL_G, \alpha)$.

(vi) It follows directly from (iv) that $(\pi, U)$ is a cross representation of $(\cL, \alpha)$
if its extension to $(\cA_1,G,\alpha)$ is a cross representation.  Conversely, assume that
$(\pi, U)$ is a cross representation of $(\cL, \alpha)$ on $\cH$, hence
$\pi(\cA) U_{\cL}(\cL) \subeq{U_{\cL}(\cL)\cB(\cH)}$ by Theorem~\ref{thm:5.1}(b)(ii).
Then
\[
\pi(\cA_1) U_{\cL}(\cL) =\big(\pi(\cA)+\C\1\big) U_{\cL}(\cL)
=\pi(\cA) U_{\cL}(\cL) + U_{\cL}(\cL) \subeq {U_{\cL}(\cL)\cB(\cH)}+U_{\cL}(\cL)=
{U_{\cL}(\cL)\cB(\cH)}
\]
from which it follows that the extension of $(\pi, U)$ to $(\cA_1,G,\alpha)$ is a cross representation.
\end{prf}

\begin{rem}
In Example~\ref{nonfullMax} below, we will see that $\Rep_\cL^\times(\alpha, \cH)$ is in general
a proper subset of $\Rep_\cL(\alpha, \cH)$. In particular,
Example~\ref{ex:non-cross-sum} shows that
infinite sums of cross representations need not be cross representations.
Unfortunately cross representations are not well behaved w.r.t. restriction to subgroups.
First, if $H\subset G$ is a closed subgroup, it is usually not clear what the host algebra
should be for~$H$. Even in the case of locally compact groups, if we take the natural
host algebras $C^*(G)$ and $C^*(H)$, then the restriction ${(\pi,U\!\restriction\!H)}$ of a cross representation
${(\pi,U)}\in\Rep_{C^*(G)}^\times(\alpha,{\cal H})$ does not in general produce a cross representation
of for $(\alpha\!\restriction\!H,C^*(H))$ (see Example~\ref{ex:noncross-restrict}).
\end{rem}

To construct a full  crossed product host,
according to Theorem~\ref{thm:5.1},
we must first  obtain the appropriate
covariant representation to construct it in.
Let  $G_d$ denote the group
$G$, endowed with the discrete topology, which is locally
compact, so that
we have a crossed product algebra $\cA \rtimes_\alpha G_d$
(cf.~Example~\ref{ex:1.1}).
Then ${\rm Rep}(\alpha,\cH)$
is in bijection with a subset of the set of representations of the
crossed product $\cA\rtimes_{\alpha} G_d$ on $\cH$.
The latter bijection preserves cyclicity and direct sums.
Any representation of $\cA\rtimes_{\alpha} G_d$ on $\cH$ can be
decomposed into cyclic components.
Thus each $(\pi,U)\in{\rm Rep}_\cL(\alpha,\cH)$ can
be decomposed into cyclic representations of
$\cA \rtimes_\alpha G_d$. It is clear that these components
are again continuous on $G$ and are $\cL$-representations of
$(\cA, G, \alpha)$ since we see from $U_\cL(G)'' = \pi(\cL)''$
(cf.~Proposition~\ref{prop:3.3}(ii)) that
subrepresentations of $\cL$-representations of $G$
are again $\cL$-representations.

\begin{defn} \mlabel{def:unicyc}
Cyclic representations of $\cA \rtimes_\alpha G_d$ are obtained from
states through the GNS-construction. Let $\fS_\cL$
denote the set of those states
$\omega$ on $\cA\rtimes_{\alpha} G_d$ which thus produce
a covariant $\cL$-representation
$(\pi_\omega,U_\omega)\in{\rm Rep}_\cL(\alpha,\cH_\omega).$
This allows us to define the {\it universal  covariant
$\cL$-representation}
$(\pi_u,U_u)\in{\rm Rep}_\cL(\alpha,\cH_u)$ by
\[
\pi_u:= \bigoplus_{\omega\in\fS_\cL}
\pi_\omega,\quad U_u:= \bigoplus_{\omega\in\fS_\cL} U_\omega
\quad\hbox{on}\quad  \cH_u=\bigoplus_{\omega\in\fS_\cL}\cH_\omega.
\]
Clearly $\cH_u = \{0\}$ if $\fS_\cL = \emptyset.$
We will use $(\pi_u,U_u)$ below to prove the existence of
crossed product hosts.  We use notation
 $U_{u,\cL} :=(\eta^*)^{-1}(U_u)\in \Rep(\cL,\cH_u)$
 for the associated representation of~$\cL$.
\end{defn}

Note that the following
theorem trivially holds if $\fS_\cL=\emptyset$
because in this case $\cC = \{0\}$ and the set of covariant
$\cL$-representations of $(\cA,\alpha)$ is empty.

\begin{thm} \mlabel{thm:exist} {\rm(Existence Theorem)} \\
Let $(\cL,\eta)$ be a host algebra for the topological group $G$
and $\alpha \: G \to \Aut(\cA)$ be a $C^*$-action. Then the following are equivalent:
\begin{description}
\item[\rm(i)] There exists a full crossed product host
$(\cC, \eta_\cA, \eta_\cL)$ for $(\alpha,\cL)$.
\item[\rm(ii)] The universal covariant $\cL$-representation $(\pi_u, U_u)$
of $(\cA,G,\alpha)$ on $\cH_u$ is a cross representation.
\item[\rm(iii)]
$\;\Rep_\cL(\alpha, \cH)=\Rep_\cL^\times(\alpha, \cH)\;$ for all Hilbert spaces $\cH$.
\end{description}
\end{thm}

\begin{prf} (i) $\Rarrow$ (ii): Let
$(\cC, \eta_\cA, \eta_\cL)$ be a full crossed product host
of  $(\alpha,\cL)$. Then
(CP5) implies the existence of a unique representation
$\rho \: \cC \to \cB(\cH_u)$ with
$\tilde\rho \circ \eta_\cA = \pi_u$ and
$\tilde\rho \circ \eta_\cL = U_{u,\cL}$. In view of Proposition~\ref{prop:a.6},
the subset $\eta_\cL(\cL)\eta_\cA(\cA)$ spans a dense subspace of $\cC$,
which implies that
\begin{align*}
\pi_u(\cA)U_{u,\cL}(\cL)
&=\rho(\eta_\cA(\cA)\eta_\cL(\cL))
\subeq \rho(\cC) \subeq
\br\rho(\eta_\cL(\cL)\eta_\cA(\cA)).\\
&= \br{U_{u,\cL}(\cL)\pi_u(\cA)}.\subseteq {U_{u,\cL}(\cL)\cB(\cH_u)}.
\end{align*}

(ii) $\Rarrow$ (iii): Since every  covariant $\cL$-representation $(\pi, U)$ is a
direct sum of cyclic ones, it is contained in some multiple of the
universal covariant $\cL$-representation $(\pi_u, U_u)$, hence a cross representation by
Lemma~\ref{lem:perm}.

(iii) $\Rarrow$ (i): If the universal covariant $\cL$-representation $(\pi_u, U_u)$ of $(\cA, G, \alpha)$
is a cross representation, then $\cC := C^*(\pi_u(\cA) U_{u,\cL}(\cL))$ is a crossed product host
by Theorem~\ref{thm:5.1}. We only need to verify that it is full,
i.e. that  $\Rep(\alpha,{\cal H})_{\eta_\times} = \Rep_\cL(\alpha,{\cal H})$
for each Hilbert space $\cH$. Since any covariant
$\cL$-representation $(\pi, U)$ of $(\cA,G,\alpha)$ can be embedded into a
multiple $(\pi_u^\kappa, U_u^\kappa)$ of the universal representation on
$\cH_u^\kappa := \ell^2(\kappa) \otimes \cH_u$,
there is a  projection $P_\cH\in \rho_u^\kappa({\cal C})' $ projecting onto
the subspace $\cH\subeq \cH_u^\kappa$. Then
 $P_\cH\cdot\rho_u^\kappa\in\Rep({\cal C},{\cal H})$ will produce
 $(\wt{P_\cH\cdot\rho_u^\kappa}) \circ \eta_\cA=P_\cH\cdot\wt{\rho_u^\kappa}\circ \eta_\cA
=P_\cH\cdot\pi_u^\kappa=\pi$ and
 $(\wt{P_\cH\cdot\rho_u^\kappa})
\circ \eta_\cL=P_\cH\cdot\wt{\rho_u^\kappa}\circ \eta_\cL=P_\cH\cdot U^\kappa_{u,\cL}=U_\cL$,
i.e. $\eta^*_\times(P_\cH\cdot\rho_u^\kappa)=(\pi, U)$.
This proves that $\eta^*_\times$ is surjective, i.e.\ that  (i) holds.
\end{prf}
\begin{rem}
(1) Below in Corollary~\ref{corAisAL} we will obtain additional equivalent conditions for
the existence of a full crossed product host.\\
(2)  In the usual case, $\cA \rtimes_\alpha G$ is a a full crossed product host
for $(\alpha,C^*(G))$. \\
(3) From condition (iii) in the preceding theorem,
we see that to prove nonexistence of a full crossed product host,
it suffices to display one covariant $\cL$-representation which is
not a cross representation (cf.\ Example~\ref{nonfullMax}).
Unlike the usual case, in general crossed product hosts need not
exist, and Theorem~\ref{thm:exist} characterizes when they do.
We will present examples below
beyond the usual case, where we do have existence.
\end{rem}

\begin{thm} \mlabel{thm:CHrelation} {\rm(Relation between crossed product hosts)}\\
Let $(\cL,\eta)$ be a host algebra for the topological group $G$,
and $(\cA, G, \alpha)$ be a $C^*$-action.
Let $(\cC, \eta_\cA, \eta_\cL)$ and
$(\hat\cC, \hat\eta_\cA, \hat\eta_\cL)$ be
crossed product hosts for $(\alpha,\cL)$,  such that
\[ \Rep(\alpha,{\cal H})_{\hat\eta_\times}\subseteq\Rep(\alpha,{\cal H})_{\eta_\times} \]
for any Hilbert space $\cH$ (cf.\ Theorem~\ref{Thm-HostCross}).
Then  there is a unique homomorphism $\Phi:\cC\to\hat\cC$ such that
$(\hat\cC, \hat\eta_\cA, \hat\eta_\cL)={\big(\Phi(\cC),\tilde\Phi\circ\eta_\cA,\tilde\Phi\circ\eta_\cL\big)}$,
i.e. $\hat\cC$ is a factor algebra of $\cC$. In particular, if $\cC$
is a full crossed product host, then all crossed product hosts are factor algebras of $\cC$.
\end{thm}

\begin{prf}
For the (faithful) universal representation $(\hat\rho_u, \hat\cH_u)$ of $\hat\cC$,
the corresponding covariant $\cL$-representation
$(\hat\pi, \hat{U})\in\Rep(\alpha,\hat\cH_u)_{\hat\eta_\times}$ of $(\cA, G, \alpha)$ is by hypothesis
in $\Rep(\alpha,\hat\cH_u)_{\eta_\times}$. Thus there is a representation $(\rho, \hat\cH_u)$ of $\cC$
such that $\eta^*_\times(\rho)=\big(\tilde\rho \circ \eta_\cA, \tilde\rho\circ \eta_G\big)
=(\hat\pi, \hat{U})$. Then by Theorem~\ref{thm:5.1}(a) we  have
\[
\rho(\al C.)=C^*\big(\hat\pi(\cA) \hat U_\cL(\cL)\big)=\hat\rho_u(\hat\cC)\cong\hat\cC\,.
\]
Thus $\Phi:=(\hat\rho_u)^{-1}\circ\rho$ is the required homomorphism.
\end{prf}
Thus the crossed product hosts for $(\alpha,\cL)$ can be partially ordered by the containments
$\Rep(\alpha,{\cal H})_{\hat\eta_\times}\subseteq\Rep(\alpha,{\cal H})_{\eta_\times}$
for any Hilbert space $\cH$.
\begin{ex}
\mlabel{FullCPH-discontAlph}
(A  full crossed product host for a discontinuous action)
Let $\cH = \ell^2(\Z)$ with the orthonormal basis
$(e_n)_{n \in \Z}$. Define a unitary representation
$U:\R\to \U(\cH)$ by $U_t e_n = e^{int}e_n$ for
$t \in \R$. Then
$\alpha_t (A) := U_t A U_t^*$ defines an action of $\R$ on
$ \cB(\cH)$.
Consider the unitary involution
$J\in \cB(\cH)$, defined by $Je_n := e_{-n}$. Then
$\alpha_t( J)e_n = e^{-2int} e_{-n}$
implies that $J$ does not have a continuous orbit map,
in particular,
the action $\,\alpha\,$ of $\,\R\,$ on $\,B(\cH)\,$ is not strongly continuous.
Define $\cA:=C^*\{\alpha_t( J)\mid t\in\R\}$, then by construction,
$\alpha \: \R \to \Aut(\cA)$ is not strongly continuous.
We will show that it has a full crossed product host,
where $\cL:=C^*(\R)$, even though
the usual crossed product is undefined.

Let $(\pi_u, U_u)$ denote the universal covariant representation
of $(\cA, \R, \alpha)$.
From the relation $U_t J  = J U_{-t}$ for $t \in \R$, we derive
that $\alpha_t(J) = J U_{-2t} = U_{2t} J$ in any covariant
representation. For $f \in C_c(\R)$ we thus obtain
\[ \pi_u(J) U_u(f) = U_u(\check f) \pi_u(J) \quad \mbox{ for }
\quad \check f(t) := f(-t) \]
and hence that
$\pi_u(J) U_{u,\cL}(\cL) \subeq
{U_{u,\cL}(\cL)\cB(\cH_u)}$. By applying $\Ad(U_s)$ to both sides, we also get that
$\pi_u(\alpha_s(J)) U_{u,\cL}(\cL) \subeq
{U_{u,\cL}(\cL)\cB(\cH_u)}$ for all $s$, hence condition (ii) in
Theorem~\ref{thm:exist} is satisfied
and so we have a full crossed product host where $G$
does not act continuously on $\cA$ .

The set
$U_u(\cL)\cB(\cH_u)$ coincides with the set $\cB(\cH_u)_c^L$ of all those operators
$Q\in \cB(\cH_u)$ for which the map
\[ \R \to \cB(\cH_u), \quad t \mapsto U_u(t)Q\]
is continuous (cf. Lemma~\ref{lem:a.1}(ii)(b) below).

The structure of $\cA$ is quite simple.
From $JU_tJ = U_{-t}$ it follows that
$\alpha_t(J) = U_t J U_{-t} = U_{2t}J$. As $J$ is contained in $\cA$,
it follows that $\cA$ is the $C^*$-algebra generated by
the unitary one-parameter group $(U_t)_{t \in \R}$ and $J$.
Hence $\cA$ is a homomorphic image of the crossed product
$C^*(\R_d) \rtimes (\Z/2\Z)$, where the $2$-element group
acts on the discrete group $\R_d$ by inversion.
\end{ex}


Our next example, is the important case of
  $\R$-actions on $\cB(\cH)$ produced
by selfadjoint operators $A$ on $\cH$ via $t\mapsto \Ad(\exp(itA))$. To characterize continuity of these actions,
we need:

\begin{prop} \mlabel{prop:c.4}
Let $(U,\cH)$ be a unitary representation of the metrizable locally compact abelian
group $G$. Then the following are equivalent:
\begin{description}
  \item[\rm(i)] The conjugation action $\alpha_g(A) := U_gAU_g^*$
of $G$ on $\cB(\cH)$ is strongly continuous.
 \item[\rm(ii)] $U$ is norm continuous.
\end{description}
\end{prop}

\begin{prf} Clearly, the norm continuity of $U$ implies the strong
continuity of the conjugation action~$\alpha$.

To see the converse, let $\Sigma := \Spec(U) \subeq \hat  G$ and assume that
$U$ is not norm continuous, i.e., that the closed subset $\Sigma$ of $\hat G$ is not compact
(Lemma~\ref{lem:c.2}). That $G$ is metrizable is equivalent to $\hat G$ being countable
at infinity, i.e., $\hat G = \bigcup_{n \in \N} K_n$ is a union of an increasing sequence
of compact subsets $K_n$ satisfying $K_n \subeq K_{n+1}^0$ for $n \in \N$
(cf. \cite[Thm~29, p.95]{Mor08}).
Note that this
implies that any compact subset of $\hat G$ is contained in some~$K_n$.
Here we may w.l.o.g.\ assume that
each $K_n$ is an identity neighborhood.

We now claim the existence of a sequence $(\chi_n)_{n \in \N}$ in $\Sigma$ and a
compact identity neighborhood $C \subeq \hat G$ such that
\begin{description}
\item[\rm(a)] The sequence $(\chi_n \chi_{n+1}^{-1})_{n \in \N}$ is not equicontinuous.
\item[\rm(b)] The compact subsets $\chi_n C$, $n \in\N$, are pairwise disjoint.
\end{description}

First we choose $C$ such that $C^{-1}C \subeq K_1$.
Then we choose the $\chi_n$ inductively by starting with
some $\chi_1 \in \Sigma$ and picking an element
\[ \chi_{n+1} \in \Sigma \setminus (\chi_1 K_1 \cup \cdots \cup \chi_n K_n) \]
which is possible because $\Sigma$ is not compact.
Then $\chi_{n+1}\chi_n^{-1} \not\in K_n$ implies that the set \break
$\{ \chi_{n+1} \chi_n^{-1} \: n \in \N\}$ is not relatively compact, hence not
equicontinuous (\cite[Prop.~7.6]{HoMo98}).
Moreover, our construction implies that, for
$m < n$, we have $\chi_n C \cap \chi_m C =\eset$, because
$\chi_n \not\in \chi_m CC^{-1} \subeq \chi_m K_1 \subeq \chi_m K_{n-1}$.

Let $P$ denote the spectral measure of $U$.
Then we obtain mutually orthogonal $U$-invariant subspaces
$\cH_n := P(\chi_n C) \cH$ and consider the closed invariant subspace
$\cK := \hat\oplus_{n \in \N} \cH_n$ of $\cH$. The restriction of $U_g$ to this subspace
can be written as $U_g\res_{\cK} = U^1_g U^2_g$, where
$U^1_gv = \chi_n(g)v$ for $v \in \cH_n$ and
$U^2$ is a continuous unitary representation whose spectral measure
is supported by the equicontinuous subset  $C \subeq \hat G$. Therefore
$U^2$ is norm continuous.

Assume that $\alpha$ is strongly continuous. Since each $\alpha_g$ is an isometry
of $\cB(\cK)$, this implies that $\alpha$ defines a continuous action
of $G$ on $\cB(\cK)$. Then
$\alpha^1_g(A) := U^2_{g^{-1}} \alpha_g(A) U^2_g$ is also strongly continuous.
From $\chi_n \in \Sigma$ it follows that $\cH_n \not=\{0\}$, so that there exists a
sequence of unit vectors  $v_n \in \cH_n$. We now consider the contraction
$A \in \cB(\cK)$, defined by $Av := \sum_n \la v, v_n \ra v_{n+1}$.
From
\[ \alpha^1_g(A) = \sum_n \chi_{n+1}(g) \chi_n(g)^{-1} \la \cdot, v_n \ra v_{n+1} \]
it now follows that
\[ \|\alpha_g^1(A) - A\|
= \sup_n |\chi_{n+1}(g) \chi_n(g)^{-1} - 1|
= \sup_n |\chi_{n+1}(g) -  \chi_n(g)|. \]
Since the sequence $(\chi_{n+1} \chi_n^{-1})_{n \in \N}$ is not equicontinuous, this leads
to the contradiction that $\limsup_{g \to \1} \|\alpha_g^1(A) -A\| > 0$.
Therefore (i) implies (ii).
\end{prf}

\begin{ex}\mlabel{ex:noncrossrep}
Let $\cH$ be an infinite-dimensional separable Hilbert space,
$\cA := \cB(\cH)$, $G := \R$, $\cL = C^*(\R)$,
$H$ be an unbounded selfadjoint operator,
$U_t := e^{itH}$ and $\alpha_t(A) := U_t A U_t^*$.
From Proposition~\ref{prop:c.4} (see \cite[Exs.~3.2.36]{BR02} for $G = \R$) we know that
$\alpha$ is not strongly continuous. Note that, as $U$ is strong operator continuous,
$(\pi,U)\in\Rep_\cL(\alpha, \cH)$ where $\pi$ is the identical representation
$\pi(A) = A$ of $\cA$.
We claim that
$(\pi,U) \in\Rep_\cL^\times(\alpha, \cH)$ if and only if $(i\1-H)^{-1}\in\cK(\cH)$.
By Lemma~\ref{lem:c.5}, we have $(i\1-H)^{-1}\not\in\cK(\cH)$ if
one point in the spectrum of $H$ is in its essential spectrum.

Assume first that $(i\1-H)^{-1}\not\in\cK(\cH)$.
Then
\[ U_{\cL}(\cL)=U_{\cL}(C^*(\R))
=\{f(H)\mid f\in C_0(\R)\}=C^*( (i\1-H)^{-1})  \not\subset\cK(\cH).      
\]
If $(\pi,U)\in\Rep_\cL^\times(\alpha, \cH)$,
then for its crossed product host
$\al C.:=C^*\big(\pi(\cA) U_{\cL}(\cL)\big)$,
Theorem~\ref{thm:5.1}(b) implies that
$\cB(\cH) \subeq M(\cC)$, so that
$\cC$ is a closed two-sided ideal of $\cB(\cH)$.
As $\1\in \cB(\cH)$, we get that
$U_{\cL}(\cL)\subseteq\cC$, hence
$\cC$ is not contained in $\cK(\cH)$.
As $\cB(\cH)$ has no proper ideal greater than $\cK(\cH)$
($\cH$ is separable), it follows that $\cC =B(\cH)$.
From $\cC = U_\cL(\cL)\cC$ we now obtain
$\lim_{t \to 0} U_t C = C$ for every
$C \in \cC$ (Lemma~\ref{lem:a.1}(ii)(b)).
For $C = \1$, this shows that the one-parameter
group $(U_t)$ is norm continuous, contradicting
the unboundedness of~$H$.
This shows that
$(\pi,U)\not\in\Rep_\cL^\times(\alpha, \cH)$.

Conversely, if  $(i\1-H)^{-1}\in\cK(\cH)$ and hence
 $U_{\cL}(C^*(\R))\subseteq\cK(\cH)$, then $\cC=\cK(\cH)$,
 and the left and right ideals generated in $\cC$ by $U_{\cL}(C^*(\R))$
 is $\cC$, hence $(\pi,U)\in\Rep_\cL^\times(\alpha, \cH)$.
 This is another example where a crossed product host exists
 for a discontinuous action.
\end{ex}

\begin{ex}
\mlabel{nonfullMax}
We construct an example where non-trivial cross representations exist, but there is no
full crossed product host. Let $(\cA_j, G, \alpha^{(j)})$, $j=1,2$,
be  automorphic $C^*$-actions
and $\cL=C^*(\R)$. Assume that we have  nontrivial
\[
(\pi_1,U_1)\in\Rep_\cL^\times(\alpha^{(1)}, \cH_1)\quad
\hbox{and}\quad (\pi_2,U_2)\in\Rep_\cL(\alpha^{(2)}, \cH_2)\big\backslash\Rep_\cL^\times(\alpha^{(2)}, \cH_2)
\]
(cf.\ Example~\ref{ex:noncrossrep}).
Let $\cA:=\cA_1\oplus\cA_2$, and the homomorphism $\alpha \: \R\to \Aut(\cA)$
be given by
\[ \alpha_t(A_1\oplus A_2):=\alpha^{(1)}_t(A_1)\oplus\alpha^{(2)}_t(A_2)
\quad \mbox{ for } \quad A_i\in\cA_i.\]
Define $\tilde\pi_j:\cA\to B(\cH_j)$
by $\tilde\pi_j(A_1\oplus A_2)=\pi_j(A_j)$, $j=1,2$.
Then $(\tilde\pi_1,U_1)\in\Rep_\cL^\times(\alpha, \cH_1)$ is a non-trivial cross representation.
On the other hand,  $(\tilde\pi_2,U_2)\in\Rep_\cL(\alpha, \cH_2)\big\backslash\Rep_\cL^\times(\alpha, \cH_2)$,
hence by Theorem~\ref{thm:exist}, there is no full crossed product host.
\end{ex}

\section{Cross representations -- special cases.}
\mlabel{XrepSC}

The existence conditions for full crossed product hosts
in Theorem~\ref{thm:exist}  need to be checked on a case-by-case basis.
However, in a few special situations, other than the usual case, it
is possible to verify these conditions more generally.

\subsection{Cross representations via compact operators}

\begin{thm} \mlabel{compactExist}
Let $(\cA, G, \alpha)$ be a  $C^*$-action and
$(\cL,\eta)$ be a host algebra for $G$.
If ${(\pi,U)}\in\Rep_\cL(\alpha,{\cal H})$ satisfies
$\pi(\cA)U_{\cL}(\cL)\subseteq\cK(\cH)$, then
 $(\pi, U)$ is a cross representation for $(\alpha,\cL)$.
 This holds in particular if $U_{\cL}(\cL)\subseteq\cK(\cH)$.
\end{thm}

\begin{prf} By   Theorem~\ref{thm:5.1}(b)(ii),
$(\pi, U)\in\Rep_\cL(\alpha,{\cal H})$
is  a cross representation for $(\alpha,\cL)$ if and only if
 $\pi(\cA) U_{\cL}(\cL) \subeq{U_{\cL}(\cL)\cB(\cH)}$. By hypothesis
 we have $\pi(\cA) U_{\cL}(\cL) \subeq\cK(\cH)$.
As the representation $U_{\cL}$ of $\cL$
is nondegenerate, Lemma~\ref{lem:a.1}(ii)(a) shows that
$\cK(\cH) \subeq U_\cL(\cL) \cB(\cH)$.
Hence  $(\pi, U)$ is  a cross representation for $(\alpha,\cL)$.
\end{prf}

It is easy to obtain examples where $U_{\cL}(\cL)\subseteq\cK(\cH)$, e.g. if
$G=\R$, $\cA=B(\cH)$, $\cL=C^*(\R)$, $\alpha_t:={\rm Ad}(e^{itH})$ for an unbounded operator
$H$ with compact resolvent, then in the defining representation
 $U_{\cL}(\cL)=C^*((i\1-H)^{-1})\subseteq\cK(\cH)$ (cf.\ Lemma~\ref{lem:c.5}
 and Example~\ref{ex:noncrossrep}).
Here is an example of the more general condition $\pi(\cA) U_{\cL}(\cL) \subeq\cK(\cH)$
for a discontinuous action:

\begin{ex}\mlabel{ex:AL-compact}
Consider the Schr\"odinger representation
on $L^2(\R)=\cH$. That is, let the operators $Q,\, P$
act on  the common invariant core consisting of the space of Schwartz functions
where $Q$ is multiplication  by the coordinate function $x$,
and  $P= i{d\over dx}$. 
Let $\cB:={(i\1-Q)^{-1}B(\cH)(i\1+Q)^{-1}}$,
and let $\alpha \: \R \to \Aut(B(L^2(R)))$
be $\alpha_t:={\rm Ad}(e^{itP^2})$. We now define
\[
\cA:=C^*(\alpha_{\R}(\cB))=C^*\Big(\mathop{\bigcup}_{t\in\R} (i\1-Q+2tP)^{-1}B(\cH)\,(i\1+Q-2tP)^{-1} \Big)
\]
which is clearly $\alpha\hbox{--invariant}$ and let  $(\pi, U)$
be the defining representation, so if $\cL=C^*(\R)$ then
$U_{\cL}(\cL)=C^*((i\1-P^2)^{-1})\not\subseteq\cK(\cH)$
(cf.\ Lemma~\ref{lem:c.5}).
However $\cA U_{\cL}(\cL)\subseteq\cK(\cH)$ hence $(\pi, U)$ is a cross representation
for $(\alpha\!\!\restriction\!\!\cA,\cL)$ by Theorem~\ref{compactExist}.
To see this, recall that $f(Q)h(P)\in \cK(\cH)$ if $f,\,h\in C_0(\R)$,
hence
${(i\1-Q)^{-1}(i\1-P^2)^{-1}}\in \cK(\cH)$ and so
\[
(i\1-P^2)^{-1}{(i\1-Q+2tP)^{-1}}={e^{itP^2}(i\1-P^2)^{-1}(i\1-Q)^{-1}e^{-itP^2}}\in \cK(\cH).
\]
However, the action  $\alpha \: \R \to \Aut(\cA)$ is not strongly continuous.
In fact, as $\cA$ contains
\[ C^*((i\1 - Q)^{-1} C_0(Q) (i\1+ Q)^{-1})
= C^*((Q^2+1)^{-1} C_0(Q)) = C_0(Q),\]
 we have
$(i\1-Q)^{-1}\in\cA$, and this element satisfies
\[
\big\|\alpha_t((i\1-Q)^{-1})-(i\1-Q)^{-1}\big\|=\|(i\1-Q+2tP)^{-1}-(i\1-Q)^{-1}\|\geq 1
\]
for all $t\not=0$ by \cite[Thm.~5.3(ii)]{BG08}.
\end{ex}

\begin{rem} \mlabel{rem:5.10}
(1) Theorem~\ref{compactExist} and Example~\ref{ex:AL-compact}
suggest that for a given algebra $\cN\subseteq \cB(\cH)$
we should study its {\it compactifier}, i.e. the $C^*$-algebra
\[
M_\cK(\cN):=\{A\in \cB(\cH)\,\mid\,A \cN\subseteq\cK(\cH)\quad\hbox{and}\quad
A^* \cN\subseteq\cK(\cH)\}
\]
which clearly contains the compacts, but it may contain more if $\cN$ is nonunital.
 So in the situation under discussion above,
 the invariant hereditary subalgebras of $\cB(\cH)$ determined by elements of the compactifier of
 $U_{\cL}(\cL)$ are algebras for which this is a cross representation.

(2) A particularly useful observation is the following. Let $G$ be a locally compact group
with a closed subgroup $H$, then the canonical homomorphism $\zeta \: C^*(H)\to M(C^*(G))$
is nondegenerate, so that
$C^*(G) = \zeta(C^*(H))C^*(G)$.
Hence,  if for ${(\pi,U)}\in\Rep_\cL(\alpha,{\cal H})$ with $\cL={C^*(G)}$,
 $\tilde U_\cL(C^*(H))$ is contained in the compactifier of $\pi(\cA)$, then ${(\pi,U)}$ is a cross
representation for ${(\alpha,\cA)}$ because
\[\pi(\cA)U_{\cL}(\cL) = \pi(\cA)U_{C^*(H)}(C^*(H)) U_{\cL}(\cL)
\subeq \cK(\cH) U_{\cL}(\cL)  \subeq \cK(\cH). \]
An easy but interesting example where
this is put to good use, is given below in Example~\ref{ex-SpX}.
\end{rem}

\begin{prop} If $G$ is a connected reductive Lie group with compact center and
$(U,\cH)$ an irreducible unitary representation, then
$U_{C^*(G)}(C^*(G)) \subeq \cK(\cH)$.
\end{prop}

\begin{prf} See \cite[Thm.~14.6.10]{Wa92} for a reference.
This result can be obtained from Remark~\ref{rem:5.10}(2)
because
for any compact subgroup $K \subeq G$ the restriction $U\res_K$
is of finite multiplicity, so that Proposition~\ref{prop:c.8} shows that
$U_{C^*(K)}(C^*(K)) \subeq \cK(\cH)$.
\end{prf}

\begin{prop} {\rm(\cite[Thm.~X.4.10]{Ne00})} If $G$ is a finite-dimensional Lie group and
$(U,\cH)$ a unitary highest weight representation, then
$U_{C^*(G)}(C^*(G)) \subeq \cK(\cH)$.
\end{prop}

\subsection{Cross representations for compact groups}

The case where $G$ is a compact group is of particular importance.
For this, the following proposition helps to identify
cross representations.
For the following proposition we recall from Definition~\ref{def:rep-concepts}
the concept of spectrum and finite multiplicity of a representation of a compact group.

\begin{prop} \mlabel{prop:6.4}
Let $(U,\cH)$ be a continuous unitary representation of the
compact group $G$.
Then the following are equivalent:
\begin{description}
\item[\rm(i)] $\Spec(U)$ is finite or $U$ is of finite multiplicity.
\item[\rm(ii)] For the identical representation $\pi$ of $\cA = \cB(\cH)$ on $\cH$,
$\alpha_g(A) = U_g A U_g^*$ and $\cL = C^*(G)$,
the pair $(\pi, U)$ is a cross representation of
 $(\alpha,\cL)$.
\end{description}
\end{prop}

\begin{prf} Let $P_\chi$, $\chi \in \hat G$, denote the projections onto the isotypic
subspaces of $\cH$ and note that each
$P_\chi\cH=\cH_\chi\otimes\cM_\chi$ w.r.t.\ which
$U_g\restriction P_\chi\cH = \chi(g) \otimes \1$
 is a tensor product. Here $\cM_\chi$ is the multiplicity space.

For $\cL = C^*(G)$, the subspace
\[
\Spann\big\{ U_{\cL}(\cL) P_\chi\,\mid\,\chi\in\hat G\big\}=\Spann\big\{ (B(\cH_\chi)\otimes\1)\, P_\chi
\,\mid\,\chi\in\Spec(U)\big\}
\]
is dense in $U_{\cL}(\cL)$.
We conclude that every
$A \in \cB(\cH)U_{\cL}(\cL)$ is also contained in
$U_{\cL}(\cL)\cB(\cH)$ if and only if,
$\cB(\cH_\chi, \cH)=\cB(\cH)\,(\cB(\cH_\chi)\otimes\1)\, P_\chi \subeq U_{\cL}(\cL)\cB(\cH)$
for all $\chi\in \Spec(U)$.
For a finite subset $F \subeq \hat G$ we write
$P_F := \sum_{\chi'\in F} P_{\chi'}$. Considering the finite subsets of $\hat G$
as a directed set with respect to inclusion, we obtain
\begin{equation}\label{eq:fcond}
\lim_{F \subeq \hat G} P_F B = B\quad\hbox{for every}\quad B \in \cB(\cH_\chi, \cH).
\end{equation}
If, conversely, this condition is satisfied for
every $B \in \cB(\cH_\chi, \cH)$, then
the closedness of $U_\cL(\cL)\cB(\cH)$ implies that $B \in U_\cL(\cL)\cB(\cH)$.

(i) $\Rarrow$ (ii): Condition \eqref{eq:fcond} is trivially satisfied if $\Spec(U)$ is finite,
which is equivalent to the norm continuity of $U$ (Proposition~\ref{prop:c.8}).
It is also satisfied if $U$ is of finite multiplicity because this implies that
all isotypical subspaces $\cH_\chi$ are finite-dimensional.
In this case every $B \in \cB(\cH_\chi, \cH)$ is of finite rank, hence Hilbert--Schmidt, i.e.,
$\|B\|_2^2 = \sum_{\chi'} \|P_{\chi'} B\|_2^2  < \infty$,
which implies in particular that
\[ \|B- P_FB\|^2
\leq \|B- P_FB\|_2^2
= \|\sum_{\chi'\not\in F} P_{\chi'} B\|_2^2
= \sum_{\chi'\not\in F} \|P_{\chi'} B\|_2^2 \to 0.\]

(ii) $\Rarrow$ (i): If, conversely, $\Spec(U)$ is infinite and some $P_\chi\cH$ is infinite-dimensional,
then we pick an orthonormal sequence $(v_n)_{n \in \N}$ in $P_\chi\cH$ and an injective map
$\eta \: \N \to \Spec(U)$. For unit vectors $w_n \in \cH_{\eta_n}$ we then obtain a bounded operator
$B \in \cB(\cH_\chi, \cH) \subeq \cB(\cH)$ by
\[ Bv := \sum_n \la v, v_n \ra w_n. \]
Since $\|P_{\eta_n} B\| = 1$ implies $\lim_{\chi'\to \infty} P_{\chi'} B \not= 0$ we see that (ii) implies (i).
\end{prf}

\begin{rem} \mlabel{rem:Compact}
If $(\cA, G, \alpha)$ is a  $C^*$-action where $G$ is compact
with host chosen as  $\cL = C^*(G)$, and if for a covariant representation ${(\pi,U)}\in\Rep_\cL(\alpha,{\cal H})$, $U$ satisfies  Proposition~\ref{prop:6.4}(i), then ${(\pi,U)}$ is a cross representation.
The converse is not true, i.e. a cross representation ${(\pi,U)}$ need not satisfy Proposition~\ref{prop:6.4}(i), as can be seen by
taking an infinite multiple of a cross representation for which $\Spec(U)$ is infinite.
\end{rem}

Proposition~\ref{prop:6.4} allows us to construct a number of interesting examples.
\begin{ex}\mlabel{ex:noncrossrep2}  (Another example of a covariant
non-cross representations)

Let $(U,\cH)$ be a continuous unitary representation of the circle group
$G = \T \cong \R/2\pi \Z$ and $H$ be its infinitesimal generator, i.e.,
$U_t = e^{itH}$ for $t \in \R$. Then $\cH$ is the orthogonal direct sum
of the eigenspaces $\cH_n = \{ v \in \cH \: H v = n v_n\}$,
$n \in \Z$. If infinitely many of these are non-zero ($H$ is unbounded)
and some $\cH_n$ is infinite-dimensional ($n$ lies in the essential spectrum of $H$),
then Proposition~\ref{prop:6.4} implies that for $\cL = C^*(\T)$,
$\cA = \cB(\cH)$, the tautological representation
$\pi$ of $\cA$ on $\cH$ and $\alpha_t(A) := U_t A U_t^*$, the covariant
representation $(\pi, U)$ is not a cross representation, i.e.,
$(\pi,U) \not\in \Rep_\cL^\times(\alpha, \cH)$.
\end{ex}

\begin{ex}\mlabel{ex:non-cross-sum}
(Infinite sums of cross representations need not be cross representations)\\
We consider the compact group $G := \T$, $\cL := C^*(G)$ and the
$C^*$-algebra
\[ \cA := \{ (A_n)_{n \in \N} \mid A_n \in \cB(L^2(G)), \sup \|A_n\| < \infty \}
\subset\prod_{n=1}^\infty \cB(L^2(G)) \]
with componentwise operations, and norm $\|(A_n)_{n \in \N}\|=\sup \|A_n\| $.
Let $G$ act componentwise on $\cA$ by ${\alpha_g((A_n))_{n \in \N}} = (U_g A_n U_g^*)_{n \in \N}$, where
$U \: G \to \U(L^2(G))$ is the regular representation.

Then, for each $k$, we have a covariant
representation $(\pi_k,U_k)$ of $(\cA, G, \alpha)$ on $L^2(G)$ defined by $\pi_k((A_n)_{n \in \N}) = A_k$ and $U_k := U$.
Since the representation
$(U, L^2(G))$ is of finite multiplicity, Proposition~\ref{prop:6.4} implies that
$(\pi_k,U_k)$ is a cross representation.

We prove that the direct sum representation $(\pi^\oplus,U^\oplus):=\big(\mathop{\oplus}\limits_{k\in \N}\pi_k,\mathop{\oplus}\limits_{k\in \N}U_k\big)$
on $\cH := \mathop{\oplus}\limits_{k\in \N} \cH_k$, $\cH_k := L^2(G)$, is not a cross representation.
Recall  that we can identify the dual group $\hat{G}=\Z$ with the orthonormal basis $\{e_n\mid n\in\Z\}$  of $L^2(G)$, where
$e_n(t) := e^{2\pi i nt}$. Then the unitary map $L^2(G)\to \ell^2(\Z)$ defined by this basis
is the Fourier transform under which $ U_\cL(\cL)$ transforms to pointwise multiplication
by $C_0(\N)$, hence $ U_\cL(\cL)e_n=\C\cdot e_n$  and $U_ze_n=z^ne_n$ for $z \in \T$  and all~$n$.
Now define  $A := (A_k)_{k \in \N} \in \cA$ where
$A_k := \la \cdot, e_1 \ra e_k \in B(L^2(G))$. Then
${\pi^\oplus(A)\cdot U^\oplus_\cL(\cL)} = \C\pi^\oplus(A)$,
 and for $z \in \T$ we have
\begin{eqnarray*}
 \|U^\oplus_z\pi^\oplus(A) - \pi^\oplus(A)\| &=& \Big\|\bigoplus_{k=1}^\infty
 \la \cdot, e_1 \ra(z^ke_k-e_k)   \Big\| = \sup_{k \in \N}\left\|
 \la \cdot, e_1 \ra(z^k-1)e_k   \right\| = \sup_{k \in \N} | z^k-1|,
 \end{eqnarray*}
 and this is equal to $2$ for some $z$ arbitrary close to $1$ (e.g. $\exp(i\pi/n)$, $n\in\N$), thus
$U^\oplus_z\pi^\oplus(A)\not\to \pi^\oplus(A)$ for $z \to 1$. Therefore
Lemma~\ref{lem:a.1}(ii)(b) implies that  $(\pi^\oplus,U^\oplus)$ is not a cross representation.
\end{ex}

\begin{ex} \mlabel{ex:noncross-restrict} (Cross representation need not restrict to
 cross representations on subgroups)

(a) Consider the compact Lie groups $H := \SO_n(\R) \subeq G := \U_n(\C)$.
We consider the natural unitary representation $U:G\to \U\big(\cF(\C^n)\big)$ on the symmetric Fock space
\[
\cH=\cF(\C^n):=\bigoplus_{k=0}^\infty\otimes_s^k\C^n\,,\quad
\otimes_s^k\C^n\equiv\hbox{symmetrized Hilbert tensor product of
$k$ copies of $\C^n$}
\]
given by
$U_g\big(v_1\otimes_s\cdots\otimes_s v_k\big):=
\big(gv_1\otimes_s\cdots\otimes_s gv_k\big)$
(cf. Example~\ref{FockCross} below).
 Since the representation of
$G$ on the subspace $\otimes_s^k\C^n$ of homogeneous elements of degree $k$ is irreducible
(\cite[Prop.~IV.1.12]{Ne00}),
the representation $U$ is multiplicity free, hence defines a cross representation
of $(\cA,G,\alpha)$ for $\cA = B(\cH)$ and $\alpha_g(A) = U_g A U_g^*$ for $\pi = \id$
(Proposition~\ref{prop:6.4}).

We claim that the representation $(\pi, U\res_H)$ is not a cross representation of
$(\cA, H, \alpha\res_H)$. In view of Proposition~\ref{prop:6.4}, it suffices to show that
$U\res_H$ has infinite multiplicities. Under the Segal--Bargmann transform
$B \: L^2(\R^n) \to \cH$ (\cite[Prop.~XII.4.3]{Ne00}), the action of
$H$ corresponds to the natural representation on $L^2(\R^n)$ given by
$(U_h f)(x) = f(h^{-1}x)$. Therefore the factorization
 $L^2(\R^n) \cong L^2(\bS^{n-1}) \otimes L^2(\R^\times_+)$ and the triviality of the
$H$-action on the second factors implies that all $H$-multiplicities in $\cH$
are infinite.

(b) If we take the Schr\"odinger representation $U$ of the $3$-dimensional
Heisenberg group $G$ on $\cH = L^2(\R)$,
let $\cA=\cB(\cH)$, then the group $U(G)$ is generated by
$\exp(itQ)$ and $\exp(itP)$, and $U(G)$ acts  on $\cA$ by conjugation,
defining a  $C^*$-action $(\cA, G, \alpha)$.
As $U_{C^*(G)}(C^*(G)) = \cK(\cH)$, this is a cross representation for
$(\alpha, C^*(G))$. However, it is not a cross representation
for the subgroup generated by $\exp(itQ)$
as $Q$ has continuous spectrum (cf.~Example~\ref{ExHHop}).
\end{ex}

\begin{ex}
\mlabel{NumberOp1}
(Number operators for bosons). Let $\cH$ be a nonzero complex Hilbert space and define a
symplectic form
$\sigma:\cH\times\cH\to\R$ by $\sigma(x,y):={\rm Im}{\langle x,y\rangle}$ where
${\langle \cdot,\cdot\rangle}$ denotes the inner product.
Then $(\cH,\sigma)$ is a symplectic space over $\R$,
and we let $\Sp(\cH,\sigma)$ denote the group of symplectic transformations of it.
Note that
unitaries on $\cH$ define symplectic transformations.
For the quantum system based on this we choose for its field algebra $\cA$  the Weyl algebra $\ccr \cH,\sigma.$
(cf.~\cite{Ma68}).
It is defined through the generators
$\{\delta_f\mid f\in \cH\}$ and the Weyl relations
\[ \delta_f^*=\delta_{-f} \quad \mbox{ and } \quad
\delta_f\delta_g=e^{-i\sigma(f,g)/2}\delta_{f+g}
\quad \mbox{ for } \quad f,g \in \cH.\]
Define a $C^*$-action $\alpha:\Sp(\cH,\sigma)\to\Aut(\cA)$ by $\alpha\s T.(\delta_x):=\delta_{T(x)}$
We are interested
in particular one-parameter subgroups of
$\Sp(\cH,\sigma)$, and it is well--known that in general the action of  these via $\alpha$ is not strongly continuous
as ${\|\delta_x-\delta_y\|}=2$
if $x\not=y$. The simplest nontrivial unitary one-parameter group one can define on $\cH$
is just multiplication by $e^{it}$, $t\in\R$. This defines an action of the circle group
$\alpha:\T\to\Aut(\cA)$ by $\alpha\s z.(\delta_x):=\delta_{zx}$, $z\in\T\subset\C$,
which is not strongly continuous.

Let $(\pi, U)$ be a covariant regular representation of
$(\cA, \T, \alpha)$ on a Hilbert space $\cK$.
Then the generator $N := -i \frac{d}{dt}|_{t = 0} U_{e^{it}}$
of the one-parameter group defined by $U$
is identified with a number operator (cf.~\cite{Ch68}).
We claim that, if  $N$ is bounded from below, then
$(\pi, U)$ is a cross representation for
$\cL = C^*(\T) \cong C_0(\Z) \cong c_0$.
Let $P_n \in U_\cL(\cL)$ denote the projection onto the
$\T$-eigenspace
\[ \cH_n  := \{ v \in \cH \: (\forall z\in \T)\, U_z v = z^n v\}.\]
Since $\Spann \{ P_n \: n \in \N\}$ is dense in
$U_\cL(\cL)$, it suffices to show via Lemma~\ref{lem:a.1}(ii)(b) that
\[ \lim_{z \to \1} U_z \pi(A)P_n = \pi(A)P_n \quad \mbox{ for } \quad
A \in \cA.\]
Since the set of elements $A \in \cA$ satisfying this relation
is a closed subspace, we may w.l.o.g.\ assume that
$A = \pi(\delta_f)$ for some $0\not= f \in \cH$.
Let $\cH = \cH_0 \oplus \cH_1$, $\cH_0 = \C f$,
denote the corresponding orthogonal decomposition of $\cH$.
The Stone--von Neumann Theorem implies that the restriction
of $\pi$ to $\cA_0 := \oline{\Delta(\cH_0,\sigma\res_{\cH_0})}$ is a multiple
of the Fock representation $(\pi_0, \cF(\cH_0))$, where
$\cF(\cH_0)$ denotes the Fock space of $\cH_0$.
Accordingly, we write
\[ \cK \cong \cF(\cH_0) \otimes \cK_1, \]
where the elements $\delta_{tf}$, $t \in \R$, act trivially
on $\cK_1$. Since the Fock representation of $\cA_0$ is covariant
with respect to the corresponding number operator, we obtain a
tensor decomposition $U_z = U_z^0 \otimes U_z^1$
and likewise
\[ N = N_0 \otimes \1 + \1 \otimes N_1 \]
for the corresponding number operators, where $\sigma(N_0) = \N_0$.
This leads to a decomposition
\[ P_n = \sum_{m \in \N_0} P_m^0 \otimes P_{n-m}^1,\]
and since $N$ is bounded from below, this sum is finite, which
in turn leads to
\[ U_z  \pi(\delta_f) P_n
= \sum_{m = 0}^M  U_z^0 \pi_0(\delta_f) P_m^0 \otimes U_z^1 P_{n-m}^1
= \sum_{m = 0}^M  z^{n-m} U_z^0 \pi_0(\delta_f)
P_m^0 \otimes P_{n-m}^1.\]
Each $P_m^0$ is a rank-one projection, hence in particular compact.
Therefore $\pi_0(\delta_f) P_m^0$ is compact, and thus
$\lim_{z \to \1}  U_z^0 \pi_0(\delta_f) P_m^0 = \pi_0(\delta_f) P_m^0$.
This implies that $\lim_{z \to 1} U_z \pi(\delta_f) P_n
= \pi(\delta_f) P_n$ for every $n \in \N$, and thus
$(\pi, U)$ is a cross representation.

Note that the boundedness of the number operator $N$
from below implies that $(\pi, U)$ is a multiple
of the Fock representation (cf.\ \cite{Ch68}).
We will return to this issue when we consider spectral conditions in \cite{GrN12}.
\end{ex}

\begin{ex}
\label{ex-SpX}
Let $(X,\sigma)$ be a nondegenerate symplectic space over $\R$,
and let $G := \Sp(X,\sigma)$ denote the group of symplectic transformations of it.
The quantum system based on this has for its field algebra $\cA$ either the Weyl algebra $\ccr X,\sigma.$
(cf.~\cite{Ma68}) above, or the Resolvent Algebra $\rsl$ (cf.~\cite{BG08}), and there is a discontinuous action
$\alpha:G\to\Aut\cA$, which for $\ccr X,\sigma.$ is the one above.

Assume that $X$ is finite-dimensional, so that $G$ is locally compact.
Let  ${(\pi,U)}\in\Rep_\cL(\alpha,{\cal H})$ with $\cL={C^*(G)}$ be the Fock representation.
Then $G$ contains the one--parameter group $H\cong \T$ generated by the number operator.
As the resolvent
of the number operator is compact for $X$ finite-dimensional, it follows from
Remark~\ref{rem:5.10}(2) that
${(\pi,U)}$ is a cross
representation for ${(\alpha,\cL)}$ (cf.\ also Lemma~\ref{lem:c.5}).
If we restrict $\alpha$ to subgroups $H' \subeq G$, then
${(\pi,U)}$ remains a cross representation for the restriction, as long as
$H \subeq H'$.
\end{ex}

\subsection{Cross representations of semidirect products}

We need tools to analyze cross representations in terms of subsystems.
It is sometimes useful to analyze the existence conditions in Theorems~\ref{thm:exist} and
\ref{thm:5.1}(b) in terms of subgroups of $G$. For instance, in the context of an action on $\cA$,
we may know that we have crossed product hosts for some subgroups generating $G$,
and then we want to conclude that $G$ itself has a crossed product host. E.g. for a Lie group $G$ we can
easily analyze these existence conditions in terms of the resolvents of the generators of the one-parameter subgroups, and
then we want use this information to establish  the conditions for $G$.

We consider a semidirect product $G=N \rtimes_\varphi  H$ of locally compact groups.
This can be concretely realized as follows (cf.~\cite[Sect.~3.3]{Wil07}).
Let $G$ be a locally compact group with a closed normal subgroup
$N$ and a closed subgroup $H$ such that $N \cap H = \{ e \}$ and $G = NH.$
Define $\varphi : H \to¨ \Aut N$ by $\varphi(h)(n) = hnh^{-1}$, then $(n, h)\to nh$ is an isomorphism
of locally compact groups between $N \rtimes_\varphi  H$ and $G$.
This allows us to identify $C_c(G)$ with $C_c(N\times H)$ by $\tilde{f}(n,h):=f(nh)$.
Then there are  morphisms ${\eta_N \: C^*(N) \to M(C^*(G))}$
 and ${\eta_H \: C^*(H) \to M(C^*(G))}$ given  by convolution of measures as follows.
Denote the Banach space of finite regular Borel measures on $G$ by $\al M.(G)$, then
convolution and
involution in $\al M.(G)$ are given by
\[
\int_G f(t)\,d(\gamma*\nu)(t):=\int_G \int_G f(st)\,
d\gamma(s)\, d\nu(t)\qquad
\hbox{and}\qquad\int_G f(t)\,d\gamma^*(t):=
\overline{\int_G\overline{f}(t^{-1})\,d\gamma(t)}
\]
for $\gamma,\;\nu\in\al M.(G)$, $f\in C_0(G)$. We identify
an $f\in L^1(N)$ with the measure $f\,d\mu_N$,  where $\mu_N$ is the Haar measure
of $N$. Then convolution of this measure with the measures $h\, d\mu_G$,  $h\in L^1(G)$,
produces the multiplier action ${\gamma_N \: C^*(N) \to M(C^*(G))}$, and likewise we define
${\gamma_H \: C^*(H) \to M(C^*(G))}$. In terms
of the $L^1\hbox{--functions}$ alone, we obtain
\[
(p*f)(x)=\int_N p(y) f(y^{-1}x)\,d\mu_N(y)\qquad\hbox{and}\qquad
(f*p)(x)=\int_N\Delta_G(y^{-1})\,f(xy^{-1})\,p(y)\,d\mu_N(y)
\]
 for all $f\in L^1(G),$ $p\in L^1(N)$ and $x\in G$
$\mu_G$-a.e.,  where $\Delta_G$ is the modular function of $G$
(cf.\ \cite[Thm.~20.9]{HR63})
 and $\gamma_H(p)\cdot f=p*f$, $\;f\cdot\gamma_H(p)=f*p$.

Note that the conjugation action of $H$ on $N$ leads to a strongly
continuous $C^*$-action ${\beta:H\to\Aut C^*(N)}$ given by
\[
\beta_h(f)(n):= \sigma(h)^{-1}f(h^{-1}nh),
\quad f\in C_c(N),\,h\in H, \, n\in N,
\]
where the continuous homomorphism $\sigma:H\to (\R^\times_+,\cdot)$
is obtained from the uniqueness of the
Haar measure on $N$ and is given by
\[
\sigma(h)\,d\mu_N(h^{-1}nh)=d\mu_N(n) \quad \mbox{ for } \quad h\in H, n \in N
\]
(cf.~\cite[Prop.~3.11]{Wil07} with the substitution $\cA=\C$).

\begin{thm} \mlabel{SemiDirect}
Consider a semidirect product $G=N \rtimes_\varphi  H$ of locally compact groups,
and let ${\gamma_N \: C^*(N) \to M(C^*(G))}$ and
${\gamma_H \: C^*(H) \to M(C^*(G))}$ be the morphisms above.
\begin{description}
\item[\rm(i)] The triple $(C^*(G), \gamma_N, \gamma_H)$
is a full crossed product host for $(\beta,C^*(H))$, where
$\beta$ is the action of $H$ on $C^*(N)$ induced from the conjugation action on~$N$.
We have \\ $C^*(G)=C^*\big(\gamma_H (C^*(H))\cdot \gamma_N (C^*(N))  \big)=\br \gamma_H (C^*(H))\cdot\gamma_N (C^*(N)).$.
\item[\rm(ii)] Let $\alpha \: G \to \Aut(\cA)$
be a $C^*$-action and denote the restricted actions of $\alpha$ to $N$ and $H$ by $\alpha^N$ and $\alpha^H$ resp.
 Let ${(\pi,U)}\in\Rep_\cL(\alpha,{\cal H})$,
 $\cL=C^*(G)$, and assume that its restrictions
${(\pi,U\!\restriction\! N)}$ and ${(\pi,U\!\restriction\! H)}$
 are cross representations for  $(\alpha^N,C^*(N))$ and $(\alpha^H,C^*(H))$ resp.
Then  ${(\pi,U)}$ is  a
cross representation for $(\alpha,C^*(G))$.
\end{description}
\end{thm}

\begin{prf}
(i) The group algebra of a semidirect product
was a motivating example for the development of
twisted group algebras, cf.~\cite{BS70}, \cite{PR92}. It is known that
$C^*(G)\cong C^*(N) \rtimes_\beta H$, hence it is a full crossed product
host for $(\beta, C^*(H))$ (Example~\ref{ex:1.1}). Its representations
correspond to the covariant representations of $(C^*(N),H,\beta)$
which in turn correspond to the covariant unitary representations
of the pair $(N,H,\phi)$, and these are in one-to-one correspondence
to continuous unitary representations of $G$.
For the sake of completeness, we take a closer look at the corresponding
isomorphism.
The isomorphism $\Phi:C^*(N) \rtimes_\beta H\to C^*(G)$ is the extension of the natural
map
\[ C_c(H,C_c(N))\to C_c(N\times H) \cong  C_c(G).\]
As $C^*(N) \rtimes_\beta H$ is a crossed product host as in
Definition~\ref{def:a.5} for the action $(C^*(N),H,\beta)$ with
choice of host $C^*(H)$, it follows from (CP4) and Proposition~\ref{prop:a.6}(i)
that
\[
C^*(N) \rtimes_\beta H= C^*\big(\eta_{C^*(N)}(C^*(N))\cdot \eta_{C^*(H)}(C^*(H))\big)
=\br \eta_{C^*(N)}(C^*(N))\cdot \eta_{C^*(H)}(C^*(H)).\,.
\]
We prove that the
isomorphism $\Phi:C^*(N) \rtimes_\beta H\to C^*(G)$ satisfies
$\tilde\Phi\circ\eta_{C^*(N)}=\gamma_N\,$ and
$\,{\tilde\Phi\circ\eta_{C^*(H)}}=\gamma_H$ where
 ${\gamma_N \: C^*(N) \to M(C^*(G))}$ and
${\gamma_H \: C^*(H) \to M(C^*(G))}$ are given above.
Now $C^*(N) \rtimes_\beta H$ is the closure of the
span of those $f\in C_c(H,C_c(N))$
of the form $f=p\otimes q$, i.e. $f(h)(n) = p(h)\,q(n)$ for $p\in C_c(H)$, $q\in C_c(N)$.
Then
\[
\eta_{C^*(N)}(r)\cdot f=(r*p)\otimes q\quad\hbox{for}\quad r\in C^*(N),
\]
 because the product in  $C^*(N)$ is convolution, and likewise
\[
f \cdot \eta_{C^*(H)}(s)=p\otimes (q*s)\quad\hbox{for}\quad s\in C^*(H).
\]
Now $\Phi(f)(n,h):= f(n)(h)=p(n)\,q(h)$ leads to
\[
\Phi\big(\eta_{C^*(N)}(r)\cdot f\big)(n,h)=(r*p)(n)\, q(h)=\big(\gamma_N(r)\cdot (p\otimes q)\big)(n,h)
=\big(\gamma_N(r)\cdot f\big)(n,h)
\]
where the second equation comes from
\begin{eqnarray*}
(\gamma_H(r)\cdot f)(n,h)&=&
(r*f)(nh)=\int_Nf(y^{-1}nh)\,r(y)\,d\mu_N(y)=
\int_Np(y^{-1}n)\,r(y)\,d\mu_N(y)\cdot q(h)\\[1mm]
&=&(r*p)(n)\, q(h).
\end{eqnarray*}
Thus $\tilde\Phi\circ\eta_{C^*(N)}=\gamma_N\,$ and likewise $\,{\tilde\Phi\circ\eta_{C^*(H)}}=\gamma_H$.

(ii) As  ${(\pi,U)}\in\Rep_\cL(\alpha,{\cal H})$, the representation $U$ is continuous, hence its restrictions to
$N$ and $H$ satisfy ${(\pi,U\!\restriction\! N)}\in\Rep_{C^*(N)}(\alpha^N,{\cal H}) $ and ${(\pi,U\!\restriction\! H)}
\in\Rep_{C^*(H)}(\alpha^H,{\cal H}) $ resp. The representations
${U_{C^*(G)}}$, ${U_{C^*(N)}}$ and ${U_{C^*(H)}}$ are all obtained from $U$ by integration, and via
the inclusions
${\gamma_N \: C^*(N) \to M(C^*(G))}$ and ${\gamma_H \: C^*(H) \to M(C^*(G))}$
we also have ${U_{C^*(N)}}=\widetilde{U_{C^*(G)}}\circ\gamma_N$ and
${U_{C^*(H)}}=\widetilde{U_{C^*(G)}}\circ\gamma_H$. Thus by (i) we get
\[
U_{C^*(G)}\big(C^*(G)\big)=\br\widetilde{U_{C^*(G)}}\big( \gamma_H (C^*(H))\cdot\gamma_N (C^*(N))\big).=
\br U_{C^*(H)}(C^*(H))\cdot U_{C^*(N)}(C^*(N))..
\]
Now the
assumption of a cross representation means by
equation~(\ref{equivCrossrep}) that
\[ \pi(\cA)U_{\cL}(\cL) \subeq
{U_{\cL}(\cL)\rho(\cC)}=\br\rho(\eta_\cL(\cL)\eta_\cA(\cA)).\]
for the appropriate host $\cL$.
So by hypothesis we have that
\[ \pi(\cA) U_{C^*(N)}(C^*(N)) \subeq \br {U_{C^*(N)}(C^*(N))\pi(\cA) }. \]
and
\[ \pi(\cA) U_{C^*(H)}(C^*(H)) \subeq \br {U_{C^*(H)}(C^*(H))\pi(\cA) }..\]
Thus
\begin{eqnarray*}
\pi(\cA) U_{C^*(G)}(C^*(G))&=&\pi(\cA)\br U_{C^*(H)}(C^*(H))\cdot U_{C^*(N)}(C^*(N))   .\\[1mm]
&\subeq&\br U_{C^*(H)}(C^*(H))\,\pi(\cA)\, U_{C^*(N)}(C^*(N)).\\[1mm]
&\subeq&\br U_{C^*(H)}(C^*(H))\, U_{C^*(N)}(C^*(N))\,\pi(\cA).\\[1mm]
&\subeq&\br U_{C^*(G)}(C^*(G))\,\pi(\cA).\\[1mm]
&\subeq& U_{C^*(G)}(C^*(G))\,B(\cH)
\end{eqnarray*}
and hence ${(\pi,U)}$ is  a
cross representation for $(\alpha,C^*(G))$.
\end{prf}

\begin{rem}
The preceding theorem is quite useful, e.g. if $G$ is an abelian finite-dimensional
connected Lie group, then it is isomorphic to $\R^n\times\T^k$ for some $n$ and $k$.
Denote the one-parameter subgroups corresponding to the factors of this product by $G_j$.
 If a representation $U:G\to \U(\cH)$ is continuous, then
$U(C^*(G_j))=C^*(R_j(\R))$,
where $R_j(t)=(it-H_j)^{-1}$, $t\in\R,$ is the resolvent of the generator
 $H_j$ of the one-parameter group $U(G_j)$. Thus $U(C^*(G))$ is generated by the products
 $\{R_1(t_1)\cdots R_n(t_n)\mid t_1,\ldots, t_n \in\R\}$. Thus if ${(\pi,U)}\in\Rep_\cL(\alpha,{\cal H})$  where $\cL=C^*(G)$,
then by Theorem~\ref{SemiDirect}, it is enough to check that
$\pi(\cA) R_j(t) \subeq\br R_j(\R)\cB(\cH).$ for each $j$ and one $t$,
to conclude that
 that ${(\pi,U)}$  is a cross representation.
\end{rem}
A close reading of the proof of Theorem~\ref{SemiDirect}(ii) shows that it has the following generalization.
\begin{thm} \mlabel{GenSemiDirect}
Let $G$ be a topological group with subgroups $N,\,H$. Assume that $G,\,N,\,H$  have hosts
$\cL_G,\,\cL_N,\,\cL_H$ respectively, such that there are  homomorphisms ${\gamma_N \: \cL_N \to M(\cL_G)}$ and
${\gamma_H \: \cL_H \to M(\cL_G)}$ satisfying
\[
\cL_G=C^*\big(\gamma_H (\cL_H)\cdot \gamma_N (\cL_N)  \big)=\br \gamma_H (\cL_H)\cdot\gamma_N (\cL_N).
\]
and $\eta^G\!\restriction\! N=\tilde\gamma_N \circ\eta^N$ and  $\eta^G\!\restriction\! H=\tilde\gamma_H \circ\eta^H$.
 Let $\alpha \: G \to \Aut(\cA)$
be a $C^*$-action and denote the restricted actions of $\alpha$ to $N$ and $H$ by $\alpha^N$ and $\alpha^H$ resp.
 Let ${(\pi,U)}\in\Rep_{\cL_G}(\alpha,{\cal H})$,
 and assume that its restrictions
${(\pi,U\!\restriction\! N)}$ and ${(\pi,U\!\restriction\! H)}$
 are cross representations for  $(\alpha^N,\cL_N)$ and $(\alpha^H,\cL_H)$ resp.
Then  ${(\pi,U)}$ is  a
cross representation for $(\alpha,\cL_G)$.
\end{thm}

\begin{prf}
Fix a ${(\pi,U)}\in\Rep_{\cL_G}(\alpha,{\cal H})$.
From the hypotheses,
Theorem~\ref{thm:a.2}(ii) implies  that both $\gamma_N$ and $\gamma_H$ are nondegenerate,
hence they extend to strictly continuous homomorphisms \\ ${\tilde\gamma_N \:M( \cL_N) \to M(\cL_G)}$ and
${\tilde\gamma_H \:M( \cL_H) \to M(\cL_G)}$.
As
\[
U\!\restriction\! N={\widetilde{U_{\cL_G}}\circ\eta^G\!\restriction\! N}=\widetilde{U_{\cL_G}}\circ\tilde\gamma_N \circ\eta^N
=(\eta^N)^*\big( \widetilde{U_{\cL_G}}\circ \tilde\gamma_N  \big)
\]
it follows from the host property of $\cL_N$ that $(U\!\restriction\! N)_{\cL_N}=
{((\eta^N)^*)^{-1} \big((U\!\restriction\! N)_{\cL_N}\big)} =
 \widetilde{U_{\cL_G}}\circ\gamma_N$.
Likewise we have that  $(U\!\restriction\! H)_{\cL_H}=\widetilde{U_{\cL_G}}\circ\gamma_H$.
Thus  we get
\[
U_{\cL_G}\big(\cL_G\big)=\br\widetilde{U_{\cL_G}}\big( \gamma_H (\cL_H)\cdot\gamma_N (\cL_N)\big).=
\br (U\!\restriction\! H)_{\cL_H}(\cL_H)\cdot (U\!\restriction\! N)_{\cL_N}(\cL_N)..
\]
Now a direct transcription of the remaining part of proof of Theorem~\ref{SemiDirect}(ii)
with the natural substitutions completes the proof.
\end{prf}
This can be further generalized to more than two subgroups.

\section{Non-cross representations}
\mlabel{NonX}

In the previous sections we showed how to construct crossed product hosts
from cross representations. Here we want to start from a non-cross representation,
and see how best  to obtain a crossed product host.
The strategy is to reduce $\cA$ to a smaller subalgebra $\cA_\cL$
for which the given representation becomes a cross representation.
The motivation for this comes from physics, where it is common
to start from a physically important representation (e.g. the Fock representation),
or one has a class of covariant $\cL$-representations
which are important (e.g.\ for the Weyl algebra, the regular representations), and it may happen that these are not
cross representations for the given pair $(\alpha,\cL)$.

One may view this situation as a compatibility question between two quantum constraints.
Constraints of a quantum system may be considered as a restriction of the representations
of its field algebra $\cA$ which the system may realize.
In our case one constraint is the restriction to the given fixed class of representations, and the other is
restriction to the class of covariant cross representations of $(\alpha,\cL)$. If
these two classes are disjoint, the constraints are incompatible.
Concretely, the nondegeneracy condition in Theorem~\ref{thm:5.1}(b) is the crucial compatibility requirement, and
 its failure indicates the incompatibility
of some elements of $\cA$ with the  constraint of restricting the representations of $\cA$
to the required class of covariant $\cL$-representations of $(\cA,\alpha)$.
We therefore seek a maximal invariant subalgebra of $\cA$ for which the nondegeneracy
condition holds for the given representations.

Let $G$ be a topological group,
let $(\cL,\eta)$ be a host algebra for $G$ and let
 $\alpha \: G \to \Aut(\cA)$ be a group homomorphism, where
$\cA$ is a $C^*$-algebra. We assume that we are given a
set of cyclic covariant $\cL$-representations, hence that
we can take their direct sum, denoted $(\pi^\oplus,U^\oplus)$.
Thus, following  the construction in Theorem~\ref{thm:5.1}(b),
we put $\al C.:=C^*\big(\pi^\oplus(\cA) U^\oplus_\cL(\cL)\big)$,
which produces
 a triple ${(\cC, \eta_\cA, \eta_\cL)}$ such that
\begin{itemize}
\item[\rm(CP1)] $\eta_\cA \: \cA \to M(\cC)$ and
$\eta_\cL \: \cL \to M(\cC)$ are morphisms of $C^*$-algebras.
\item[\rm(CP4)] $\eta_\cA(\cA) \eta_\cL(\cL) \subeq \cC$
and $\cC$ is generated by this set as a $C^*$-algebra.
\end{itemize}

As $(\pi^\oplus,U^\oplus)$ need not be a cross representation,
(CP2)and (CP5) will fail in general.
Note that if (CP2) fails, then Theorem~\ref{thm:a.2}(iii) is lost,
and hence the covariance requirement (CP3) does not make sense,
as it uses the multiplier extension $\tilde\eta_\cL \: M(\cL)\to M(\cC)$.
Covariance will have to be expressed differently,
and our first task is to obtain an adequate covariance condition
to replace (CP3). We need this, as
  $(\cC, \eta_\cA, \eta_\cL)$ is defined from a covariant $\cL$-representation,
and we are only interested in obtaining covariant $\cL$-representations from $\cC$.

\begin{defn} \mlabel{def:5.1}
Assume (CP1) and (CP4) for $(\cC, \eta_\cA, \eta_\cL)$ in the context above.
For any Hilbert space $\cH$, a (nondegenerate) representation $\rho\in\Rep(\cC, \cH)$ is called an
{\it $\cL$-representation} if ${\tilde\rho\circ\eta_\cL:\cL\to \cB(\cH)}$ is a
nondegenerate representation of $\cL$.  We write
$\Rep_\cL(\cC, \cH)$ for the set of
$\cL$-representations of $\cC$ on $\cH$.
\end{defn}

If the defining representation $(\pi^\oplus,U^\oplus)$ of $(\cC, \eta_\cA, \eta_\cL)$ is a
covariant $\cL$-representation, then the defining representation of $\cC$ is
an $\cL$-representation. If (CP2) holds, then all non-degenerate
representations of $\cC$ are $\cL$-representations
(cf. Theorem~\ref{thm:a.2}(vii)).

\begin{prop} \mlabel{LrepsofC}
Assume (CP1) and (CP4) for
$(\cC, \eta_\cA, \eta_\cL)$ in the context above. Then the following
assertions hold:
\begin{itemize}
\item[\rm(i)] Any subrepresentation of an $\cL$-representation $\rho$
of $\cC$
is again an $\cL$-representation of $\cC$, and any orthogonal sum of $\cL$-representations
is again an $\cL$-representation of $\cC$.
\item[\rm(ii)] Let $\rho\in\Rep_\cL(\cC, \cH)$.
Then there is a unitary $\cL$-representation of $G$, denoted
$U^\rho:G\to \U(\cH)$ uniquely specified by
$U^\rho(g)\cdot(\tilde\rho\circ\eta_\cL)(L)
=(\tilde\rho\circ\eta_\cL)(\eta(g)L)$ for all $g\in G,$ $L\in\cL$.
Moreover, $U^\rho(G)''=\big((\tilde\rho\circ\eta_\cL)(\cL))''$.
\item[\rm(iii)] Let $\rho\in\Rep_\cL(\cC, \cH)$  satisfy the
covariance relation:
\begin{equation}
  \label{eq:corel}
U^\rho(g)\cdot(\tilde\rho\circ\eta_\cA)(A)\cdot U^\rho(g)^*
=(\tilde\rho\circ\eta_\cA)\big(\alpha_g(A)\big)
\qquad\hbox{for all}\qquad g\in G,\; A\in\cA.
\end{equation}
Then $\tilde\rho\circ\eta_\cA:\cA\to \cB(\cH)$ is nondegenerate, and ${(\tilde\rho\circ\eta_\cA,\,U^\rho)}
\in\Rep_\cL(\alpha, \cH)$.
\item[\rm(iv)] Let $\rho\in\Rep_\cL(\cC, \cH)$ be faithful and covariant, i.e. it satisfies \eqref{eq:corel}.
Then all $\cL$-representations of $\cC$ are covariant.
 If, in addition, $\cA$ is unital, then $\eta_\cA:\cA\to M(\cC)$
 is nondegenerate.
 \item[\rm(v)] Assume that there exists a
faithful covariant $\cL$-representation of $\cC$.
For each Hilbert space $\cH$ define a map
 $$ \eta^*_\times \: \Rep_\cL(\cC,{\cal H}) \to \Rep_\cL(\alpha, \cH),  \quad\hbox{given by}\quad
\eta^*_\times(\rho):=\big(\tilde\rho \circ \eta_\cA, U^\rho\big).$$
Then $ \eta^*_\times$ is injective. Moreover $ \eta^*_\times$ takes cyclic (resp irreducible)
 representations $\rho\in\Rep_\cL(\cC,{\cal H})$
to cyclic (resp. irreducible) covariant representations ${(\pi,\,U)}
\in\Rep_\cL(\alpha, \cH)$.
\end{itemize}
\end{prop}

\begin{prf}
(i) By definition ${(\tilde\rho\circ\eta_\cL)(\cL)\subset \rho(\cC)'' \subset \cB(\cH)}$, hence it preserves any $\rho(\cC)$-invariant subspace.
In particular, if there is a nonzero $v$ in a $\rho(\cC)$-invariant
subspace such that
${(\tilde\rho\circ\eta_\cL)(\cL)v}=\{0\}$, then ${(\tilde\rho\circ\eta_\cL)(\cL)}$ is degenerate,
which contradicts the assumption that $\rho$ is an $\cL$-representation. Thus all
subrepresentations of $\cL$-representations of
$\cC$ are $\cL$-representations. Likewise, if
$\rho\in\Rep(\cC, \cH)$ is an orthogonal sum of $\cL$-representations, and there is a nonzero $w\in\cH$  such that
${(\tilde\rho\circ\eta_\cL)(\cL)w}=\{0\}$, then for each of its components $w_i$ we must have
${(\tilde\rho\circ\eta_\cL)(\cL)w_i}=\{0\}$ which contradicts the assumption that $\rho$ is an orthogonal sum of $\cL$-representations.

(ii) As $\rho\in\Rep_\cL(\cC, \cH)$, the representation
${\tilde\rho\circ\eta_\cL:\cL\to \cB(\cH)}$ is
nondegenerate. By the host property of $\cL$, it defines
a unitary representation
$U^\rho := \eta^*(\tilde \rho \circ \eta_\cL) :G\to \U(\cH)$.
It is uniquely specified by
$U^\rho(g)\cdot(\tilde\rho\circ\eta_\cL)(L)
=(\tilde\rho\circ\eta_\cL)(\eta(g)L)$ for all $g\in G,$ $L\in\cL$.
That $U^\rho(G)''=(\tilde\rho\circ\eta_\cL)(\cL)''$
follows directly from Remark~\ref{vNalgHost}.

(iii) As $(\cL,\eta)$ is a host algebra for $G$,
the subgroup $\eta(G)$ spans a strictly dense subspace of $M(\cL)$
(Proposition~\ref{prop:3.3}(ii)).
Since ${\tilde\rho\circ\eta_\cL}$ extends to $M(\cL)$ by strictly continuity, we have
for $A \in \cA$ and $L \in \cL$:
\[
\rho\big(\eta_\cL(L)\eta_\cA(A)\big)
\in \Big({{(\tilde\rho\circ\eta_\cL)}(\Spann(\eta(G)))\cdot(\tilde\rho\circ\eta_\cA)(A)}\Big)^{\hbox{-- s-op}}
\subeq \lbr(\tilde\rho\circ\eta_\cA)(\cA)\cdot(\tilde\rho\circ\eta_\cL)(M(\cL))\rbr^{\hbox{-- s-op}}, \]
using covariance for the last inclusion, where closures are w.r.t.\ strong operator topology, and square brackets indicate span.
This implies that
\[
\rho\big(\eta_\cA(\cA)\eta_\cL(\cL)\big)\subeq
{\lbr(\tilde\rho\circ\eta_\cL)(M(\cL))\cdot(\tilde\rho\circ\eta_\cA)(\cA)\rbr^{\hbox{-- s-op}}}.
\]
If $(\tilde\rho\circ\eta_\cA)(\cA)v=0$ for some nonzero $v\in\cH$, then clearly
$\rho\big(\eta_\cL(\cL)\eta_\cA(\cA)\big) v=0$. On the other hand from the inclusion above we get
that $\rho\big(\eta_\cA(\cA)\eta_\cL(\cL)\big)v=0$ and hence that ${\rho(\cC)v=0}$
which contradicts the nondegeneracy of
$\rho$. Thus the representation $\tilde\rho\circ\eta_\cA:\cA\to \cB(\cH)$ is nondegenerate. By assumption of the covariance relation,
it follows that ${(\tilde\rho\circ\eta_\cA,\,U^\rho)}$ is a covariant representation, and by (ii) $U^\rho$ is an
$\cL$-representation, hence ${(\tilde\rho\circ\eta_\cA,\,U^\rho)}\in\Rep_\cL(\alpha, \cH)$.

(iv) Given that $\rho\in\Rep_\cL(\cC, \cH)$ is faithful and covariant, we define automorphisms
$\gamma_g\in\Aut\cC$ for each $g\in G$ as follows. For each $A \in \cA$ and $L \in \cL$ let
\[
\rho\big(\gamma_g\big(\eta_\cL(L)\eta_\cA(A)\big)\big):=
U^\rho_g\rho\big(\eta_\cL(L)\eta_\cA(A)\big)(U^\rho_g)^*=\rho\big(\eta_\cL(\beta_g(L))\eta_\cA(\alpha_g(A))\big)
\]
where $\beta_g(L):= \eta(g)L\eta(g)^*$ is the
automorphism of $\cL$ corresponding to $g \in G$.
In other words,
\[
\gamma_g\big(\eta_\cL(L)\eta_\cA(A)\big):=
\eta_\cL(\beta_g(L))\eta_\cA(\alpha_g(A)).
\]
Since $\rho$ extends to a faithful representation of $M(\cC)$, we also know that the extension
of $\gamma_g$ to $\tilde\gamma_g\in\Aut M(\cC)$ is given by $\tilde\rho\circ\tilde\gamma_g
=\Ad(U_g^\rho)\circ\tilde\rho$. Thus
\[
\tilde\gamma_g(\eta_\cA(A))= \eta_\cA(\alpha_g(A))\qquad\hbox{and}\qquad
\tilde\gamma_g(\eta_\cL(L))=\eta_\cL(\beta_g(L))\,.
\]
We also obtain identities such as
\begin{equation}
  \label{eq:ident1}
\eta_\cL(L)\eta_\cA(A)\eta_\cL(L')
=\eta_\cL(L\eta(g)^*)\eta_\cA(\alpha_g(A))\eta_\cL(\eta(g)L')
\end{equation}
from the definition of $U^\rho_g$ and the covariance condition.

Let  $\rho_1\in\Rep_\cL(\cC, \cH_1)$. We want to prove that it is covariant.
As ${\tilde\rho_1\circ\eta_\cL:\cL\to \cB(\cH_1)}$ is nondegenerate,
it suffices to prove that
\begin{align*}
& \Big((\tilde\rho_1\circ\eta_\cL)(L)\,v,\;
(U^{\rho_1}_g)^*(\tilde\rho_1\circ\eta_\cA)\big(\alpha_g(A)\big)U^{\rho_1}_g
\cdot(\tilde\rho_1\circ\eta_\cL)(L')\,w\Big)\\
&=\Big((\tilde\rho_1\circ\eta_\cL)(L)\,v,\;
(\tilde\rho_1\circ\eta_\cA)(A)
\cdot(\tilde\rho_1\circ\eta_\cL)(L')\,w\Big)
\end{align*}
for all $L,\,L'\in\cL$, $v,\,w\in\cH_1$, $g\in G$ and $A\in\cA$. The left hand side rearranges with \eqref{eq:ident1} to
\[
\Big(v,\;\tilde\rho_1\Big(\eta_\cL(L^*\eta(g)^*)\,
\eta_\cA\big(\alpha_g(A)\big)\eta_\cL(\eta(g)L')\Big)\,w\Big)
=\Big(v,\;\tilde\rho_1\big(\eta_\cL(L^*)\,
\eta_\cA(A)\eta_\cL(L')\big)\,w\Big).
\]
As this rearranges to the right hand side of our desired condition,
covariance is proven for $\rho_1$.

Since
 ${\big(\tilde\rho\circ\eta_\cA,U^\rho  \big)}$
 is a covariant representation, both
$\tilde\rho\circ\eta_\cA:\cA\to \cB(\cH)$ and
$\tilde\rho\circ\eta_\cL$
 are nondegenerate, and as $\cA$ is unital, we
 get that $(\tilde\rho\circ\eta_\cA)(\1)=\1$. Since $\tilde\rho:M(\cC)\to \cB(\cH)$ is faithful,
it follows that $\eta_\cA(\1)= \1$. In particular,
$\eta_\cA\: \cA\to M(\cC)$ is nondegenerate.

 (v) That the map $\eta^*_\times$ is well-defined follows from (iii) and (iv).
 If  $\rho_1,\,\rho_2\in\Rep_\cL(\cC,{\cal H})$ satisfy
 $\eta^*_\times(\rho_1)= \eta^*_\times(\rho_2)$,
i.e. ${\big(\tilde\rho_1 \circ \eta_\cA, U^{\rho_1}\big)}
 ={\big(\tilde\rho_2 \circ \eta_\cA, U^{\rho_2}\big)}$, then
 \[
 \rho_1(\eta_\cA(A)\eta_\cL(L))
= (\tilde\rho_1 \circ \eta_\cA)(A)\, U^{\rho_1}_\cL(L)
 =(\tilde\rho_2 \circ \eta_\cA)(A)\, U^{\rho_2}_\cL(L) = \rho_2(\eta_\cA(A)\eta_\cL(L))
 \]
for $A \in \cA, L \in \cL$, and this implies that
$\rho_1 = \rho_2$.  Thus $ \eta^*_\times$ is injective.

Recall that
 $(\pi,U)\in{\rm Rep}_\cL(\alpha,\cH)$ is cyclic (resp. irreducible)
if and only if it defines
a cyclic (resp. irreducible) representation of $\cA \rtimes_\alpha G_d$, which we denote by $\pi\times U.$
 Let $\eta^*_\times(\rho)=(\pi,U)$.
Since for any $\psi\in\cH$ we have
\[
\overline{\rho(\cC)\psi} = \overline{C^*\big(\pi(\cA)U_\cL(\cL)\big)\psi}
= \overline{C^*\big(\pi(\cA)U_\cL(\cL)''\big)\psi} = \overline{C^*\big(\pi(\cA)U(G)''\big)\psi}
=\overline{(\pi\times U)(\cA \rtimes_\alpha G_d)\psi},
\]
$\rho$ is cyclic if and only if $\pi\times U$ is cyclic.
A representation is irreducible if and only if every
nonzero vector is cyclic, so it is clear that $\rho$ is irreducible if and only if $\pi\times U$ is irreducible.
\end{prf}

The situation in (iv) is the natural one, where we construct $(\cC, \eta_\cA, \eta_\cL)$ from a covariant $\cL$-representation
of $(\cA, G, \alpha)$ as in Theorem~\ref{thm:5.1}(b).
In view of (i),  it makes sense to define the
universal $\cL$-representation of $\cC$:
\begin{defn} \mlabel{def:univLrepC}
Assume (CP1) and (CP4) for $(\cC, \eta_\cA, \eta_\cL)$ in the context above.
 Let $\fS_\cL(\cC):=\{\omega\in\fS(\cC)\,\mid\,\rho_\omega\in
\Rep_\cL(\cC, \cH_\omega)\}$
denote the set of those states,
whose GNS--representations $\rho_\omega$ are  $\cL$-representations of $\cC$.
We define the {\it universal
$\cL$-representation of $\cC$},  $\rho_u^\cC\in{\rm Rep}_\cL(\cC,\cH_u^\cC)$ by
\[
\rho_u^\cC:= \bigoplus_{\omega\in\fS_\cL(\cC)}
\rho_\omega,
\quad\hbox{on}\quad  \cH_u^\cC=\bigoplus_{\omega\in\fS_\cL(\cC)}\cH_\omega.
\]
Clearly $\cH_u^\cC = \{0\}$ if $\fS_\cL(\cC) = \emptyset.$
When there is no danger of confusion, we will omit the superscript $\cC$.
\end{defn}

\begin{lem} \mlabel{UnivLrep}
Let
$(\pi_u,U_u)\in{\rm Rep}_\cL(\alpha,\cH_u)$ be the the universal  covariant $\cL$-representation,
and define $\al C.:=C^*\big(\pi_u(\cA) U_{u,\cL}(\cL)\big)$ and $\eta_\cA$ and $ \eta_\cL$  as in
Theorem~\ref{thm:5.1}(b), producing the triple $(\cC, \eta_\cA, \eta_\cL)$
satisfying (CP1) and (CP4).
Then the defining representation of $\cC\subset\cB(\cH_u)$ is equivalent to the
universal $\cL$-representation $(\rho_u^\cC, \cH_u^\cC)$ of $\cC$
in Definition~\ref{def:univLrepC}.
\end{lem}
\begin{prf}
By Proposition~\ref{LrepsofC}(v) we have the injective map
 $ \eta^*_\times \: \Rep_\cL(\cC,{\cal H}_u^\cC) \to \Rep_\cL(\alpha, \cH_u^\cC) $
 which preserves cyclic components.
Hence $ \eta^*_\times (\rho_u^\cC)=\big(\tilde\rho_u^\cC \circ \eta_\cA, U^{\rho_u^\cC}\big)$
is the direct sum over the same cyclic components as $\rho_u^\cC$ with the same multiplicities.
The cyclic components of $\rho_u^\cC$ are the GNS-representations of $\fS_\cL(\cC)$,
and under $ \eta^*_\times$ these become GNS-representations of states in~$\fS_\cL$.
(Recall that the cyclic components of $(\pi_u,U_u)$ are the GNS-representations of states
in~$\fS_\cL$.)
Conversely, given $\omega\in\fS_\cL$, then there is a vector in $\cH_u$ which
will reproduce it, hence in the defining representation for $\cC$ this vector
will produce a state in $\fS_\cL(\cC)$. Thus we have the same cyclic components
with the same multiplicities,
so $(\pi_u,U_u)$ is equivalent to $ \eta^*_\times (\rho_u^\cC)$ and
hence $\rho_u^\cC$ is equivalent to $ (\eta^*_\times)^{-1} (\pi_u,U_u)$.
\end{prf}

Since in our natural examples, the defining representation of $\cC$ is an $\cL$-representation, it is natural
to require the universal $\cL$-representation  $\rho_u$ of $\cC$
to be faithful.
We can now state our covariance assumption:
\begin{defn} \mlabel{def:CP3'}
Given a triple $(\cC, \eta_\cA, \eta_\cL)$ satisfying (CP1) and (CP4) as above, then the {\it covariance condition} is
given by:
\begin{description}
\item[\rm(CP3')] The universal $\cL$-representation  $\rho_u$ is faithful and satisfies
\[
U^{\rho_u}(g)\cdot (\tilde\rho_u\circ\eta_\cA)(A)
\cdot U^{\rho_u}(g)^*=(\tilde\rho_u\circ\eta_\cA)\big(\alpha_g(A)\big)
\qquad\hbox{for all}\qquad g\in G,\; A\in\cA.
\]
\end{description}
\end{defn}

Note that by Proposition~\ref{LrepsofC}(iv), any triple $(\cC, \eta_\cA, \eta_\cL)$ obtained
as in Theorem~\ref{thm:5.1}(b) from a covariant $\cL$-representation will satisfy condition (CP3').
Moreover, if we replace $\rho_u$ by any other faithful $\cL$-representation, the resulting covariance condition
will be equivalent to (CP3').

\begin{rem}
\mlabel{actionMC}
If $(\cC, \eta_\cA, \eta_\cL)$ is a crossed product host,
then the action $\alpha \: G \to \Aut(\cA)$
extends to an action $\tilde\alpha \: G \to \Aut(M(\cC))$ by $\tilde\alpha_g=\Ad(\eta_G(g))$ where $\eta_G = \tilde\eta_\cL \circ \eta$.
On the other hand, if $(\cC, \eta_\cA, \eta_\cL)$ is not a crossed product
host, but it satisfies (CP1), (CP3') and (CP4), then via
$\Ad(U^{\rho_u}(g))$ we can define (in $\tilde\rho_u$) an automorphism of $\cB(\cH_u)$ which preserves $\cC$, $\eta_\cA(\cA)$ and $\eta_\cL(\cL)$
and coincides on $\eta_\cA(\cA)$ with $\eta_\cA\circ\alpha_g$ (see
proof of Proposition~\ref{LrepsofC}(iv) above).
Moreover, the automorphism on $\cC$ extends uniquely
to  $M(\cC)$, so it also preserves $M(\cC)$. On $\cC$ the automorphism $\gamma_g\in\Aut\cC$ is determined by
\[
\gamma_g\big(\eta_\cL(L)\eta_\cA(A)\big):=
\eta_\cL(\beta_g(L))\eta_\cA(\alpha_g(A))\qquad\hbox{for all}
\quad A \in \cA,\;L \in \cL,
\]
where $\beta_g(L):= \eta(g)L\eta(g)^*$ is the automorphism of $\cL$
determined by the multiplier action of $G$ on $\cL$. Moreover,  $\gamma_g$
also preserves the subalgebra $M(\cC)_\cL:=\eta_\cL(\cL) M(\cC)\eta_\cL(\cL)$, and if
the multiplier action $\eta \: G \to \U(M({\cal L}))$ is strictly continuous, then the restriction
of $\gamma_g$ to $M(\cC)_\cL$ is a strongly continuous action, as
${\gamma_g\big(\eta_\cL(L_1) N\eta_\cL(L_2)\big)}={\eta_\cL(\eta(g)L_1) N\eta_\cL(L_2\eta(g^{-1}))}$ for $L_i\in\cL$, $N\in M(\cC)$.
\end{rem}
By Proposition~\ref{LrepsofC}(iii) and (iv), if condition (CP3') holds, then each
$\cL$-representation $\rho\in\Rep_\cL(\cC, \cH)$ will produce a covariant pair
${(\tilde\rho\circ\eta_\cA,\,U^\rho)}\in\Rep_\cL(\alpha, \cH)$. Thus, given such a triple $(\cC, \eta_\cA, \eta_\cL)$,
it makes sense to seek a subalgebra which can ``carry'' the $\cL$-representations. Note that a representation
of $\cC$ is an $\cL$-representation if and only if
it is nondegenerate when we restrict it to the
$C^*$-subalgebra $\cC_\cL := \eta_\cL(\cL)\cC\eta_\cL(\cL)$ of $\cC$.
Moreover, an $\cL$-representation $\rho$ of $\cC$ is uniquely determined by its restriction to  $\eta_\cL(\cL)\cC\eta_\cL(\cL)$
via the relation
\[
\rho(C)=\slim_i\slim_j\rho(\eta_\cL(E_i)C\eta_\cL(E_j))\qquad\hbox{for any approximate identity $(E_j)$ of $\cL$.}
\]
This leads us to a closer analysis of the hereditary subalgebra
$\cC_\cL$ of $\cC$ (cf.\ Proposition~\ref{prop:a.3}).

In the context of assuming (CP1) and (CP4) for $(\cC, \eta_\cA, \eta_\cL)$,
consider the hereditary $C^*$-subalgebra of $M(\cC)$ generated by $\eta_\cL(\cL)\subset M(\cC)$.
It is
\[
M(\cC)_\cL:=\eta_\cL(\cL) M(\cC)\eta_\cL(\cL)=\eta_\cL(\cL) M(\cC)\cap M(\cC)\eta_\cL(\cL)
\]
(cf.~\cite[Cor.~II.5.3.9 and Prop.~II.5.3.2]{Bla06}),
hence it is characterized by
\begin{equation}
\label{HerApprox}
M(\cC)_\cL=\big\{M\in M(\cC)\mid\hbox{
$\eta_\cL(E_j)M \to M $ and $ M\eta_\cL(E_j) \to M$}\big\}
\end{equation}
 for any approximate identity $\{E_j\}$ of $\cL$.
The intersection of $M(\cC)_\cL$ with $\cC$  is again hereditary
by the following lemma.

\begin{lem} \mlabel{lem:a.3b} Let
$\cB_0 \subeq \cB$ be a hereditary subalgebra of the $C^*$-algebra $\cB$.
Then for any $C^*$-subalgebra ${\cal D} \subeq \cB$, the
intersection ${\cal D}\cap \cB_0$ is hereditary in ${\cal D}$.
\end{lem}

\begin{prf} That $\cB_0$ is hereditary means that, for
$B \in \cB$ and $B_0 \in \cB_0$, the relation
$0 \leq B \leq B_0$ implies $B \in \cB_0$, i.e.,
the positive cone in $\cB_0$ is a face of the positive cone
in $\cB$ (\cite[Sect.~3.2]{Mu90}).
For $D \in {\cal D}$ and $D_0 \in \cB_0 \cap {\cal D}$ with
$0 \leq D \leq D_0$ we thus obtain $D \in \cB_0$, so that
$D \in {\cal D} \cap \cB_0$.
\end{prf}

We put $\cC_\cL:=M(\cC)_\cL\cap\cC=\eta_\cL(\cL) \cC\eta_\cL(\cL)$.
Then $\eta_\cL$ restricts to a homomorphism
\[
\eta_\cL^\circ\: \cL \to M(\cC_\cL), \quad \eta_\cL^\circ(L)C :=
\eta_\cL(L)C,\;\; C\in\cC_\cL
\]
which is nondegenerate by the characterization of $M(\cC)_\cL$ in (\ref{HerApprox}).

The condition (CP2), i.e. nondegeneracy of $\eta_\cL \: \cL \to M(\cC)$, is equivalent to
the condition
$\| \eta_\cL(E_j) C-  C \| \to 0$ for all $C\in\cC$, and any approximate identity
$\{E_j\}$ of $\cL$ (cf.  Theorem~\ref{thm:a.2}(iv)).
This  means that the hereditary subalgebra $\cC_\cL$ must be all of $\cC$.
In the case that this does not hold, our task is to restrict  the system in order to obtain
one for which the crossed product host can be defined. It is natural
to try to build  a crossed product host from $\cC_\cL$,
where we already have the nondegenerate action $\eta_\cL^\circ\: \cL \to M(\cC_\cL)$.
However, the action $\eta_\cA \: \cA \to M(\cC)$ need not restrict to $\cC_\cL$, so we define:
\begin{defn} \mlabel{def:al}
Given a triple $(\cC, \eta_\cA, \eta_\cL)$ satisfying (CP1) and (CP4)
and $\cC_\cL$ as  above, let
\[
\cA_\cL:=\big\{A\in \cA\mid \eta_\cA(A)\cC_\cL\subseteq\cC_\cL\quad\hbox{and}\quad \eta_\cA(A^*)\cC_\cL\subseteq\cC_\cL\big\}\,.
\]
Note that this includes the commutant of $\eta_\cL(\cL)$ in $\cA$. Thus  $\eta_\cA$ restricts to a homomorphism
\[
\eta_\cA^\circ\: \cA_\cL \to M(\cC_\cL), \quad \eta_\cA^\circ(A)C :=
\eta_\cA(A)C,\;\; C\in\cC_\cL,\; A\in\cA_\cL.
\]
\end{defn}

To put the algebra $\cA_\cL$ into context, recall the following framework. Let $\cC^{**}$
be the enveloping  $W^*$-algebra
which contains $M(\cC)$ as a subalgebra.
Then the weak limit
\begin{equation}\label{eq:appident}
P := \lim \eta_\cL(E_j)
\end{equation}
exists in $\cC^{**}$ for any approximate identity $\{E_j\}$ of $\cL$,
and defines a projection. From condition (\ref{HerApprox}) we see that
 $\cC_\cL = M(\cC)_\cL\cap\cC=P \cC^{**} P \cap \cC$, hence $P$ is the open
projection of  $\cC_\cL$ in Pedersen's terminology
(cf.~\cite[Prop.~3.11.9 and 3.11.10]{Pe89} and Proposition~\ref{prop:a.3}(iii)). We have:

\begin{prop} \mlabel{propALchar}
Given a $C^*$-action $(\cA, G, \alpha)$
and a triple $(\cC, \eta_\cA, \eta_\cL)$ satisfying
(CP1) and (CP4), then with $\cA_\cL$ and $P$ as defined above;-
\begin{itemize}
\item[\rm(i)]
$\cA_\cL=\big\{A\in \cA\mid [\eta_\cA(A),P]=0\big\}=\big\{A\in \cA\mid \eta_\cA(A)=P\eta_\cA(A)P+(\1-P)\eta_\cA(A)(\1-P)\big\}$
\\
i.e. $\eta_\cA(A)$ consists of diagonal elements w.r.t.\ matrix decomposition of elements of $M(\cC)$ w.r.t.\ $P$ and $\1-P$.
\item[\rm(ii)] $\cA_\cL = \{ A \in \cA \mid
\eta_\cA(A)\eta_\cL(\cL) \subeq \eta_\cL(\cL)\cC\;\;\hbox{and}\;\;
\eta_\cA(A^*)\eta_\cL(\cL) \subeq \eta_\cL(\cL)\cC\}$,
\item[\rm(iii)] For any (hence all) approximate identities $\{E_j\}$ of $\cL$ we have\\
$\cA_\cL =\left\{ A \in \cA \,\big|\,
\big\| \big(\eta_\cL(E_j) -\1\big)\eta_\cA(B) \eta_\cL(L)
 \big\| \to 0 \;\;\mbox{ for} \;\;B = A\;\hbox{and}\;A^*,\;\hbox{and for all}\;\;L \in \cL  \right\}$.
\item[\rm(iv)]
In addition, assume (CP3'). Then $\cA_\cL$ is an  $\alpha_G\hbox{--invariant}$ subalgebra of $\cA$.
Moreover $\cA_\cL$ contains all the elements of $\cA$ which are invariant
w.r.t.\ $\alpha_G$.
\end{itemize}
\end{prop}
\begin{prf}
(i) Since $PC=C$ for all $C\in\cC_\cL$, the definition of $\cA_\cL $ implies
for $A \in \cA_\cL$ the relations
$P\eta_\cA(A)C =\eta_\cA(A)C$  and
$P\eta_\cA(A^*)C =\eta_\cA(A^*)C$.
Thus $(\1-P)\eta_\cA(A)C=0$ for all $C\in\cC_\cL$,
hence by substituting for $C$ an approximate identity for $\cC_\cL$,
we get  $(\1-P)\eta_\cA(A)P=0$. From the analogous relation for the adjoint, we
also get $P\eta_\cA(A)(\1-P)=0$, hence
\[ \eta_\cA(A)=P\eta_\cA(A)P+(\1-P)\eta_\cA(A)(\1-P),
\quad \mbox{ i.e. } \quad [\eta_\cA(A),P]=0.\]
Conversely, if  $[\eta_\cA(A),P]=0$, then $\eta_\cA(A)C=\eta_\cA(A)PCP=P\eta_\cA(A)CP\in P \cC^{**} P \cap \cC=\cC_\cL$
for all $C\in\cC_\cL$, hence $\eta_\cA(A)\cC_\cL\subseteq\cC_\cL$. Likewise $\eta_\cA(A^*)\cC_\cL\subseteq\cC_\cL$,
and hence $A\in\cA_\cL$.

(ii) Let $A \in \cA$ satisfy
$\eta_\cA(A)\eta_\cL(\cL) \subeq \eta_\cL(\cL)\cC$ and
$\eta_\cA(A^*)\eta_\cL(\cL) \subeq \eta_\cL(\cL)\cC$.
Then
\[
\eta_\cA(A)\cC_\cL=\eta_\cA(A)\eta_\cL(\cL) \cC\eta_\cL(\cL)\subeq
\eta_\cL(\cL)\cC\cdot \cC\eta_\cL(\cL)=\cC_\cL
\]
and likewise $\eta_\cA(A^*)\cC_\cL\subeq\cC_\cL$, i.e. $A\in\cA_\cL $.
Conversely, let $A\in\cA_\cL $, i.e.
$\big[\eta_\cA(A),P\big]=0$. Then
\[
\eta_\cA(A)\eta_\cL(\cL)=\eta_\cA(A)P\eta_\cL(\cL)=P\eta_\cA(A)\eta_\cL(\cL)
\subeq\br\eta_\cL(\cL)  \eta_\cA(A)\eta_\cL(\cL).\subeq M(\cC)_\cL\,.
\]
However (CP4) implies that $\eta_\cA(A)\eta_\cL(\cL)\subset\cC$, hence
$\eta_\cA(A)\eta_\cL(\cL)\subset M(\cC)_\cL\cap\cC=\cC_\cL\subeq \eta_\cL(\cL)\cC$.
Likewise $\eta_\cA(A^*)\eta_\cL(\cL)\subeq \eta_\cL(\cL)\cC$.

(iii) Let $A\in\cA_\cL $. Then by (ii) $\eta_\cA(A)\eta_\cL(\cL) \subeq \eta_\cL(\cL)\cC$
and so $\big\| \big(\eta_\cL(E_j) -\1\big)\eta_\cA(A) \eta_\cL(L)
 \big\| \to 0 $
for any approximate identity $\{E_j\}$ of $\cL$, since
$\big\| \big(\eta_\cL(E_j) -\1\big) \eta_\cL(L)
 \big\| \to 0 $. The same is true for $A^*$, hence the condition in (iii) follows.
 Conversely, the condition in (iii) implies immediately that
 $\eta_\cA(A)\eta_\cL(\cL) \subeq \eta_\cL(\cL)\cC\supseteq\eta_\cA(A^*)\eta_\cL(\cL)$,
 using the fact that  $\eta_\cL(\cL)\cC$ is closed. Thus by (ii) we have $A\in\cA_\cL $.

 (iv) By {\rm(CP3')} the universal $\cL$-representation  $\rho_u$ is faithful, hence it
has a faithful extension to $M(\cC)$. Moreover via the covariance condition  ${\rm Ad}(U^{\rho_u}(g))$ defines an automorphism
of $M(\cC)$ which preserves both $\rho_u(\cC)$, $\tilde\rho_u(\eta_\cA(\cA))$ and $\tilde\rho_u(\eta_\cL(\cL))$ for all $g\in G$
(see the proof of Proposition~\ref{LrepsofC}(iv)). Therefore these automorphisms preserve $\rho_u(\cC_\cL)=\tilde\rho_u(\eta_\cL(\cL))\rho_u(\cC)
\tilde\rho_u(\eta_\cL(\cL))$, and hence the idealizer of it in $\tilde\rho_u(M(\cL))$ and hence the intersection of this idealizer
with $\tilde\rho_u(\eta_\cA(\cA))$. By faithfulness of $\tilde\rho_u$, the latter set is $\tilde\rho_u(\eta_\cA(\cA_\cL))$,
and by the covariance condition ${\rm Ad}(U^{\rho_u}(g))$ implements $\alpha_g$ on $\cA\supset\cA_\cL$. It follows by
faithfulness of $\tilde\rho_u$, and the fact that $\cA_\cL$ is defined via $\eta_\cA$ that
$\cA_\cL$ is an  $\alpha_G\hbox{--invariant}$ subalgebra of $\cA$.

For the last statement, notice that if $A\in\cA$ is $\alpha_G\hbox{--invariant}$, then by the covariance condition
$\tilde\rho_u(\eta_\cA(A))$ commutes with  $U^{\rho_u}(G)$. From Proposition~\ref{LrepsofC}(ii)
we get that
$\tilde\rho_u(\eta_\cL(\cL))'= U^{\rho_u}(G)'$, hence that $\tilde\rho_u(\eta_\cA(A))\in \tilde\rho_u(\eta_\cL(\cL))'$.
Since $\tilde\rho_u$ is  faithful, we find that ${\big[\eta_\cA(A)),\eta_\cL(\cL)\big]}=0$, and hence by condition~(iii)
that $A\in\cA_\cL$.
\end{prf}
\noindent
As the commutant of an operator is a von Neumann algebra, it follows from
(i) that $\cA_\cL$ is in fact closed in $\cA$
w.r.t.\ the weak operator topology of the universal representation of $\cC$.
\begin{cor}\mlabel{corAisAL}
With $\alpha \: G \to \Aut(\cA)$ and a triple $(\cC, \eta_\cA, \eta_\cL)$ satisfying
(CP1) and (CP4) as above, we have
\begin{itemize}
\item[\rm(i)]
 $\cA_\cL=\cA$ if and only if $\;\cC_\cL=\cC$
\item[\rm(ii)] If  $(\cC, \eta_\cA, \eta_\cL)$ is
 constructed from the universal covariant $\cL$-representation
of $(\cA, G, \alpha)$, then a full crossed product host exists
if and only if   $\cA_\cL=\cA$.
\end{itemize}
\end{cor}

\begin{prf}
(i) That $\cC_\cL=\cC$ implies $\cA_\cL=\cA$ is trivial,
so we prove the reverse.
Let $\cA_\cL=\cA$. Then,
by Proposition~\ref{propALchar}(ii), we have
$\eta_\cA(A)\eta_\cL(\cL) \subeq \eta_\cL(\cL)\cC$ and
$\eta_\cL(\cL)\eta_\cA(A) \subeq \cC\eta_\cL(\cL)$ for all $A\in\cA$.
This implies that $\eta_\cA(\cA)\eta_\cL(\cL) \subeq
\eta_\cL(\cL) M(\cC)\cap M(\cC)\eta_\cL(\cL)\cap\cC=M(\cC)_\cL\cap\cC=\cC_\cL\subeq\cC$
and so as the first term generates $\cC$ we get that $\cC=\cC_\cL$.

(ii) From  Proposition~\ref{propALchar}(iii),
we note that the existence condition in Theorem~\ref{thm:exist} for a full crossed product host
becomes $\cA = \cA_\cL$ for the triple $(\cC, \eta_\cA, \eta_\cL)$.
\end{prf}

For the triple $(\cC_\cL, \eta_\cA^\circ, \eta_\cL^\circ)$,
 the action $\eta_\cL^\circ\: \cL \to M(\cC_\cL)$ is nondegenerate (i.e. we have (CP2)),
but we lost (CP4) by the restriction  from $\cA$
to $\cA_\cL$, i.e. it need not be true that $\cC_\cL$ is generated by
$\eta_\cL(\cL)\eta_\cA(\cA_\cL)$. This leads us to define
\[
\cC_\cL^\vee:=C^*\big(\eta_\cL(\cL)\eta_\cA(\cA_\cL)\big)\subset\cC\,.
\]
Note that as $\eta_\cL(\cL)=P\eta_\cL(\cL)P$, and $\eta_\cA(\cA_\cL)$ commutes with $P$ we have
$\eta_\cL(\cL)\eta_\cA(\cA_\cL)\subset P \cC^{**} P \cap \cC=\cC_\cL$, hence
$\cC_\cL^\vee\subseteq \cC_\cL$. By the analogous proof to that of Theorem~\ref{thm:5.1}(b),
we see that $\eta_\cL(\cL)$ is in the idealizer of $\cC_\cL^\vee$
in $\cC_\cL$, and by condition
(\ref{HerApprox}) the multiplier action which this produces on $\cC_\cL^\vee$
is in fact nondegenerate. Thus $\eta_\cL$
restricts to the nondegenerate action
\[
\eta_\cL^\vee\: \cL \to M(\cC_\cL^\vee), \quad \eta_\cL^\vee(L)C :=
\eta_\cL(L)C \quad \mbox{ for } \quad C\in\cC_\cL^\vee.
\]
Likewise, by the analogous proof to that of Theorem~\ref{thm:5.1}(b), we also restrict the action $\eta_\cA$:
\[
\eta_\cA^\vee\: \cA_\cL \to M(\cC_\cL^\vee), \quad \eta_\cA^\vee(A)C :=
\eta_\cA(A)C,\;\; C\in\cC_\cL^\vee,\; A\in\cA_\cL.
\]
It follows that the triple  $(\cC_\cL^\vee, \eta_\cA^\vee, \eta_\cL^\vee)$
satisfies (CP1), (CP2) and (CP4) for the restricted action \break $\alpha:G\to\Aut\cA_\cL$.
\begin{prop} \mlabel{restXprod}
Given a group homomorphism $\alpha \: G \to \Aut(\cA)$ and a triple $(\cC, \eta_\cA, \eta_\cL)$ satisfying
(CP1), (CP3') and (CP4), then
the triple  $(\cC_\cL^\vee, \eta_\cA^\vee, \eta_\cL^\vee)$
for the restricted action \break $\alpha:G\to\Aut\cA_\cL$
satisfies (CP1)-(CP4), hence  $(\cC_\cL^\vee, \eta_\cA^\vee, \eta_\cL^\vee)$ is a
crossed product host (not necessarily full) for ${(\alpha\!\restriction \!\cA_\cL,\,\cL)}$.
\end{prop}

This will be called the {\it restricted crossed product host}.

\begin{prf}
We only have to prove (CP3), as we already know from Proposition~\ref{propALchar}(iv) that $\cA_\cL$
is an $\alpha_G\hbox{--invariant}$ subalgebra, and from the preamble that  (CP1), (CP2) and (CP4)
are satisfied.
As  (CP2) holds, Theorem~\ref{thm:a.2}(iii) applies, hence the multiplier extension $\tilde\eta_\cL^\vee \: M(\cL)\to M(\cC_\cL^\vee)$
is defined, and the covariance requirement (CP3) makes sense.
Hence we have to prove that
\[ \tilde\eta_\cL^\vee(\eta(g))\eta_\cA^\vee(A)\tilde\eta_\cL^\vee(\eta(g))^*
= \eta_\cA^\vee(\alpha_g (A))\qquad\hbox{for }\qquad g\in G,\; A\in\cA_\cL.\]
Abbreviate the notation to $\eta_G^\vee:= \tilde\eta_\cL^\vee\circ\eta:G\to M(\cC_\cL^\vee)$.  Starting now from (CP3'), i.e.
\[
U^{\rho_u}(g)\cdot(\tilde\rho_u\circ\eta_\cA(A))
\cdot U^{\rho_u}(g)^*=\tilde\rho_u\circ\eta_\cA\big(\alpha_g(A)\big)
\qquad\hbox{ for }\qquad g\in G,\; A\in\cA,
\]
recall that
$U^{\rho_u}:G\to \U(\cH_u)$
is  just obtained from the strictly continuous extension of
the representation $\tilde\rho_u\circ\eta_\cL$ from $\cL$ to $M(\cL)$.
Let $P_\cL\in \cB(\cH_u)$ be the projection onto the essential subspace of $\tilde\rho_u\circ\eta_\cA(\cA_\cL)$,
i.e. $P_\cL=\slim\limits_j\tilde\rho_u\circ\eta_\cA(E_j)$ for an approximate identity $(E_j)$ of $\cA_\cL$.
By (CP3') we see that all $U^{\rho_u}(g)$ commute with $P_\cL$, hence $U^{\rho_u}(G)$ restricts to
$P_\cL\cH_u$. In fact, by Proposition~\ref{LrepsofC}(ii) we get that
$P_\cL\in\tilde\rho_u(\eta_\cL(\cL))'= U^{\rho_u}(G)'$, and so it follows that $P_\cL\cH_u$
is also the essential subspace of $\rho_u(\cC_\cL^\vee)=C^*\big(\tilde\rho_u(\eta_\cL(\cL)\eta_\cA(\cA_\cL))\big)$.

Let $\rho_u^\vee:\cC_\cL^\vee\to B(P_\cL\cH_u)$ be the restriction of $\rho_u(\cC_\cL^\vee)$ to its essential
subspace. Clearly $\rho_u^\vee$ is faithful, and it satisfies (CP3') (restricted to $P_\cL\cH_u$), and in fact
 $U^{\rho_u}(G)$ restricted to $P_\cL\cH_u$ is  again the strictly continuous extension of
the representation $\tilde\rho_u^\vee\circ\eta_\cL^\vee$ from $\cL$ to $M(\cL)$.
(Recall that $\cC^\vee_\cL\subset\cC$, and that $\eta_\cL^\vee(\cL)$ is just the restriction of
$\eta_\cL(\cL)\subset M(\cC)$ to $\cC^\vee_\cL$).
On the other hand,
$\tilde\rho_u^\vee \circ \eta_G^\vee:=\tilde\rho_u^\vee \circ (\tilde\eta_\cL^\vee\circ\eta)$ is also a strictly continuous extension of
the representation $\tilde\rho_u^\vee\circ\eta_\cL^\vee$ from $\cL$ to $M(\cL)$, hence they must be equal, i.e.
${U^{\rho_u}(g)\restriction P_\cL\cH_u}
=\tilde\rho_u^\vee \circ \eta_G^\vee(g)$. Thus (CP3') becomes
\[
\tilde\rho_u^\vee \big( \eta_G^\vee(g)\big)\cdot\tilde\rho_u^\vee\big(\eta_\cA^\vee(A)\big)\cdot \tilde\rho_u^\vee \big( \eta_G^\vee(g)\big)^*
=\tilde\rho_u^\vee\big(\eta_\cA^\vee\big(\alpha_g(A)\big)\big)
\qquad\hbox{for all}\qquad g\in G,\; A\in\cA_\cL,
\]
and as $\tilde\rho_u^\vee$ is faithful on $M(\cC_\cL^\vee)$, this implies (CP3).
\end{prf}

If  $(\cC, \eta_\cA, \eta_\cL)$ is given explicitly by some covariant $\cL$-representation
$(\pi,U)\in{\rm Rep}(\alpha,\cH)$, then we conclude that the restriction of it to
$\cA_\cL$ is a cross representation of ${(\alpha\!\restriction \!\cA_\cL,\,\cL)}$.
The natural question is when the crossed product
host for ${(\alpha\!\restriction \!\cA_\cL,\,\cL)}$ is full, and we will
analyze this question below in specific contexts.

\begin{prop} \mlabel{FullrestX}
Given a $C^*$-action $(\cA,G,\alpha)$,
let $(\cC, \eta_\cA, \eta_\cL)$  be constructed from
the universal  covariant $\cL$-representation
$(\pi_u,U_u)\in{\rm Rep}_\cL(\alpha,\cH_u)$ where
$\al C.:=C^*\big(\pi_u(\cA) U_{u,\cL}(\cL)\big),$ and $\eta_\cA$ and $ \eta_\cL$
are as in Theorem~\ref{thm:5.1}(b).
Suppose that either of the following conditions hold:
\begin{itemize}
\item[(i)] $\fS_\cL$ separates $\al C.^\vee:=C^*\big(\pi_u^\vee(\cA_\cL) U_{u,\cL}^\vee(\cL)\big)$
where $(\pi_u^\vee,U_u^\vee)\in{\rm Rep}_\cL(\alpha,\cH_u^\vee)$ denotes
the universal  covariant $\cL$-representation of
$(\cA_\cL, G, \alpha)$.
\item[(ii)]  $\fS_\cL\restriction\big(\cA_\cL \rtimes_\alpha G_d\big)=\fS_\cL^\vee$,
where $\fS_\cL^\vee$
denotes the set of those states
$\omega$ on $\cA_\cL\rtimes_{\alpha} G_d$ which produce
covariant $\cL$-representations
$(\pi_\omega,U_\omega)\in{\rm Rep}_\cL(\alpha\!\restriction \!\cA_\cL,\cH_\omega)$.
\end{itemize}
Then the restricted crossed product host $(\cC_\cL^\vee, \eta_\cA^\vee, \eta_\cL^\vee)$ for
 ${(\alpha\!\restriction \!\cA_\cL,\,\cL)}$ is full.
\end{prop}

\begin{prf}
If $\al C.^\vee\cong\al C._\cL^\vee$ then as  $(\cC_\cL^\vee, \eta_\cA^\vee, \eta_\cL^\vee)$ is a
crossed product host, it follows that $\al C.^\vee$ is a
crossed product host (Proposition~\ref{restXprod}), i.e. the
universal  covariant $\cL$-representation of $(\cA, G, \alpha)$ is a
cross representation,
and thus $\al C.^\vee$ is full by Theorem~\ref{thm:exist}. So it suffices to show that each of the two conditions imply
that $\al C.^\vee\cong\al C._\cL^\vee$.

Recall that $\al C._\cL^\vee=C^*\big(\pi_u(\cA_\cL) U_{u,\cL}(\cL)\big)\subseteq\cC$ and
 $\al C.^\vee=C^*\big(\pi_u^\vee(\cA_\cL) U_{u,\cL}^\vee(\cL)\big)$ where
 ${(\pi_u,U_u)}$ is constructed from the GNS-representations of $\fS_\cL$ on $\cA \rtimes_\alpha G_d$
 and ${(\pi_u^\vee,U_u^\vee)}$ is constructed from the GNS-representations of $\fS_\cL^\vee$ on $\cA_\cL \rtimes_\alpha G_d
 \subset\cA \rtimes_\alpha G_d$. Now restrictions of $\cL$-representations of $\cA$
 to $\cA_\cL$ are still $\cL$-representations, so as
 subrepresentations of $\cL$-representations are $\cL$-representations,
 we have that
 $\fS_\cL\restriction\big(\cA_\cL \rtimes_\alpha G_d\big)\subseteq\fS_\cL^\vee$.
 Now the GNS-representation of $\omega\in\fS_\cL$ on $\cA \rtimes_\alpha G_d$ need not be cyclic when
 we restrict it to $\cA_\cL \rtimes_\alpha G_d$, but it decomposes into cyclic components. Each of these
 is the GNS-representation of a state in $\fS_\cL^\vee$ as subrepresentations of $\cL$-representations are $\cL$-representations.
Thus $(\pi_u, U_u)$ is a direct sum of subrepresentations,
each of which occurs in $(\pi_u^\vee, U_u^\vee)$.
Hence there is a $*$-homomorphism
$\Phi:\al C.^\vee\to\al C._\cL^\vee$.

If we assume (ii), then the same subrepresentations occur in both $\al C._\cL^\vee$ and $\al C.^\vee$, hence $\Phi$ is an isomorphism.
The condition (i) makes sense, since we have that
 $\fS_\cL\restriction\big(\cA_\cL \rtimes_\alpha G_d\big)\subseteq\fS_\cL^\vee$, so any
 any $\omega\in\fS_\cL$ can be realized as a vector state on $\al C.^\vee$,
 and via $\Phi$ it coincides with the lifting of this state from $\al C._\cL^\vee$ by $\Phi$.
Thus $\omega(\ker\Phi)=0$ and so if we assume condition (i) we must have $\ker\Phi=\{0\}$ i.e. $\al C.^\vee\cong\al C._\cL^\vee$.
\end{prf}
It is also natural to construct the following subalgebra:

\begin{prop}
\mlabel{NContAlg} Let $(\cA,G, \alpha)$ be a $C^*$-action,
and suppose that $(\cC, \eta_\cA, \eta_\cL)$ satisfies
(CP1), (CP3') and (CP4), and that the
multiplier action $\eta \: G \to \U(M({\cal L}))$ is
strictly continuous. Let $P \in \cC^{**}$ be the projection from
\eqref{eq:appident}. Then the subspace
 \[
 \cA^{(\cL)}_0:=\{A\in\cA\mid\eta_\cA(A)\in\eta_\cL(\cL) M(\cC)\eta_\cL(\cL)\}
 =\{A\in\cA\mid\eta_\cA(A)=P\eta_\cA(A)P\}
 \] is a closed two-sided $\alpha_G$-invariant ideal of
 $\cA_\cL$.   Moreover, the action $\alpha^\eta\: G \to \Aut(\eta_\cA(\cA^{(\cL)}_0))$
 defined by
$\alpha^\eta_g(\eta_\cA(A)):=\eta_\cA(\alpha_g(A))$, $A\in \cA^{(\cL)}_0$ is strongly continuous.
\end{prop}

\begin{prf}
Any $A\in \cA^{(\cL)}_0$ must satisfy
$\eta_\cA(A)=P\eta_\cA(A)P$,
hence by Proposition~\ref{propALchar}(i) it is clear that it is in $\cA_\cL$. Conversely, if an $A\in\cA_\cL$
satisfies $\eta_\cA(A)=P\eta_\cA(A)P$ then it must be in $\cA^{(\cL)}_0$ because
$\eta_\cA(\cA^{(\cL)}_0)= M(\cC)_\cL\cap\eta_\cA(\cA)=P \cC^{**} P \cap \eta_\cA(\cA)$, thus $\cA^{(\cL)}_0$ is precisely
the elements of $\cA_\cL$ satisfying $\eta_\cA(A)=P\eta_\cA(A)P$. Since $\cA_\cL$ is in the commutant of $P$, it is clear
that $\cA_\cL\cdot\cA^{(\cL)}_0=\cA^{(\cL)}_0$ hence that  $\cA^{(\cL)}_0$ is a closed two-sided ideal of
 $\cA_\cL$. That $\cA^{(\cL)}_0$ is $\alpha_G\hbox{--invariant}$ follows immediately from
 the fact that both  $\eta_\cA(\cA_\cL)$ and $M(\cC)_\cL$ are
$\alpha_G$-invariant subalgebras of $M(\cC)$
 w.r.t.\ the extended action (cf. Proposition~\ref{propALchar}(iv)
and Remark~\ref{actionMC}). The strong continuity of $\alpha^\eta$ also follows immediately from the strong
continuity on $M(\cC)_\cL$ (cf. Remark~\ref{actionMC}).
\end{prf}
Whilst the ideal $\cA^{(\cL)}_0$ has very nice properties, below we will
see in Example~\ref{QM-Weyl} that it can be zero, even when $\cA_\cL=\cA$.

 \begin{ex}
 \mlabel{ExHHop}
 We want to obtain the restricted crossed product  $(\cC_\cL^\vee, \eta_\cA^\vee, \eta_\cL^\vee)$
 explicitly in an example.
We continue Example~\ref{ex:noncrossrep}, and recall the details.
Let $\cH$ be an infinite-dimensional separable Hilbert space,
$\cA := \cB(\cH)$, $G := \R$, $\cL = C^*(\R)$,
$H$ be an unbounded selfadjoint operator,
$U_t := e^{itH}$ and $\alpha_t(A) := U_t A U_t^*$.
Then $\alpha$ is not strongly continuous (Proposition~\ref{prop:c.4}). Moreover
$(\pi,U)\in\Rep_\cL(\alpha, \cH)$ where $\pi$ is the identical representation
$\pi(A) = A$ of $\cA$. We have $U_\cL(\cL)=\{h(H)\,\mid\,h\in C_0(\R)\}= C^*((i\1-H)^{-1}).$
We showed that
$(\pi,U) \in\Rep_\cL^\times(\alpha, \cH)$ if and only if $(i\1-H)^{-1}\in\cK(\cH)$.

Choose an $H$ so that $(\pi,U)\not\in\Rep_\cL^\times(\alpha, \cH)$, i.e. one point
in the spectrum of $H$ is not in its essential spectrum, hence $U_\cL(\cL)$ contains
noncompact elements.
By definition
$\cC=C^*\big(\pi(\cA) U_\cL(\cL) \big)\subset \cB(\cH)=\pi(\cA)=\eta_\cA(\cA)$
and as $\eta_\cA(\cA)\subseteq M(\cC)$, we see that $\cC$ is a closed two-sided ideal
of $\cB(\cH)$ and as it contains noncompact elements, we must have
$\cC=B(\cH)=\pi(\cA)=\eta_\cA(\cA)$. Now $\cC_\cL:=\eta_\cL(\cL) \cC\eta_\cL(\cL)
=U_\cL(\cL)\cB(\cH)U_\cL(\cL)$
is an hereditary subalgebra of the von Neumann algebra $\cB(\cH)=\cC$. (As
$(\pi,U)\in\Rep_\cL(\alpha, \cH)$,  hence $U_\cL(\cL)$ is nondegenerate, the strong closure
of $\cC_\cL$ is $\cB(\cH)$). Notice that as $\cC_\cL$ is hereditary in $\cB(\cH)$, it contains all
projections which it dominates, e.g.\ it contains the spectral projections $P[a,b]$ of $H$
(in fact it is the closed linear span of its projections by
\cite[Cor.~4.1.14]{Mu90}).

By Definition~\ref{def:al}
$\cA_\cL$ is the idealizer of $\cC_\cL$ in $\cA$,
and we know it is proper by
$(\pi,U)\not\in\Rep_\cL^\times(\alpha, \cH)$. In this case  we have $\cC_\cL\subseteq\cA_\cL$
as all algebras are contained in $\cA$, and it is clear that $\cC_\cL$ is an ideal of $\cA_\cL$.
Moreover (cf. Proposition~\ref{NContAlg}), here we have $\cA^{(\cL)}_0=\cC_\cL$.
We want to express $\cA_\cL$ directly
in terms of the properties of $H$.
From Proposition~\ref{contSpace}(iii) below, we have that
\[ \cA_\cL = \big\{ A \in \cA \;\big|\;
\lim_{t \to 0} \eta_\cA(\alpha_t(B))\eta_\cL(L) = \eta_\cA(B)\eta_\cL(L)
\mbox{ for all } L \in \cL\mbox{ and }B \in\{ A,\, A^*\}\big\}.\]
If $\rho\in{\rm Rep}_\cL(\cC,\cH)$ is the (faithful) defining representation of $\cC$,
obtained from $(\pi,U)$, then by applying it to the last
characterization of $\cA_\cL$ we see:
\begin{eqnarray*}
\cA_\cL &=& \big\{ A \in \cA \;\big|\;
\lim_{t\to 0} U_t\pi(B)U_{-t}U_\cL(L) = \pi(B)U_\cL(L)
\mbox{ for all } L \in \cL\mbox{ and }B \in\{ A,\, A^*\}\big\}
\\[1mm]
&=& \big\{ A \in \cA\;\big|\;
\lim_{t\to 0} \big(U_t-\1\big)\pi(B)U_\cL(L) = 0
\mbox{ for all } L \in \cL\mbox{ and }B \in\{ A,\, A^*\}\big\}
\\[1mm]
&=& \big\{ A \in \cA \;\big|\;
\lim_{t\to \infty} \big(P[-t,t]-\1\big)\pi(B)U_\cL(L) = 0
\mbox{ for all } L \in \cL\mbox{ and }B \in\{ A,\, A^*\} \big\}
\end{eqnarray*}
where in the second equality we used the fact that $t\to U_{-t}U_\cL(L)$ is norm continuous
and that $(U_t)$ is bounded, and in the last equality we used Lemma~\ref{lem:a.1}(ii)(c).

To obtain the restricted crossed product  $(\cC_\cL^\vee, \eta_\cA^\vee, \eta_\cL^\vee)$
for  ${(\alpha\!\restriction \!\cA_\cL,\,\cL)}$, note that
\[
\cC_\cL^\vee:=C^*\big(\eta_\cL(\cL)\eta_\cA(\cA_\cL)\big)\subseteq \cC_\cL=\eta_\cL(\cL) \cC\eta_\cL(\cL)
=U_\cL(\cL)\cB(\cH)U_\cL(\cL)\unlhd\cA_\cL
\]
where the last relation indicates inclusion as a closed two-sided ideal. Since the $C^*$-algebra
$\cC_\cL^\vee$ contains
$\eta_\cL(\cL)\cC_\cL=\cC_\cL$ it follows that $\cC_\cL^\vee=\cC_\cL=U_\cL(\cL)\cB(\cH)U_\cL(\cL)$.
\end{ex}
\begin{rem}
In a physics context, the group homomorphism $\alpha \: G \to \Aut(\cA)$ represents
a physical transformation or symmetry group acting on the field algebra $\cA$, the choice of
host algebra $(\cL,\eta)$  for $G$ is a constraint specifying the type of unitary representations
allowed to implement $\alpha$, and a given $(\pi,U)\in{\rm Rep}_\cL(\alpha,\cH)$ is a
representation with an important physical interpretation which the system should preserve.
In this case $\cA_\cL$ is the algebra of observables which are compatible with these constraints,
so physical observables are required to be contained in it.
Then by Proposition~\ref{restXprod}, the restricted system ${(\alpha\!\restriction \!\cA_\cL,\,\cL)}$
has a crossed product host, $(\cC_\cL^\vee, \eta_\cA^\vee, \eta_\cL^\vee)$, i.e.
$(\pi\!\restriction \!\cA_\cL,U)$ is a cross representation for it.
\end{rem}

\section{Crossed products for discontinuous actions}
\label{CPH-discont}

In this section we consider the special case where
$G$ is a locally compact group, the host algebra is $\cL = C^*(G)$,
and the action $\alpha \: G \to \Aut(\cA)$  need not be strongly continuous.
This is an important case for physics, and has already been analyzed with different tools
by Borchers~\cite{Bo96, Bo84}.

For the standard case, where  the action $\alpha$ is strongly continuous, the
usual crossed product algebra $\cA \rtimes_\alpha G$ is a crossed
product host (Example~\ref{ex:1.1}).
For the case of discontinuous actions, we have already
seen in Examples~\ref{FullCPH-discontAlph},
and \ref{nonfullMax} that crossed product hosts can exist,
and we will see more examples below. Although the multiplier action
of $G$ on any crossed product host  $\cC$
is strictly continuous (cf.\ Remark~\ref{rem:2.5}(a)),
the fact that $\eta_\cA(\cA) \subeq M(\cC)$
need not be contained in $\cC$ leaves room for
discontinuity for the action of $G$ on $\cA$.
We will see below in Corollary~\ref{corContX} that a weaker form of
continuity for $\alpha$ is still required for the existence of a
crossed product host.

The following lemma will be needed below.
\begin{lem} \mlabel{lem:prod}
Let $\cF$ be a unital $C^*$-algebra, $G$ be a topological group,
$\eta \:  G \to \U(\cF)$ be a group homomorphism and
$\alpha_g(F) := \eta(g)F\eta(g)^*$ for $F \in \cF$.
Let $ B \in \cF$ satisfy
$\lim\limits_{g \to \1} \eta(g)B = B$. Then
for an $ A \in \cF$ we have
 $\lim\limits_{g \to \1} \eta(g)AB = AB\;$
if and only if
$\;\lim\limits_{g \to \1} \alpha_g(A)B = AB$.
\end{lem}

\begin{prf} If $\lim\limits_{g \to \1}  \alpha_g(A)B = AB$,
then
\[ \eta(g)AB - AB = \alpha_g(A)B - AB + \eta(g)A(B - \eta(g^{-1})B) \to 0\]
follows from $\|\eta(g)A\| \leq \|A\|$.
Conversely, if  $\lim\limits_{g \to \1} \eta(g)AB= AB$, then we similarly obtain
\[  \alpha_g(A)B - AB  = \eta(g)AB - AB + \eta(g)A(\eta(g^{-1})B-B) \to 0.
\qedhere\]
\end{prf}

\noindent We first want to characterize $\cA_\cL$ directly in terms of the action $\alpha \: G \to \Aut(\cA)$.

\begin{prop}
\mlabel{contSpace}
Let  $(\cL,\eta)$ be a host algebra for a topological group $G$ such that
the multiplier action $\eta \: G \to \U(M({\cal L}))$
is strictly continuous
and let $(\cA, G, \alpha)$
be a $C^*$-action for which we have a
triple $(\cC, \eta_\cA, \eta_\cL)$ satisfying
(CP1), (CP3') and (CP4). We define $\cA_\cL$ as above.
Then the following assertions hold:
\begin{description}
\item[\rm(i)] The subspace $\eta_\cL(\cL)\cC$ is contained in the closed
right ideal of $\cC$:
\[ \cC_c^L := \{ C \in \cC \mid \lim_{g \to \1} U^{\rho_u}(g)\cdot\rho_u(C)= \rho_u(C)\}.\]
\item[\rm(ii)] For $A \in \cA$ we have
$\eta_\cA(A) \eta_\cL(\cL) \subeq \cC_c^L$
if and only if for each $L \in \cL$ the map
\[ G \to \cC, \quad g \mapsto \eta_\cA(\alpha_g(A))\eta_\cL(L)\]
is continuous.
\item[\rm(iii)] In addition, let $G$ be locally compact and
$\cL = C^*(G)$.
Then $\eta_\cL(\cL) \cC = \cC_c^L$ and
\[ \cA_\cL = \big\{ A \in \cA \mid
\lim_{g \to \1} \eta_\cA(\alpha_g(B)) \eta_\cL(L) = \eta_\cA(B) \eta_\cL(L)
\mbox{ for all } L \in \cL\mbox{ and }B \in\{ A,\, A^*\}\big\}.\]
\end{description}
\end{prop}
\begin{prf} (i) follows
from the strict continuity of the multiplier action $\eta \: G \to \U(M({\cal L}))$ via
\[
U^{\rho_u}(g)\cdot\rho_u(\eta_\cL(L)C)
=U^{\rho_u}(g)\cdot(\tilde\rho_u\circ\eta_\cL)(L)\cdot\rho_u(C)
=(\tilde\rho_u\circ\eta_\cL)(\eta(g)L)\cdot\rho_u(C)
\quad\hbox{for}\quad L\in\cL,\,C\in\cC.
\]
(ii) If $\eta_\cA(A) \eta_\cL(L) \in \cC_c^L$, then the map
\[ G \to \cC, \quad g \mapsto  U^{\rho_u}(g)\cdot
\rho_u\big( \eta_\cA(A)\eta_\cL(L)\big)
= (\tilde\rho_u\circ\eta_\cA)(\alpha_g(A))\cdot
(\tilde\rho_u\circ\eta_\cL)\big(\eta(g)L\big) \]
is continuous at $\1$. Since
$g \mapsto \eta_\cL(\eta(g)L)$ is continuous at $\1$, it follows that the
map
\[ G \to \cC, \quad g \mapsto \eta_\cA(\alpha_g(A))\eta_\cL(L) \]
also is continuous at $\1$ using the fact that $\tilde\rho_u$ is faithful (cf.\ Lemma~\ref{lem:prod}).

If, conversely, the latter map is continuous at $\1$, then by
Lemma~\ref{lem:prod} we get that $\eta_\cA(A)\eta_\cL(L) \in \cC_c^L$.

(iii) From Proposition~\ref{LrepsofC}(iii) and (iv) we get the covariant
pair ${(\tilde\rho_u\circ\eta_\cA,\,U^{\rho_u})}\in\Rep_\cL(\alpha, \cH_u)$ hence
that $\rho_u(\cC)=C^*\big((\tilde\rho_u\circ\eta_\cA)(\cA)\cdot
(\tilde\rho_u\circ\eta_\cL)(\cL)\big)$.
In view of (i), Lemma~\ref{lem:a.1}(b)  entails that $\cC_c^L = \eta_\cL(\cL) \cC$.
The assertion on $\cA_\cL$ now follows by combining this
with (ii), using Proposition~\ref{propALchar}(ii).
\end{prf}

\begin{rem}
 Note that in the case of Proposition~\ref{contSpace}(iii),
$\cA_\cL$ contains the well-known subalgebra
\[
\cA_c:= \{ A \in \cA \mid  \lim_{g \to \1} \alpha_g(A) = A\},
\]
though $\cA_\cL$ can be strictly larger. In the case that $\eta_\cA$ is faithful, the containment
$\cA^{(\cL)}_0\subset\cA_c$ is usually proper, e.g.\ if $\cA$ is unital, then the identity is in
$\cA_c$, but not in $\cA^{(\cL)}_0$.
\end{rem}
\begin{cor}\mlabel{corContX}
Let $G$ be locally compact, $\cL = C^*(G)$, and
$(\cA, G, \alpha)$ be a $C^*$-action.  If  $(\cC, \eta_\cA, \eta_\cL)$ is
 constructed from the universal covariant $\cL$-representation
of $(\cA, G, \alpha)$, then the following are equivalent:
\begin{description}
\item[\rm(i)]  A full crossed product host exists.
\item[\rm(ii)] $\lim_{g \to \1} \eta_\cA(\alpha_g(A)) \eta_\cL(L) = \eta_\cA(A) \eta_\cL(L)$
for $A\in\cA$ and $L \in \cL$.
\item[\rm(iii)]  The conjugation action of $G$ on $\cC$ is strongly continuous.
\end{description}

The conditions imply that the maps $G \to M(\cC),
g\to\eta_\cA(\alpha_g(A))$ are continuous w.r.t.\ the strict topology of $\cC$
for all $A\in\cA$.
If $\cA$ is unital, then (i)-(iii) are  equivalent to
\begin{description}
  \item[\rm(iv)] For every $A \in \cA$, the map
$G \to M(\cC), g \mapsto \eta_\cA(\alpha_g(A))$ is strictly continuous.
\end{description}
\end{cor}

\begin{prf} By Proposition~\ref{contSpace}(iii), condition (ii)
means that  $\cA=\cA_\cL$, and by
Corollary~\ref{corAisAL}(ii), in the current context, this is equivalent to
the existence of a full crossed product host, i.e., (i) and (ii) are equivalent.

Since multiplier action of $G$ on $\cL$ is strictly continuous,
the same holds for the conjugation action of $G$ on $\cL$, which is a linear action
by isometries. Therefore the limit of
$\eta_\cA(\alpha_g(A))\eta_\cL(L)$ for $g \to \1$ exists if and only if
the limit of $\eta_\cA(\alpha_g(A))\eta_\cL(\eta(g)L\eta(g)^{-1})$ for $g \to \1$ exists,
and in this case the two limits coincide. Since $\cC$ is generated by
$\eta_\cA(\cA) \eta_\cL(\cL)$, the latter condition means that the
conjugation action of $G$ on $\cC$ is strongly continuous. Therefore
(ii) and (iii) are also equivalent.

Assume (i). Then
$\cC=\cC_\cL=\eta_\cL(\cL) \cC\eta_\cL(\cL)$ by Corollary~\ref{corAisAL},
and so $\lim\limits_{g \to \1} \eta_\cA(\alpha_g(A)) C = \eta_\cA(A) C$ for all $A\in\cA$
and $C\in\cC$ by (ii), i.e. (iv) holds.
Conversely, if $\cA$ is unital, then by Proposition~\ref{LrepsofC}(iv) we have that $\eta_\cA(\1)=\1$, and hence
$\eta_\cL(\cL)\subset\eta_\cL(\cL)\eta_\cA(\cA)\subset\cC$, so any net $(B_\lambda)\subset M(\cC)$ which
converges strictly to $B \in M(\cC)$, will also satisfy $\lim\limits_\lambda B_\lambda\eta_\cL(L)=B\eta_\cL(L)$
for all $L\in\cL$. This gives the sufficiency of (iv) by setting $B_\lambda=\eta_\cA(\alpha_{g_\lambda}(A))$
where $(g_\lambda)\subset G$ is a net converging to $\1$.
\end{prf}

\begin{rem}
For the construction of the crossed product in the conventional sense, one requires strong continuity of the action,
i.e. $\lim\limits_{g \to \1} \alpha_g(A)=A$ for all $A\in\cA$. Here we see that
to obtain a full crossed product host,
we simply replace this by condition~(ii) in Corollary~\ref{corContX}.
Or, for unital $\cA$, we replace the strong continuity
of $\alpha$ on $\cA$ by the requirement of the continuity of the
maps $G \to M(\cC), g \mapsto \alpha_g(A)$ w.r.t.\ the strict topology.
\end{rem}
In the case of a one-parameter group, i.e. $G=\R$, $\cL=C^*(\R)$, we can obtain more
information from Proposition~\ref{contSpace}.

\begin{cor}\mlabel{corRCont}
For  $G=\R$ and  $\cL=C^*(\R)$, let $(\cA, \R, \alpha)$
be a $C^*$-action for which we have a
triple $(\cC, \eta_\cA, \eta_\cL)$ satisfying
(CP1), (CP3') and (CP4).  For a Hilbert space $\cH$, consider
the injection (cf. Proposition~\ref{LrepsofC}(iv)):
 $$ \eta^*_\times \: \Rep_\cL(\cC,{\cal H}) \to \Rep_\cL(\alpha, \cH),  \quad\hbox{given by}\quad
\eta^*_\times(\rho):=\big(\tilde\rho \circ \eta_\cA, U^\rho\big)=: \big(\pi^\rho, U^\rho\big).$$
Let $\rho\in \Rep_\cL(\cC,{\cal H})$ be faithful, and denote the
spectral measure of $U^\rho$, resp., its infinitesimal generator
by $P_\rho$. Then an $A \in \cA$ is in $\cA_\cL$ if and only if
\[
 \lim_{t \to \infty} P_\rho([-t,t])\pi^\rho(B)P_\rho([-s,s]) = \pi^\rho(B)P_\rho([-s,s])
\mbox{ for all } s \in \R_+\mbox{ and }B \in\{ A,\, A^*\}.\]
\end{cor}
\begin{prf}
As $\rho\in \Rep_\cL(\cC,{\cal H})$ is faithful, its extension $\tilde\rho$ to $M(\cC)$ is faithful.
As we have
 \[
 \cA_\cL = \{ A \in \cA \mid
\eta_\cA(B)\eta_\cL(\cL) \subeq \eta_\cL(\cL)\cC\;\mbox{ for }B \in\{ A,\, A^*\}\}
\]
(cf. Proposition~\ref{propALchar}(ii)), and the condition takes place in $M(\cC)$, we get that
 \[
 \cA_\cL = \{ A \in \cA \mid
\pi^\rho(B)U^\rho_\cL(\cL) \subeq U^\rho_\cL(\cL)\rho(\cC)\;\mbox{ for }B \in\{ A,\, A^*\}\}.
\]
 Since $\rho$ is an $\cL$-representation of $\cC$, we have
 $\pi^\rho(B)U^\rho_\cL(\cL) \subeq U^\rho_\cL(\cL)\rho(\cC) \subeq U^\rho_\cL(\cL)\cB(\cH)$.
 On the other hand,
if $\pi^\rho(B)U^\rho_\cL(\cL) \subeq U^\rho_\cL(\cL)\cB(\cH)$, then as the left hand side
 is in $\rho(\cC)$, we get $\pi^\rho(B)U^\rho_\cL(\cL) \subeq U^\rho_\cL(\cL)\cB(\cH)\cap\rho(\cC)
 =\rho(\eta_\cL(\cL)\cC)$.
Thus $\pi^\rho(B)U^\rho_\cL(\cL) \subeq U^\rho_\cL(\cL)\rho(\cC)$ if and only if
 $\pi^\rho(B)U^\rho_\cL(\cL) \subeq U^\rho_\cL(\cL)\cB(\cH)$.

In Lemma~\ref{lem:a.1}(ii)(c),  it is shown that, for $D\in \cB(\cH)$, the condition
$D U^\rho_\cL(\cL) \subeq U^\rho_\cL(\cL)\cB(\cH) = \cB(\cH)_c^L$ is equivalent to
\[  \lim_{t \to \infty} P([-t,t])DP([-s,s]) = DP([-s,s])\quad
\hbox{for all}\; s>0.\]
By substituting  $D=\pi^\rho(B)$, we obtain the claim of the corollary.
\end{prf}
\begin{rem}
An application of this Corollary to Example~\ref{ExHHop}, produces the convenient formula
\[ \cA_\cL = \big\{ A \in \cB(\cH) \;\big|\;
(\forall s \in \R_+) \lim_{t \to \infty} P([-t,t])AP([-s,s]) = AP([-s,s])
\big\}, \]
to calculate $\cA_\cL$ in the identical representation of $\cA=\cB(\cH)$.
\end{rem}
\begin{defn}
With Borchers we write
$\cA^*_c$ for the closed subspace of the topological dual
$\cA^*$ of $\cA$, consisting of all
elements $\phi \in \cA^*$ for which the map
$G \to \cA^*, g \mapsto \phi \circ \alpha_g$
is norm continuous.
\end{defn}
Now $\cA^*_c$ is a folium, i.e., the predual of a $W^*$-algebra
(cf. \cite[Thm.~II.2.2]{Bo96} and \cite{Bo93}). A folium in $\cA^*$
determines
a representation of $\cA$
such that the folium is the set of normal functionals of the representation,
and this representation is unique up to quasi--equivalence (cf. \cite[Prop.~10.3.13]{KR86}).
 In fact, Borchers shows in
\cite[Thm.~III.2]{Bo83} that $\cA^*_c$ is the folium of functionals associated
to a covariant representation of $(\cA,G, \alpha)$. We can now identify that
representation in our picture:

\begin{prop}
\mlabel{BorchContS}
Let $(\cA,G,\alpha)$ be a
$C^*$-action, let $\cL$ be a full host algebra for $G$ for which the multiplier action
$\eta \: G \to \U(M({\cal L}))$ is strictly continuous.
Let the triple $(\cC, \eta_\cA, \eta_\cL)$  be constructed from
the universal  covariant $\cL$-representation
$(\pi_u,U_u)\in{\rm Rep}_\cL(\alpha,\cH_u)$,
where
$\al C.:=C^*\big(\pi_u(\cA) U_{u,\cL}(\cL)\big),$ and $\eta_\cA$ and $ \eta_\cL$ are defined as in
Theorem~\ref{thm:5.1}(b). We also recall
$\cC_\cL:=\eta_\cL(\cL) \cC\eta_\cL(\cL)\subseteq\cC$.
\begin{description}
\item[\rm(i)]  Then
there is a continuous surjection $\nu:(\cC_\cL)^*\to\cA^*_c$
given by $\nu(\varphi)(A):=\lim\limits_j\varphi\big(P\eta_\cA(A)F_jP\big)$ for $A\in\cA$, $\varphi\in(\cC_\cL)^*$, where
${(F_j)}$ is an approximate identity of $\cC_\cL$ and
$P$  
is the projection from \eqref{eq:appident}.
\item[\rm(ii)] $\cA^*_c$ is the predual of $\pi_u(\cA)''\subseteq \cB(\cH_u)$
where
$(\pi_u,U_u)\in{\rm Rep}_\cL(\alpha,\cH_u)$ is the the universal  covariant $\cL$-representation.
\end{description}
\end{prop}

\begin{prf} Since $(\cC, \eta_\cA, \eta_\cL)$ is constructed from a covariant $\cL$-representation it satisfies condition (CP3')
as well as (CP1) and (CP4). Thus by Proposition~\ref{LrepsofC}(iv) all $\cL$-representations of $\cC$ are covariant and produce covariant pairs $(\pi,U)\in{\rm Rep}_\cL(\alpha,\cH)$ for some $\cH$ by Proposition~\ref{LrepsofC}(iii).
An $\cL$-representation $\rho$ of $\cC$ is uniquely determined by its restriction to  $\cC_\cL=\eta_\cL(\cL)\cC\eta_\cL(\cL)$
via the relation
\[
\rho(C)=\slim_i\slim_j\rho(\eta_\cL(E_i)C\eta_\cL(E_j)),\;\; C\in\cC,\qquad\hbox{for any approximate identity $(E_j)$ of $\cL$.}
\]
Thus the $\cL$-representations of $\cC$ are unique extensions (on the same space) of representations of $\cC_\cL$,
and each of these produce a covariant pair $(\pi,U)\in{\rm Rep}_\cL(\alpha,\cH)$ and
$\pi(\cA)\cup U(G)\subset \rho(\cC_\cL)''$.

(i) As $\cC_\cL=\eta_\cL(\cL) \cC\eta_\cL(\cL)=P \cC^{**} P \cap \cC$
is a hereditary subalgebra of $\cC$,
it has the unique extension property for positive functionals, which implies
that the map
$\xi:(\cC_\cL)^*\to\cC^*$ by $\xi(\varphi)(C):=\varphi(PCP)$ for $C\in\cC$ is a linear injection
which is isometric on the selfadjoint functionals by Lemma~\ref{lem:a.2.1}. Since $\cC$ is a closed two-sided ideal
in $M(\cC)$, it is hereditary also in $M(\cC)$, so again positive functionals extend uniquely
(with the same norm), producing
another linear injection $\zeta:\cC^*\to M(\cC)^*$ by $\zeta(\varphi)(M)= \lim\limits_j\varphi\big(MF_j\big)$, $M\in M(\cC)$,
which is isometric on the selfadjoint functionals.
The restriction of functionals on $M(\cC)$ to $\eta_\cA(\cA)$ is a linear norm-reducing map.
The composition of these three maps is $\nu:(\cC_\cL)^*\to\cA^*$, and it is clear that it is linear and
norm continuous on the selfadjoint part of $(\cC_\cL)^*$. We have to establish that its range is $\cA^*_c$.

Recall from Remark~\ref{actionMC} that there is an action $\gamma \: G \to \Aut(M(\cC))$
which preserves $\cC$, $\eta_\cA(\cA)$ and $\eta_\cL(\cL)$ and
satisfies $\gamma(g) \circ \eta_\cA = \eta_\cA \circ \alpha_g$ for
$g \in G$.
It also preserves $M(\cC)_\cL$ hence $\cC_\cL=M(\cC)_\cL\cap\cC$ and defines a
strongly continuous $G$-action on
$\cC_\cL$. Since the GNS-representation of each state of $\cC_\cL$ is covariant by the observation above,
it follows that $\gamma$ defines a strongly
continuous action on the Banach space
$(\cC_\cL)^*$. Since the maps defining $\nu$ are all
covariant w.r.t.\ $\gamma$, it follows that $\nu$ is covariant w.r.t.\ $\gamma$ hence that its image consists of functionals with
norm continuous orbit maps
w.r.t.\ $\gamma$. Since the restriction of $\gamma$ to
$\eta_\cA(\cA)$ is compatible with $\alpha$, this implies that
the image of $\nu$ is $\alpha$-continuous, i.e. $\nu\big((\cC_\cL)^*\big)\subseteq\cA^*_c$.

For the converse inclusion, recall from Borcher's
Theorem \cite[Thm.~III.2]{Bo83} that $\cA^*_c$ is the folium of functionals associated to a
covariant representation $(\pi,U)\in{\rm Rep}(\alpha,\cH)$. As $\cL$ is a full host algebra, we have in fact that it is
a covariant $\cL\hbox{--representation}$.
We thus obtain an $\cL\hbox{--representation}$ $\rho$ of $\cC$,
hence of $\cC_\cL$, and it satisfies
\[ \rho(\eta_\cA(A)\eta_\cL(L)) = \pi(A) U_\cL(L)
\quad \mbox{ for } \quad A \in \cA, L \in \cL.\]
Since by Proposition~\ref{LrepsofC}(v) the map  $ \eta^*_\times \: \Rep_\cL(\cC,{\cal H}) \to \Rep_\cL(\alpha, \cH)$
by $\eta^*_\times(\rho):=\big(\tilde\rho \circ \eta_\cA, U^\rho\big)$ is injective, we get $\tilde\rho \circ \eta_\cA=\pi$.
Thus the normal states of $\pi$ are
the unique extensions of the normal states of $\rho$ restricted to $\cC_\cL$, hence
these are in the range of $\nu$.
Thus by linearity the whole predual of $\pi_u(\cA)''$
is in the range of $\nu$,
and as the predual of $\pi_u(\cA)''$ is $\cA^*_c$,
we conclude that we have the converse
inclusion, hence $\nu\big((\cC_\cL)^*\big)=\cA^*_c$.

(ii) The universal representation of $\cC_\cL$ produces the universal  $\cL$--representation
$\rho_u$,  and by Lemma~\ref{UnivLrep} we have $\eta^*_\times(\rho_u)={\big(\tilde\rho_u \circ \eta_\cA, U^{\rho_u}\big)}
={(\pi_u,U_u)}$. By Proposition~\ref{LrepsofC}(iii), the representation  $\pi_u=\tilde\rho_u \circ \eta_\cA$ of $\cA$
is nondegenerate. Moreover $\pi_u(\cA)\cup U_u(G)\subset \rho_u(\cC_\cL)''$, and
the predual of $\rho_u(\cC_\cL)''$ is $(\cC_\cL)^*$
(cf.\ \cite[ Prop.~III.5.2.6]{Bla06}).
The map  $\nu:(\cC_\cL)^*\to\cA^*_c$  in the previous part is just the restriction of the normal functionals of
$\rho_u(\cC_\cL)''$ to $\pi_u(\cA)\subset\rho_u(\cC_\cL)''$, hence $\cA^*_c=\nu((\cC_\cL)^*)=\rho_u(\cC_\cL)''_*\restriction
\pi_u(\cA)=\pi_u(\cA)''_*$ because $\rho_u(\cC_\cL)''_*\!\restriction\!\pi_u(\cA)''=\pi_u(\cA)''_*$ and normal functionals
of $\pi_u(\cA)''$ are uniquely determined by their restrictions to $\pi_u(\cA)$.
Here we used the fact that if  $\cN\subset\cM$ is an inclusion of
von Neumann algebras, then
$\cN_*=\cM_*\!\!\restriction\!\cN$ (\cite[Thm~4.2.10]{Mu90})
to extend normal functionals on $\cN$ to normal functionals on $\cM$.
\end{prf}

\section{Examples}
\label{SectExmp}

We have already seen in Examples~\ref{FullCPH-discontAlph}
and \ref{nonfullMax} that crossed product hosts can exist
for discontinuous actions. Here we want to develop further examples.

We will analyze examples which have an interest for the physics
of bosonic particles. Let $(X,\sigma)$ be a nondegenerate symplectic space over $\R$,
and let $\Sp(X,\sigma)$ denote the group of symplectic transformations of it.
The quantum system based on this has field algebra $\cA$ being either the Weyl algebra $\ccr X,\sigma.$
(cf.~\cite{Ma68} and Example~\ref{NumberOp1} above), or the Resolvent algebra $\rsl$ (cf.~\cite{BG08}).
Both are defined through generators satisfying a set of the relations. Let
$\{\delta_x\mid x\in X\}$ (resp. $\{R(\lambda,x)\mid x\in X,\, \lambda\in
\R^\times\}$ )
denote the generators of  $\ccr X,\sigma.$ (resp. $\rsl$). Then we define
an automorphic action $\alpha:\Sp(X,\sigma)\to\Aut(\cA)$ by $\alpha\s T.(\delta_x):=\delta_{T(x)}$
(resp. $\alpha\s T.(R(\lambda,x)):=R(\lambda,T(x))$ for $T\in \Sp(X,\sigma)$.
In the following we will be interested
in particular one-parameter subgroups of
$\Sp(X,\sigma)$ whose corresponding action $\alpha$
is not strongly continuous.

\begin{ex}
\mlabel{QM-Weyl}
 (A free quantum particle in one dimension - Weyl algebra)\\
Let $X=\R^2$ with symplectic form $\sigma((x_1,y_1),(x_2,y_2))=x_1y_2-x_2y_1$ and
fix the symplectic basis ${(p ,q)}= {((1,0),(0,1))}$.
Let  $\cA=\ccr X,\sigma.$ and let $\pi:\cA\to B(L^2(\R))$ be the Schr\"odinger representation
w.r.t.\ this basis, i.e. $\pi(\delta_{tp})=e^{itP}$, $\pi(\delta_{tq})=e^{itQ}$
where $Q,\,P$ are the usual operators of multiplication by $x$ and $i{d\over dx}$ on the appropriate domain in
$L^2(\R)=:\cH$.
We consider the action $\alpha:\R\to\Aut(\cA)$ determined by
\[
\pi\big(\alpha_t(\delta_x)\big)=e^{itQ^2}\pi(\delta_x)e^{-itQ^2}=\pi(\delta\s T_t(x).)\quad\hbox{for}\quad
t\in\R,\; x\in X
\]
where $T_t(sq+rp)= sq+r(p+2tq)$, $s,r\in\R$. Note that $T_t\in\Sp(X,\sigma)$, and that a  Fourier transform
converts this to the usual  action of the free Hamiltonian $P^2$.
The action $\alpha$ is not strongly
continuous because $t\mapsto\delta_{p+tq}=\alpha_t(\delta_p)$ is not norm continuous, as ${\|\delta_x-\delta_y\|}=2$
if $x\not=y$. If we take $\cL=C^*(\R)\cong C_0(\R)$ then the $\cL\hbox{--representations}$ of $\R$ are just the strong operator continuous
unitary representations.
Thus, for $U_t = e^{itQ^2}$, the pair  $(\pi,U)$
is in $ {\rm Rep}_\cL(\alpha,\cH)$,
and so it defines a
triple $(\cC, \eta_\cA, \eta_\cL)$ satisfying (CP1), (CP3') and (CP4).
Now $U_\cL(f)=\hat{f}(Q^2)$ for $f\in L^1(G)$, where $\hat{f}$ denotes the Fourier transform of $f$.
Thus $U_\cL(\cL)=\{h(Q^2)\,\mid\,h\in C_0(\R)\}$ and, as the spectrum of $Q^2$ is positive, this is clearly a factor algebra of
$C_0(\R)\cong C^*(\R)$. Hence $\eta_\cL=U_\cL$ is not faithful though $\eta_\cA=\pi$ is faithful (as $\cA$ is simple).
 Next, we want to determine the algebra $\cA_\cL\subseteq\cA$. By Proposition~\ref{propALchar}(iv) we already have that
$\delta_{tq}\in \cA_\cL$. To check whether $A=\delta_{tp}\in \cA_\cL$, we consider the condition
\[
\big\| \big(\eta_\cL(E_j) -\1\big)\eta_\cA(A) \eta_\cL(L)
 \big\| \to 0 \;\;\mbox{ for all} \;\;L \in \cL
 \]
from Proposition~\ref{propALchar}(iii).
Consider the approximate identity $(E_n)_{n\in\N}$ of $C_0(\R)\cong C^*(\R)$ where
$E_n$ is a smooth bump function which is $1$ on ${[-n,n]}$ and zero outside ${[-n-1,n+1]}$.
Let $\U_\cL(L)=f(Q^2)$ for some $f\in C_0(\R)$. Then
\[ \eta_\cA(A) \eta_\cL(L)=e^{itP}f(Q^2)=f((Q-t\1)^2)e^{itP}
\in \{h(Q)\mid h\in C_0(\R)\}\cdot e^{itP} \]
and so
\[
\big\| \big(\eta_\cL(E_n) -\1\big)\eta_\cA(A) \eta_\cL(L)
 \big\| = \big\| \big(E_n(Q^2) -\1\big)f((Q-t\1)^2)e^{itP}
 \big\|\to 0 \,.
 \]
 Thus $\delta_{tp}\in \cA_\cL$ and as $\cA$ is generated as a $C^*$-algebra by $\delta_{tq}$ and $\delta_{tp}$ for all $t$,
 it follows that $\cA_\cL=\cA$, hence that $(\pi,U)$ is  a cross representation for $(\alpha,\cL)$.
 Note that the ideal $\cA^{(\cL)}_0$ from Proposition~\ref{NContAlg}
 is zero as $\cA$ is simple, and the
 action $\alpha$ is not strongly continuous. Moreover, by irreducibility $\pi$ determines
 the implementers $U:\R\to \U(\cH)$ up to a $\T\hbox{--valued}$ multiplier. As the second Moore
cohomology  $H^2(\R,\T)$  of $\R$ is trivial (\cite[Thm.~7.38]{Var07}),
it follows that the multiplier is a coboundary, i.e.
 the implementers can be rewritten to produce the representation $U:\R\to \U(\cH)$ given above.
 The scalar factors involved leave $U_\cL(\cL)$ unchanged, so the Schr\"odinger representation
 remains a cross representation,
regardless of the choice of unitary implementers.
Now any regular representation
of $\cA$ is a direct sum of copies of the Schr\"odinger
representation, hence a cross representation
by Lemma~\ref{lem:perm}(iii), and therefore
$\cA_\cL=\cA$ in these.

 This is one situation referred to in the introduction of Section~\ref{NonX}. That is, from a physics
 point of view one is only interested in regular representations, so a full crossed product host
 will not be useful if it produces  covariant $\cL\hbox{--representations}$
 $(\pi,U)\in {\rm Rep}_\cL(\alpha,\cH)$ for which $\pi$ is nonregular.
 Ideally one seeks a crossed product host whose universal
$\cL$-representation   $\rho_u\in{\rm Rep}_\cL(\cC,\cH_u)$ produces a regular representation
$\tilde\rho_u\circ\eta_\cA$ of $\cA$.
 \end{ex}

 \begin{ex} (A free quantum particle in one dimension - Resolvent algebra)\\
 We repeat the previous example, with the only difference that we replace
 the Weyl algebra $\ccr X,\sigma.$ with the Resolvent algebra $\rsl=:\cA$.
 Thus for the the
Schr\"odinger representation (which is faithful on $\cA$) we have
\begin{eqnarray*}
\pi\big( R(\lambda, sq+rp) \big)&=& \big(i\lambda\1-(sQ+rP)\big)^{-1}\quad \hbox{for}\quad r,\,s\in\R
\\[1mm]
\hbox{and}\qquad\qquad
\pi\big(\alpha_t(R(\lambda, x) )\big)&=&e^{itQ^2}\pi(R(\lambda, x))e^{-itQ^2}=\pi(R(\lambda,T_t( x)))\quad\hbox{for}\quad
t\in\R,\; x\in X
\end{eqnarray*}
with the $T_t\in\Sp(X,\sigma)$ we had above.
By \cite[Thm.~5.3(ii)]{BG08},
the action  $\alpha$ is not strongly continuous.
As above, the pair $(\pi,U)$ is in $ {\rm Rep}_\cL(\alpha,\cH),$
and so it defines a triple $(\cC, \eta_\cA, \eta_\cL)$ satisfying (CP1), (CP3') and (CP4).
Now $U_\cL(f)=\hat{f}(Q^2)$ for $f\in L^1(G)$, where $\hat{f}$ denotes the Fourier transform of $f$.
Thus $U_\cL(\cL)=\{h(Q^2)\,\mid\,h\in C_0(\R)\}= C^*((\1+Q^2)^{-1}).$
However $(\1+Q^2)^{-1}=-(i\1-Q)^{-1}(i\1+Q)^{-1}=-\pi\big(R(1,q)R(1,-q)\big)\in\pi(\cA)$, hence
$U_\cL(\cL)\subseteq\pi(\cA)$ and so $\eta_\cL(\cL)\subseteq\eta_\cA(\cA)$ and thus
$\cC\subseteq\eta_\cA(\cA)$. Since $\eta_\cA(\cA)\subseteq M(\cC)$ it follows that $\cC$ is a two-sided
ideal of $\eta_\cA(\cA)$. By \cite[Thm.~3.8]{BG08}, we have that
\begin{eqnarray*}
\lbr R(\lambda, x)\rsl\rbr&=&
\lbr\rsl R(\lambda, x)\rbr\quad\hbox{for all}\quad x\in X, \mbox{ hence } \\
\lbr\eta_\cL(\cL)\eta_\cA(\cA)\rbr
&=& \eta_\cA\Big(\lbr R(1,q)R(1,-q)\rsl\rbr  \Big)
\\[1mm]
&=& \eta_\cA\Big(\lbr\rsl R(1,q)R(1,-q)\rbr \Big)
=\lbr\eta_\cA(\cA)\eta_\cL(\cL)\rbr=\cC
\end{eqnarray*}
and thus
\[ \eta_\cA(\cA)\eta_\cL(\cL) \subeq \cC=\eta_\cL(\cL)\cC.\]
By Proposition~\ref{propALchar}(ii) we therefore obtain $\cA_\cL=\cA$, so that
 $(\pi,U)$ is  a cross representation for $(\alpha,\cL)$.
In this case we have for the ideal
\[
\cA^{(\cL)}_0:=\{A\in\cA\mid\eta_\cA(A)\in\eta_\cL(\cL) M(\cC)\eta_\cL(\cL)\}
=\{A\in\cA\mid\eta_\cA(A)\in\eta_\cL(\cL) \eta_\cA(\cA)\eta_\cL(\cL)\}
=\eta_\cA^{-1}(\cC)
\]
since $\cC$ is an ideal of $\eta_\cA(\cA)$, and it is
\[
\cC=\lbr\eta_\cL(\cL)\eta_\cA(\cA)\rbr
= \eta_\cA\Big(\lbr R(1,q)R(1,-q)\rsl\rbr \Big)
=\lbr\eta_\cA(\cA)\eta_\cL(\cL)\rbr
=\eta_\cL(\cL)\eta_\cA(\cA)\eta_\cL(\cL)\,.
\]
Now $\eta_\cA$ is faithful because $\pi$ is faithful and $U_\cL(\cL)$ is
nondegenerate, hence
\[ \cA^{(\cL)}_0=\lbr R(1,q)R(1,-q)\rsl\rbr.\]
 By the analogous reasoning as in
the preceding example,
we find that every regular representation $\pi$ of $\cA$ is a cross representation
${(\pi,U)}\in {\rm Rep}_\cL^\times(\alpha,\cH)$.
\end{ex}

 \begin{ex}
 \label{FockCross}
 (The Fock representation for a bosonic quantum field).\\
 We continue with the setting of Example~\ref{NumberOp1}.
 Thus $(\cH,\sigma)$ consists of a nonzero complex Hilbert space
 $\cH$ and
$\sigma:\cH\times\cH\to\R$ is $\sigma(x,y):={\rm Im}{\langle x,y\rangle}$.
We take the Weyl algebra $\cA=\ccr \cH,\sigma.$, and the action
$\alpha:\Sp(\cH,\sigma)\to\Aut(\cA)$ by $\alpha\s T.(\delta_x):=\delta_{T(x)}$.
Note that $\U(\cH)\subset \Sp(\cH,\sigma)$.
We briefly recall the  Fock representation $\pi_F:\cA\to\cB(\cF(\cH))$. The bosonic Fock space is
\[
\cF(\cH):=\bigoplus_{n=0}^\infty\otimes_s^n\cH\,,\quad
\otimes_s^n\cH\equiv\hbox{symmetrized Hilbert tensor product of
$n$ copies of $\cH$}
\]
with a convention $\otimes_s^0\cH:=\C.$ The finite particle space
$\cF_0(\cH):={\Spann\{\otimes_s^n\cH\mid n=0,1,\cdots\}}$ is dense in $\cF(\cH)$.
For each $f\in\cH$ we define  on $\cF_0(\cH)$ a (closable)
 creation operator $a^*(f)$ by
\[
a^*(f)\, \big(v_1\otimes_s\cdots\otimes_s v_n\big):= \sqrt{n+1}\,
S\big(f\otimes v_1\otimes\cdots\otimes v_n\big)
=:\sqrt{n+1}\,f\otimes_s v_1\otimes_s\cdots\otimes_s v_n
\]
where $S$ is the symmetrizing operator. Define on $\cF_0(\cH)$ an essentially selfadjoint
operator by $\phi(f):=\big(a^*(f)+a(f)\big)/\sqrt{2}$
 where $a(f)$ is the adjoint of $a^*(f)$.
The Fock representation $\pi_F:\cA\to
\cB\big(\cF(\cH)\big)$ is then defined by $\pi_F(\delta_f)=\exp(i\overline{\phi(f)})$,
for all $f\in\cH,$ and it is irreducible. Given a strong operator continuous
one-parameter group
$t\to U_t=\exp(itA)\in \U(\cH)$, where $A$ is selfadjoint,
define a unitary group
$t\to\Gamma(U_t)\subset U\big(\cF(\cH)\big)$ by $\Gamma(U_t)\big(v_1\otimes_s\cdots\otimes_s v_n\big):=
\big(U_tv_1\otimes_s\cdots\otimes_s U_tv_n\big)$ which is strong operator continuous, with
generator given on $\cF_0(\cH)$ by
\[
\dd\Gamma(A)\big(v_1\otimes_s\cdots\otimes_s v_n\big)=
Av_1\otimes_sv_2\otimes_s\cdots\otimes_s v_n+\cdots+v_1\otimes_s\cdots\otimes_s v_{n_1}\otimes_s Av_n.
\]
We then have covariance $\pi_F(\alpha\s U_t.(A))=\Gamma(U_t)\pi_F(A)\Gamma(U_t)^*$.

In this example we want to prove that if $A\geq 0$, where zero is isolated in its spectrum, then  ${(\pi_F,\Gamma(U))}$ is  a
cross representation for $(\alpha,C^*(\R))$.
Such systems are analyzed in \cite{Bo96},
and are important for quantum field theory.

To start, we assume first that
\[ \sigma(A) = \sigma_p(A) \subeq \{n\kappa\mid n=1,\,2,\ldots  \}
\quad \mbox{ for some } \quad \kappa > 0.\]
Let $\cL=C^*(\R)$, then we  prove that ${(\pi_F,\Gamma(U))}$ is  a
cross representation with this choice of $A$. Now
\[
U_\cL(\cL)=\{f(A)\mid f\in C_0(\R)\}=\oline{\Spann\{P_\lambda \mid
\lambda \in \kappa\N\}},
\]
where $P_{\lambda}$ denotes the projection on the
$\lambda$-eigenspace of $A$.
As $\cH$ has a Hilbert basis of eigenvectors of $A$, by considering the spaces
$P\s\lambda_{1}.\!\!\cH\otimes_s\cdots\otimes_s P\s\lambda_{n}.\!\!\cH$,
we see that $\cF(\cH)$
has a Hilbert basis of eigenvectors for $\dd\Gamma(A)$, and that
$\sigma(\dd\Gamma(A)) = \sigma_p(\dd\Gamma(A))\subeq \kappa \N_0$.

Let $\mu\in\sigma(\dd\Gamma(A))$ and $\hat P_\mu$ denote the projection onto
its eigenspace. Then
\[
\Gamma(U)_\cL(\cL)=\{f(\dd\Gamma(A))\mid f\in C_0(\R)\}=
\oline{\Spann\big\{\hat P_\mu\mid \mu\in \sigma(\dd\Gamma(A))\big\}}.
\]
By Lemma~\ref{lem:a.1}(ii)(b),  ${(\pi_F,\Gamma(U))}$ will be  a
cross representation for $(\alpha,C^*(\R))$ if
\[
\lim_{t\to 0}\Gamma(U_t)\,\pi_F(B)\,\hat P_\mu
=\pi_F(B)\,\hat P_\mu \quad \mbox{ for all } \quad
B\in\cA\;\mbox{ and } \;\mu\in\sigma(\dd\Gamma(A)).\]
 Let $B=\delta_f$ where $f$ is a nonzero eigenvector:
$Af=\lambda f$. Let $\cH=\cH_0\oplus\cH_1$ where $\cH_0=\C f$, then
$\cA=\ccr \cH,\sigma.\cong \cA_0\otimes\cA_1$ where $\cA_i=\ccr \cH_i,\sigma.$,
$i=0,1$ (cf.~\cite[3.4.1]{Ma68}).
Moreover $\pi_F=\pi_0\otimes\pi_1$ on $\cF(\cH)\cong\cF(\cH_0)\otimes\cF(\cH_1)$ where
$\pi_i:\cA_i\to \cB(\cF(\cH_i))$ are the respective
Fock representations with respective second quantization
maps $\Gamma_i$. Then
\[
\Gamma(U_t)=\Gamma_0(U_t\!\restriction\!\cH_0)\otimes\Gamma_1(U_t\!\restriction\!\cH_1)\quad\hbox{and}\quad
\dd\Gamma(A)=\dd\Gamma_0(A\!\restriction\!\cH_0)\otimes\1+\1\otimes
\dd\Gamma_1(A\!\restriction\!\cH_1)\,.
\]
Since $\sigma\big(\dd\Gamma_0(A\!\restriction\!\cH_0)\big)
=\{m\lambda\mid m=0,1,\ldots\}$, we have
\begin{equation}
  \label{hatPdecomp}
\hat P_\mu=\sum\Big\{\hat P^0_{m\lambda}\otimes \hat P^1_{\mu-m\lambda}\;\Big|\; m=0,1,\ldots\quad\hbox{such that}\quad
\mu-m\lambda\in \sigma_p\big(  \dd\Gamma_1(A\!\restriction\!\cH_1)   \big)\Big\}
\end{equation}
where $\hat P^i_\nu$ denotes the projection onto the $\nu$-eigenspace
of $\dd\Gamma_i(A\!\restriction\!\cH_i)$.
As $\sigma(\dd\Gamma_1(A\!\restriction\!\cH_1)) \subeq \N_0 \kappa$
and $\lambda\geq\kappa>0$,
this sum is finite for each $\mu$.
\begin{eqnarray*}
\Gamma(U_t)\,\pi_F(\delta_f)\,\hat P_\mu&=&
\sum_{m=0}^M\Gamma_0(U_t\!\restriction\!\cH_0)\pi_0(\delta_f) \hat P^0_{m\lambda}\;\otimes \;\Gamma_1(U_t\!\restriction\!\cH_1)\hat P^1_{\mu-m\lambda} \\[1mm]
&=&
\sum_{m=0}^Me^{it(\mu-m\lambda)}\Gamma_0(U_t\!\restriction\!\cH_0)\pi_0(\delta_f)
\hat P^0_{m\lambda}\;\otimes \;\hat P^1_{\mu-m\lambda}  \\[1mm]
&\longrightarrow & \sum_{m=0}^M\pi_0(\delta_f) \hat P^0_{m\lambda}\otimes
\hat P^1_{\mu-m\lambda}
= \pi_F(\delta_f)\,\hat P_\mu
\end{eqnarray*}
as $t\to 0$, since each $\hat P^0_{m\lambda}$ has finite rank,
using Theorem~\ref{compactExist}.
Thus
\[
\lim_{t\to 0}\Gamma(U_t)\,\pi_F(\delta_f)\,\hat P_\mu
=\pi_F(\delta_f)\,\hat P_\mu\quad\hbox{i.e.}\quad \pi_F(\delta_f)\in\cA_\cL
\]
for all eigenvectors $f$ of $A$ by Proposition~\ref{propALchar}(ii).
As $\cA_\cL$ is a $C^*$-algebra, this means that
\begin{equation}
  \label{eq:lastrel}
\lim_{t\to 0}\Gamma(U_t)\,\pi_F(B)\,\hat P_\mu
=\pi_F(B)\,\hat P_\mu\quad\hbox{for all }\quad B\in C^*\{\delta_f\mid f\quad
\hbox{an eigenvector of}\;\; A\}.
\end{equation}
So for this part it will suffice
to show that
\begin{equation}
  \label{eq:lastrel2}
\pi_F\big(C^*\{\delta_f\mid f\quad
\hbox{an eigenvector of}\;\; A\}\big)\hat P_\mu
=\pi_F(\cA)\hat P_\mu\quad \mbox{ for all } \quad \mu\in\sigma(\dd\Gamma(A)).
\end{equation}
 Now $\hat P_\mu$ commutes with the projections $\tilde P_n$ onto the $n$-particle spaces
(as $\Gamma(U_t)$ preserve these spaces). Note that as $\sigma(A)\geq\kappa>0$,
each eigenspace $\hat P_\mu\cF(\cH)$ of $\dd\Gamma(A)$ (which is made up of the spaces
$P\s\lambda_{1}.\!\!\cH\otimes_s\cdots\otimes_s P\s\lambda_{n}.\!\!\cH$),
contains $n$-particle spaces of order at most $\mu/\kappa$.
So by the finite sums above for $\hat P_\mu$ in (\ref{hatPdecomp}),
we see that $\hat P_\mu\tilde P_n\not=0$ for only finitely many $n$. Thus it suffices to show that \eqref{eq:lastrel}
holds if we replace $\hat P_\mu$ by $\tilde P_n$. In fact, as the eigenvectors
of $A$ span a dense subspace of
$\cH$, for (\ref{eq:lastrel2}) it suffices to prove that if $v_n\to 0$ in $\cH$, then
$\big\|\pi_F(\delta_{v_n}-\1)\tilde P_n\big\|\to 0$. Now the $n\hbox{-particle}$ vectors are analytic vectors
for $\phi(f)$ (cf. proof of \cite[Thm.~X.41]{ReSi75}). Thus
\[
\pi_F(\delta_{v_n}-\1)\tilde P_n=\sum_{k=1}^\infty{(i\phi(v_n))^k\over k!}\,\tilde P_n
\]
and using the estimates $\|\phi(v)^k\psi\|\leq 2^{k/2}((n+k)!)^{1/2}\|v\|^k\|\psi\|$ if
$\psi\in\tilde P_n\cF(\cH)$ (cf. proof of \cite[Thm~X.41]{ReSi75}) we conclude that
\[
\Big\|\pi_F(\delta_{v_n}-\1)\tilde P_n\Big\|\leq
\sum_{k=1}^\infty{(2^{k/2}((n+k)!)^{1/2}\over k!}\,\|v_n\|^k
\leq\|v_n\|\sum_{k=1}^\infty{(2^{k/2}((n+k)!)^{1/2}\over k!}
\]
if $\|v_n\|\leq 1$. As the last series converges, it is clear that
$\big\|\pi_F(\delta_{v_n}-\1)\tilde P_n\big\|\to 0$ as $v_n\to 0$.
We conclude that
\[
\lim_{t\to 0}\Gamma(U_t)\,\pi_F(B)\,\hat P_\mu
=\pi_F(B)\,\hat P_\mu
\]
for all $B\in\cA$ and $\mu\in\sigma_p(\dd\Gamma(A))$, i.e.  ${(\pi_F,\Gamma(U))}$ is  a
cross representation for $(\alpha,C^*(\R))$.

Next, assume that $A$ is a positive operator on $\cH$ with strictly positive spectrum
$\sigma(A)\subseteq[a,\infty)$, $a>0$, in the same context as above, then we want to show that
here too,  ${(\pi_F,\Gamma(U))}$ is  a cross representation for $(\alpha,C^*(\R))$.
Given the spectral resolution $A=\int_{\sigma(A)}\lambda\,dP(\lambda)$, let
$\varepsilon>0$
and define
\[ B_\varepsilon:=\int_{\sigma(A)}f_\varepsilon(\lambda)\,dP(\lambda)
\quad \mbox{ where } \quad
f_\varepsilon(\lambda)=\eps
+ \lfloor \lambda/\varepsilon\rfloor \varepsilon, \]
i.e. the right continuous step function with steps of size $\varepsilon$.
As $\|\id_{\R_+}-f_\varepsilon\|_\infty\leq\varepsilon$ we  have $\|A-B_\varepsilon\|\leq\varepsilon$.
However $\sigma(B_\varepsilon)=\sigma_p(B_\varepsilon)\subeq \eps\N$, so
$\sigma_p(\dd\Gamma(B_\varepsilon))=\sigma(\dd\Gamma(B_\varepsilon))$. Thus by the previous part
 ${(\pi_F,\Gamma(V))}$ is  a
cross representation for $(\beta,C^*(\R))$ where $V_t:=\exp(itB_\varepsilon)$ and
$\beta:\R\to\Aut\cA$ is given by $\beta_t(\delta_v):=\delta\s V_tv.$.
Now as $-1$ is not in the spectra of $\dd\Gamma(A)$ and $\dd\Gamma(B_\varepsilon)$
we have
\begin{eqnarray*}
\Gamma(U)_\cL(\cL)=C^*\big((\1+ \dd\Gamma(A))^{-1}\big)\quad &\hbox{and}&\quad
\Gamma(V)_\cL(\cL)=C^*\big((\1+ \dd\Gamma(B_\varepsilon))^{-1}\big)
 \\[1mm]
 (\1+ \dd\Gamma(A))^{-1} - (\1+ \dd\Gamma(B_\varepsilon))^{-1} &=&
(\1+ \dd\Gamma(A))^{-1} \dd\Gamma(B_\varepsilon-A)(\1+ \dd\Gamma(B_\varepsilon))^{-1} \,.
\end{eqnarray*}
and the factors in this last expression commute.
If $\dd\Gamma(B_\varepsilon-A)(\1+ \dd\Gamma(B_\varepsilon))^{-1}$ is bounded,  this implies that
$(\1+ \dd\Gamma(B_\varepsilon))^{-1}\in\Gamma(U)_\cL(\cL)\cB(\cF(\cH))$.
Now the factors in the products above are positive, and commute, therefore
their products are positive. Thus as $0\leq \dd\Gamma(B_\varepsilon-A)\leq \dd\Gamma(B_\varepsilon)$
it follows that
\[ 0\leq \dd\Gamma(B_\varepsilon-A)(\1+ \dd\Gamma(B_\varepsilon))^{-1}\leq
\dd\Gamma(B_\varepsilon)(\1+ \dd\Gamma(B_\varepsilon))^{-1}\in \cB(\cF(\cH))\]  hence
$\dd\Gamma(B_\varepsilon-A)(\1+ \dd\Gamma(B_\varepsilon))^{-1}$ is bounded, so
\[
\Gamma(V)_\cL(\cL)=C^*((\1+ \dd\Gamma(B_\varepsilon))^{-1})\subset\Gamma(U)_\cL(\cL)\,\cB(\cF(\cH)).
\]
Next observe that by strict positivity $\sigma(A)\subseteq[a,\infty)$, $a>0$,
we can find an $\eps>0$ such that  $B_\varepsilon-A\leq A$   (e.g. $\eps< a/2$). Thus by
$0\leq \dd\Gamma(B_\varepsilon-A)\leq \dd\Gamma(A)$ we get
\[ 0\leq \dd\Gamma(B_\varepsilon-A)(\1+ \dd\Gamma(A))^{-1}\leq
\dd\Gamma(A)(\1+ \dd\Gamma(A))^{-1}\in \cB(\cF(\cH))\]
and so, as above, we get from the resolvent identity that
\[\Gamma(U)_\cL(\cL)\subset
\Gamma(V)_\cL(\cL)\,\cB(\cF(\cH)).\]
Together, these two inclusions produce:
\begin{eqnarray*}
\pi_F(\cA)\,\Gamma(U)_\cL(\cL)&\subseteq&
\pi_F(\cA)\,\Gamma(V)_\cL(\cL)\,\cB(\cF(\cH))\\[1mm]
&\subseteq&\Gamma(V)_\cL(\cL)\,\cB(\cF(\cH))\qquad\hbox{
as ${(\pi_F,\Gamma(V))}$ is  a
cross rep. for $\beta$}\\[1mm]
&\subseteq&\Gamma(U)_\cL(\cL)\,\cB(\cF(\cH)).
\end{eqnarray*}
Thus ${(\pi_F,\Gamma(V))}$ is  a
cross representation for $(\alpha,C^*(\R))$.

Finally, we assume that  $A\geq 0$, where zero is isolated in its spectrum.
Decompose $\cH=\cH_0\oplus\cH_1$ where $\cH_0$ is the kernel of $A$,
and hence $A$ is strictly positive on $\cH_1$. Then as above
$\cA=\ccr \cH,\sigma.\cong \cA_0\otimes\cA_1$ where $\cA_i=\ccr \cH_i,\sigma.$,
$i=0,1$ and
 $\pi_F=\pi_0\otimes\pi_1$ on $\cF(\cH)\cong\cF(\cH_0)\otimes\cF(\cH_1)$ where
$\pi_i:\cA_i\to \cB(\cF(\cH_i))$ are the respective
Fock representations with respective second quantization
maps $\Gamma_i$. Now $\Gamma_0(U_t)=\1$ so  $\Gamma(U_t)=\1\otimes\Gamma_1(U_t)$ hence
$\Gamma(U)_\cL(\cL)=\1\otimes\Gamma_1(U)_{\cL}(\cL)$. Let $f=f_0+f_1\in \cH_0\oplus\cH_1$,
then for each $L\in\cL$ we have
\begin{eqnarray*}
\Gamma(U_t)\,\pi_F(\delta_f)\,\Gamma(U)_{\cL}(L) &=& \pi_0(\delta_{f_0})\otimes\Gamma_1(U_t)\,\pi_1(\delta_{f_1})\,\Gamma_1(U)_{\cL}(L)
\\[1mm]
&\to& \pi_0(\delta_{f_0})\otimes\pi_1(\delta_{f_1})\,\Gamma_1(U)_{\cL}(L)=\pi_F(\delta_f)\,\Gamma(U)_{\cL}(L)
\end{eqnarray*}
as $t\to 0$ by the previous part for strictly positive $A$. Thus  ${(\pi_F,\Gamma(U))}$ is  a
cross representation for $(\alpha,C^*(\R))$.

We do not know whether  ${(\pi_F,\Gamma(U))}$ is still a
cross representation for $(\alpha,C^*(\R))$ if one assumes only that  $A\geq 0$.
\end{ex}

\section{Discussion}

Above, we extended crossed products to singular actions
 $\alpha \: G \to \Aut(\cA)$, relative to the choice of
a host algebra $\cL$. There is still much further to be explored,
in particular, we need to analyze in detail crossed product hosts when
we have a spectral condition included. Host algebras have been constructed
explicitly for group representations subject to a spectral condition
(cf.~\cite{Ne00, Ne08}), and these can now easily be included in the
constructions above of crossed product hosts (work is in progress on this
subject~\cite{GrN12}).
This naturally will have to connect with the deep work of Borchers~\cite{Bo96}
and the spectral theory of Arveson~\cite{Ar74} which deal with these
topics.

Further directions concern the development
of host algebras for non--locally compact groups.
There are also numerous dynamical systems for physical systems which need to be analyzed to
establish whether the representations used are cross representations or not.

\section*{Appendix}
\appendix

\section{Some facts on multiplier algebras} \mlabel{app:1}

\begin{lem}\mlabel{lem:a.1}
Let $\cS$ and $\cT$ be $C^*$-algebras
and $(E_i)_{i \in I}$ be an  approximate identity in $\cS$.
For a morphism  $\zeta \: \cS \to M(\cT)$ of $C^*$-algebras,
the following assertions hold:
\begin{description}
\item[\rm(i)] $\cT^L_c :=  \{ T \in \cT \mid
\lim \zeta(E_i)T = T \} = \zeta(\cS)\cT$ is a closed
right ideal of $\cT$.
\item[\rm(ii)]  Suppose that $\cT = \cB(\cH)$ and that
$\zeta$ is a non-degenerate representation of
$\cS$ on $\cH$, then
\begin{description}
\item[\rm(a)] $\cB(\cH)^L_c = \zeta(\cS) \cB(\cH)\supeq \cK(\cH)$.
\item[\rm(b)] If $G$ is a locally compact group,
$\cS = C^*(G)$ and $U \: G \to \U(\cH)$ is the
unitary representation defined by $\zeta$, then
\[ \cB(\cH)^L_c = \zeta(C^*(G))\,\cB(\cH)
= \{ A \in \cB(\cH) \mid \lim_{g \to \1} U_g A = A\}.\]
\item[\rm(c)] If, in addition, $G = \R$,
and $P$ is the spectral measure corresponding to $U$, then
\[ \cB(\cH)^L_c=\zeta(C^*(\R))\,\cB(\cH)
= \{ A \in \cB(\cH) \mid \lim_{t \to \infty} P([-t,t])A = A\}.\]
Moreover, for
$A \in \cB(\cH)$, we have
\[ A \zeta(C^*(\R)) \subeq \cB(\cH)_c^L
\quad \Longleftrightarrow \quad
A P([-t,t]) \in \cB(\cH)_c^L\quad \mbox{ for each} \quad
t > 0.\]
\end{description}
\end{description}
\end{lem}

\begin{prf} (i) In view of \cite[Thm.~5.2.2]{Pa94},
the subset $\zeta(\cS)\cT$ is a closed subspace of
$\cT$ and it obviously is a right ideal.
It clearly is contained $\cT^L_c$.
Conversely, every element of the form
$\zeta(E_i)T$ is contained in $\zeta(\cS)\cT$,
so that $\zeta(E_i)T \to T$ implies that
$T \in \zeta(\cS)\cT$. This proves~(i).

(ii) (a) In view of the closedness of $\cB(\cH)^L_c$,
it suffices to show that it contains all rank
one operators of the form
$A = \la \cdot, x \ra y$ for $x,y \in \cH$.
Since the representation $\zeta$ is non-degenerate,
$\zeta(E_i) y \to y$, and therefore $\zeta(E_i)A \to A$.

(b) Since the action of $G$ on $C^*(G)$ by left multipliers
is continuous, $A \in \zeta(C^*(G))\cB(\cH)$
implies $\lim_{g \to \1} U_g A = A$.

Conversely, by the fact that $G$ acts on $\cB(\cH)$ by
unitary multipliers, the subspace
\[ \cB := \{ A \in \cB(\cH) \mid \lim_{g \to \1} U_g A = A \} \]
is closed and we thus obtain a strongly continuous
representation
$\beta(g)A = U_g A$ of $G$ on $\cB$.
Let $\hat\beta \: L^1(G) \to B(\cB)$ be the corresponding
integrated representation. For $f \in L^1(G)$, $v,w \in \cH$
we then have
\[ \la \hat\beta(f)A v,w \ra =
\int_G f(g) \la U_g A v,w \ra\, dg =
\la \zeta(f)Av,w\ra\quad\hbox{for all}\;\; A\in\cB,
\]
showing that $\zeta(f)A = \hat\beta(f)A$ for $f \in L^1(G)$ and $A\in\cB$.

For each open $\1$-neighborhood $V$ in $G$, let  $E_V
= E_V^* \geq 0$ be an element of $L^1(G)$ supported in $V$ with
$\int_G E_V(g)\, dg = 1$. Then $(E_V)$ with the partial order of reverse inclusion of the $V$, is an approximate
identity in $L^1(G)$ and for each $B \in \cB$,
the relation
$\lim\limits_{g \to \1} \beta(g)B = B$ now implies that
$\lim\limits_{V \to \1}  \hat\beta(E_V)B = B$, and
we thus obtain $\cB \subeq \zeta(C^*(G))\cB(\cH)= \cB(\cH)_c^L$.

(c) For $t > 0$ there exists an
$f \in C_0(\R)$ with $f(x) = 1$ for $|x| \leq t$.
Then the operator $P(f)$ defined by the spectral
integral is contained in $\zeta(C^*(\R))$, and for $A \in \cB(\cH)$,
we have
\[ P([-t,t])A
=  P(f) P([-t,t])A \in \zeta(C^*(\R))\cB(\cH)
\subeq \cB(\cH)_c^L.\]
Conversely, $P([-t,t]) \in \cB(\cH)_c^L$
shows that $A = \lim_{t \to \infty} P([-t,t])A$
entails $A \in \cB(\cH)_c^L$.

If $A \in \cB(\cH)$ satisfies
$A \zeta(C^*(\R)) \subeq \cB(\cH)_c^L$
and $t > 0$, then we choose $f \in C_0(\R)$
with $f(x) = 1$ for $|x| \leq t$. Then
$A P([-t,t]) \in A P(f) P([-t,t])  \in \cB(\cH)_c^L$
because $\cB(\cH)_c^L$ is a right ideal.
If, conversely, $A P([-t,t])\in \cB(\cH)_c^L$
holds for every $t > 0$ and
$f \in C_c(\R)$ is supported in $[-t,t]$,
then $A P(f) = A P([-t,t]) P(f) \in \cB(\cH)_c^L$,
and we conclude that $A \zeta(C^*(G)) \subeq \cB(\cH)_c^L$.
\end{prf}

The following theorem provides a list of characterizations of
non-degeneracy for homomorphisms $\zeta \: \cS \to M(\cT)$
(cf.\ Definition~\ref{def:1.1c}(ii)).
It will be a crucially important tool for our constructions.

\begin{thm}\mlabel{thm:a.2} For $C^*$-algebras $\cS$ and $\cT$
and a morphism  $\zeta \: \cS \to M(\cT)$ of $C^*$-algebras,
the following are equivalent:
\begin{description}
\item[\rm(i)] $\zeta$ is non-degenerate.
\item[\rm(ii)] $\zeta(\cS)\cT\zeta(\cS) = \cT$.
\item[\rm(iii)] $\zeta$ extends to a strictly continuous homomorphism
$\tilde\zeta \: M(\cS) \to M(\cT)$ of unital $C^*$-algebras.
\item[\rm(iv)] For any  approximate identity
$(E_i)_{i \in I}$ in $\cS$ we have that $\lim\zeta(E_i)T=T$
for all $T\in\cT$.
\item[\rm(v)] $\zeta(\cS)\cT$ is strictly dense in $\cT$.
\item[\rm(vi)] $\zeta(\cS)\cT$ is weakly  dense in $\cT$.
\item[\rm(vii)] For each $\pi \in \Rep(\cT,\cH)$, the representation
 $\tilde\pi \circ \zeta$ of $\cS$ on $\cH$ is non-degenerate.
\end{description}
\end{thm}

\begin{prf} (i) $\Rarrow$ (ii): By  \cite[Th.~5.2.2]{Pa94} we have
$\overline{\rm span}(\zeta(\cS)\cT)=\zeta(\cS)\cT\subseteq\cT$, so as
$\zeta(\cS)\cT$ is dense, then $\cT=\zeta(\cS)\cT$. Thus
$\cT=\cT^*=\cT \zeta(\cS)$ and so
$\cT=\zeta(\cS)\cT\zeta(\cS)$.

(ii) $\Rarrow$ (iii): Clearly, (ii) implies (i).
Now (iii) follows from
\cite[Prop.~10.3]{Ne08}, which asserts the
existence of a unique extension $\tilde\zeta  \: M(\cS) \to M(\cT)$,
and that $\tilde\zeta$ is strictly continuous.

(iii) $\Rarrow$ (iv): For any  approximate identity
$(E_i)_{i \in I}$ in $\cS$ we have $E_i \to \1$ in the strict topology
of $M(\cS)$, so that we obtain for any $T \in \cT$ that
$\zeta(E_i) T = \tilde\zeta(E_i)T \to \tilde\zeta(\1)T = T$
via strict continuity of $\tilde\zeta$.

(iv)$\Rarrow$ (v):
This implies that $\zeta(\cS)\cT$ is dense in $\cT$.
In particular, this implies the density in the strict topology,
which is weaker.

(v) $\Rarrow$ (vi): We only have to show that
the weak topology on $\cT$ is weaker than the strict topology.
Suppose that $T_i \to T$ for the strict topology. Since
any functional $\phi \in \cT^*$ is a difference of two positive functionals,
it suffices to show that $\phi(T_i) \to \phi(T)$ holds for
any state $\phi \in \fS(\cT)$.
Let $\pi \: \cT \to \cB(\cH_\phi)$ be the
GNS-representation of $\phi$, and $v \in \cH_\phi$ with
$\phi(T) = \la \pi(T)v,v\ra$ for $T \in \cT$.
Then $\pi$ is non-degenerate, hence extends to a
 representation $\tilde\pi$ of $M(\cT)$ which is continuous
with respect to the strict topology on $M(\cT)$ and the
strong topology on $\cB(\cH)$ (\cite[Prop.~10.4]{Ne08}).
It follows in particular
that $\pi(T_i)v \to \pi(T)v$, and hence that
$\phi(T_i) \to \phi(T)$.

(vi) $\Rarrow$ (vii): For
$\pi \in \Rep(\cT,\cH)$, let $v \in \cH$ and assume that
$\tilde\pi(\zeta(\cS))v = \{0\}$. We have to show that
$v = 0$. Let $\phi(T) := \la \pi(T)v,v\ra$ be the corresponding
positive functional on $\cT$. Then
\[ \phi(\zeta(\cS)\cT)
= \la \tilde\pi(\zeta(\cS))\pi(\cT)v,v\ra
= \la \pi(\cT)v,\tilde\pi(\zeta(\cS))v\ra = \{0\},\]
so that the weak density of $\zeta(\cS)\cT$ in $\cT$
implies that $\phi = 0$, and hence that $v = 0$ because
$\pi$ is non-degenerate.

(vii) $\Rarrow$ (vi): Let $(S_i)_{i \in I}$ be a symmetric
approximate identity in $\cS$. Since each element of $\cT^*$ is a difference
of positive functionals, it suffices to show that
$\phi(\zeta(S_i)T) \to \phi(T)$ holds for each $T \in \cT$
and each positive functional $\phi \in \cT^*$.

Let $\pi \: \cT \to \cB(\cH_\phi)$ be the
GNS-representation of $\phi$,  and $v \in \cH_\phi$ with
$\phi(T) = \la \pi(T)v,v\ra$ for $T \in \cT$.
Then $\pi$ is non-degenerate and (vi) implies that
$\beta := \tilde\pi \circ \zeta$ is non-degenerate.
This implies that $\beta(S_i)w \to w$ for every $w \in \cH$
because the set of all elements for which this is the case is a closed
subspace containing the dense subspace $\beta(\cS)\cH$.
We conclude that
\[ \phi(\zeta(S_i)T)
= \la \beta(S_i) \pi(T)v,v \ra
\to \la \pi(T)v,v \ra  = \phi(T).\]

(vi) $\Rarrow$ (i) follows from the general fact that a subspace
of a locally convex space is weakly dense if and only if it is
dense  (cf.~\cite[Thm~1.4, Ch.~V]{Co97}).
\end{prf}

\begin{rem} If $\zeta(\cS)$ is strictly dense in $M(\cT)$,
then, for each $T \in \cT$, the closure of $\zeta(\cS)T$
contains $T$. In particular, $\zeta$ is non-degenerate.
In this case it would make sense to call the pair
$(\cT,\zeta)$ a host algebra for $\cS$ because,
for each Hilbert space $\cH$, we obtain an injective map
\[ \zeta^* \: \Rep(\cT,\cH) \to \Rep(\cS,\cH), \quad
\pi \mapsto \tilde\pi \circ \zeta.\]
\end{rem}
\begin{lem} \mlabel{lem:a.3a}
Let $\cS$ be a $C^*$-algebra
and $(S_j)_{j \in J}$ be a approximate identity.
Then for every homomorphism $\pi \: \cS \to \cB(\cH)$
the net
$\pi(S_j)$ converges strongly to the projection onto the essential subspace
\[ \cH_{\rm ess} := \{ v \in \cH \mid \pi(\cS)v = \{0\}\}^\bot
=  \oline{\Spann(\pi(\cS)\cH)} = \pi(\cS)\cH. \]
\end{lem}

\begin{prf}  First we note that the last equality follows from
\cite[Th.~5.2.2]{Pa94}.

Since $(S_j)$ is bounded, the subspace
\[ \cH_c := \{ v \in \cH \mid \pi(S_j)v \to v \}
\subseteq \oline{\Spann(\pi(\cS)\cH)} = \pi(\cS)\cH\]
of $\cH$ is closed. Clearly, it is contained in
$\cH_{\rm ess}$. From $\pi(\cS)\cH \subeq \cH_c$ we therefore
obtain $\cH_c = \cH_{\rm ess}$. On the other hand, we have
 $\pi(\cS)v = \{0\}$ for $v \in \cH_c^\bot$.
This proves the lemma.
\end{prf}

\begin{prop}\mlabel{prop:a.3}  If
$\zeta \: \cS \to M(\cT)$ is a morphism of $C^*$-algebras
and $(E_j)_{j \in J}$ is an approximate identity in $\cS$,
then the following assertions hold:
\begin{description}
\item[\rm(i)] $\cT_\cS := \zeta(\cS) \cT\zeta(\cS)$
is a  $\cS$-biinvariant $C^*$-subalgebra of $\cT$.
It is maximal with respect to the property that the corresponding
morphism
\[ \zeta^\sharp \: \cS\to M(\cT_\cS), \quad \zeta^\sharp(S)T
:= \zeta(S)T\]
is non-degenerate.
\item[\rm(ii)] $\cT_\cS$ is a hereditary subalgebra of $\cT$.
\item[\rm(iii)] In $\cT^{**}$ the weak limit
$P := \lim \zeta(E_j)$ exists and satisfies
\[ \cT_\cS = P \cT^{**}P \cap \cT
\quad \mbox{ and } \quad P \cT^{**} \cap \cT = \zeta(\cS) \cT.\]
\end{description}
\end{prop}

\begin{prf} (i) Note that
$\cT_\cS=\{T \in \cT\mid
 \zeta(E_j) T \to T \quad \mbox{ and } \quad
T \zeta(E_j) \to T\}$.
From this description and the boundedness of $(E_j)$ it follows easily
that $\cT_\cS$ is a closed $*$-subalgebra of $\cT$.

Since $\cT_\cS$ is $\cS$-biinvariant, we also obtain a
homomorphism
\[ \zeta^\sharp \: \cS \to M(\cT_\cS), \quad \zeta^\sharp(S)T
:= \zeta(S)T\]
of $C^*$-algebras. The definition of $\cT_\cS$ shows that
$\zeta^\sharp(E_j) T \to T$ for $T \in \cT_\cS$,
so that $\cT_\cS = \zeta^\sharp(\cS)\cT_\cS.$
We conclude with Theorem~\ref{thm:a.2} that $\zeta^\sharp$
is non-degenerate. If $\cT_0\subset\cT$ is another
$\cS$-biinvariant $C^*$-subalgebra for which the corresponding
map $\zeta^\sharp \: \cS \to M(\cT_0)$ is nondegenerate,
then $\zeta^\sharp(E_j) T \to T$ for $T \in \cT_0$,
so that $\cT_0 = \zeta^\sharp(\cS)\cT_0=\zeta(\cS)\cT_0$. Then
$\cT_0 = \cT_0^*=\cT_0\zeta(\cS)$ hence
$\cT_0 =\zeta(\cS)\cT_0\zeta(\cS)\subset\zeta(\cS)\cT\zeta(\cS)=\cT_\cS $.
Thus $\cT_\cS $ is maximal in the claimed sense.

(ii) A subalgebra $\cT_0 \subeq \cT$ is hereditary if
$\cT_0 \cT \cT_0 \subeq \cT_0$ (\cite[Thm.~3.2.2]{Mu90}).
This condition is trivially satisfied for $\cT_\cS$:
\[ \cT_\cS \cT \cT_\cS
= \zeta(\cS) \cT \zeta(\cS)\cT \zeta(\cS) \cT \zeta(\cS)
\subeq \zeta(\cS) \cT \zeta(\cS) = \cT_\cS.\]

(iii) We consider the enveloping $W^*$-algebra
$\cT^{**}$ of $\cT$, which coincides, as a Banach space, with the bidual
of $\cT$ and note that $M(\cT)$ is a $C^*$-subalgebra of $\cT^{**}$.
We realize $\cT^{**}$ as a von Neumann algebra on some
Hilbert  space $\cH$. Then Lemma~\ref{lem:a.3a}
implies the existence of a projection
$P := \lim E_j \in \cT^{**}$ (weak limit) such that
the range of $P$ is the essential subspace
for $\cS$ in any non-degenerate representation of $\cT$.

For any $T \in \cT_\cS$ we now find that
$T P = P T = T$. Suppose, conversely, that $T \in \cT$ satisfies
$PTP = T$. We claim that $T \in \cT_\cS$.
As $\lim E_j T = PT = T$ holds $\sigma(\cT^{**}, \cT^*)$-weakly
in $\cT^{**}$, is also
holds $\sigma(\cT, \cT^*)$-weakly in $\cT$. This implies that
$T$ is contained in the weak closure of the norm closed subspace
\[ \zeta(\cS)\cT = \{ T \in \cT \mid \zeta(E_j) T \to T\}. \]
Since closed subspaces of Banach spaces are also weakly closed
(cf.~\cite[Thm~1.4, Ch.~V]{Co97}),
we obtain $T \in \zeta(\cS)\cT$. Likewise $TP = T$
implies $T \in \cT\zeta(\cS)$, so that
$T \in \zeta(\cS)\cT \zeta(\cS) = \cT_\cS$.
This proves the identity
\[ \cT_\cS = P \cT^{**} P \cap \cT.\]

To verify that
$P \cT^{**} \cap \cT = \zeta(\cS) \cT$, note that both
sides are closed right ideals of $\cT$ and
the corresponding hereditary subalgebras
\[ \zeta(\cS)\cT \cap (\zeta(\cS)\cT)^*
= \zeta(\cS)\cT \cap \cT \zeta(\cS)
= \zeta(\cS)\cT \zeta(\cS)\]
and
\[ (P \cT^{**} \cap \cT) \cap (P\cT^{**} \cap \cT)^*
= P \cT^{**} \cap \cT^{**}P \cap \cT
= P \cT^{**} P \cap \cT = \cT_\cS\]
coincide. This implies that
$P \cT^{**} \cap \cT = \zeta(\cS) \cT$
(\cite[Thm.~3.2.1]{Mu90}).
\end{prf}

\begin{rem}
$P$ is called the open
projection of  $\cT_\cS$ in Pedersen's terminology,
cf.~\cite[Prop.~3.11.9 and 3.11.10]{Pe89}.

The preceding proposition has an interesting consequence
for the connection between representations of $\cS$ and $\cT$.
Each non-degenerate representation $\pi \: \cT_\cS \to \cB(\cH)$
defines a non-degenerate representation
$\tilde\pi \circ \zeta^\sharp$ of $\cS$
(Theorem~\ref{thm:a.2}), but since every state of $\cT_\cS$
has a unique extension to $\cT$ (Lemma~\ref{lem:a.2.1}),
every cyclic representation
of $\cT_\cS$ embeds into a cyclic representation of $\cT$.

Conversely, if $v$ is a cyclic vector for a representation $\pi$ of $\cT$,
then the closure of
\[ \pi(\cT_\cS)v
= \tilde\pi(\zeta(\cS))\pi(\cT)\tilde\pi(\zeta(\cS))v
= P \tilde\pi(\zeta(\cS))\pi(\cT)\tilde\pi(\zeta(\cS))Pv
\subeq P \pi(\cT)Pv \]
is the representation space of a cyclic representation of $\cT_\cS$.
\end{rem}

\begin{lem} A representation
$\pi \: \cA \to \cB(\cH)$ of a $C^*$-algebra $\cA$ is non-degenerate
if and only if it is non-degenerate as a homomorphism
$\cA \to M(\cK(\cH))$.
\end{lem}

\begin{prf} If the representation of $\cA$ on $\cH$ is non-degenerate
and $(A_i)$ is an approximate identity in $\cA$,
then $\pi(A_i)v \to v$ holds for every $v \in \cH$ and hence
$\pi(A_i)F \to F$ for every finite rank operator $F$ on $\cH$.
Since the finite rank operators form a dense subspace of $\cK(\cH)$,
it follows that $\pi(\cA)\cK(\cH)$ is dense in $\cK(\cH)$.

If, conversely, $\pi(\cA)\cK(\cH)$ is dense in $\cK(\cH)$ and
$v \in \cH$, then $\pi(\cA)v = \{0\}$ implies
\[ \la \pi(\cA)\cK(\cH)v, v\ra = \{0\}\] and hence that
$\la \cK(\cH)v,v\ra = \{0\}$, which leads to $v = 0$. This means that
$\pi$ is a non-degenerate representation.
\end{prf}

\begin{lem} Let $G$ be a locally compact group, $\cS := C^*(G)$,
$\cT$ be a $C^*$-algebra and
$\eta_G \: G \to \U(M(\cT))$ be a homomorphism.
We write $\cT_c$ for the closed right ideal of $\cT$ consisting
of all elements $T$ for which the $G$-orbit map
$G \to \cT, g \mapsto \eta_G(g)T$ is continuous.
Then there exists a non-degenerate homomorphism
$\eta_\cS \: \cS \to M(\cT)$ which is compatible
with $\eta_G$ in the sense that
\begin{equation}
  \label{eq:cstarcond}
\eta_G(g) \eta_\cS(S)T = \eta_\cS(gS)T \quad \mbox{ for } \quad
S \in C^*(G), T \in \cT
\end{equation}
if and only if $\cT = \cT_c$.
\end{lem}

\begin{prf} Since the canonical homomorphism
$\eta \: G \to M(C^*(G))$ is strictly continuous
($G$ acts continuously by translations on $L^1(G)$),
\eqref{eq:cstarcond} implies that
$\eta_\cS(\cS)\cT \subeq \cT_c.$
We conclude that $\cT = \cT_c$ whenever $\eta_\cS$ is non-degenerate.

If, conversely, $\cT = \cT_c$, then the strong continuity of the
$G$-action $(g,T) \mapsto \eta_G(g)T$ on $\cT$ leads to a
morphism $L^1(G) \to M(\cT), f \mapsto \int_G f(g)\eta_G(g)\, dg$
of Banach-$*$-algebras and the universal property of $C^*(G) = \cS$
further implies the existence of an extension
$\eta_\cS \: \cS \to M(\cT)$
satisfying \eqref{eq:cstarcond}. If $(\delta_j)_{j \in J}$
is an approximate identity in $L^1(G)$, then the continuity of
the $G$-action on $\cT$ implies $\eta_G(\delta_j)T \to T$
for every $T \in \cT$, and therefore $\eta_\cS$ is non-degenerate.
\end{prf}

\section{Functionals and hereditary $C^*$-subalgebras} \mlabel{app:2}

\begin{lem} \mlabel{lem:a.2.1}
Suppose that $\cB$ is a hereditary subalgebra of the $C^*$-algebra
$\cA$. Then each positive functional
$\phi$ on $\cB$ has a unique positive extension
$\tilde\phi$ to $\cA$ with $\|\tilde\phi\| = \|\phi\|$,
and this assignment extends to a linear embedding
\[ \cB^* \into \cA^*, \quad \phi \mapsto \tilde\phi\]
with $\tilde\phi\res_\cB = \phi$ for $\phi \in \cB^*$ which is
isometric on the real subspace of selfadjoint functionals.
\end{lem}

\begin{prf} Any selfadjoint $\phi \in \cB^*$ (the topological dual)
has a unique (Jordan) decomposition $\phi = \phi_+ - \phi_-$ with
$\phi_\pm$ positive and
$\|\phi\| = \|\phi_+\| + \|\phi_-\|$ (cf.\ \cite[Thm.~3.2.5]{Pe89}).

According to \cite[Thm.~3.3.9]{Mu90}, each positive functional
$\phi \in \cB^*$ has a unique extension $\tilde\phi$ to $\cA$
with the same norm. Furthermore, if
$(B_j)_{j \in J}$ is an approximate identity
in $\cB$, then
\[ \tilde\phi(A) = \lim_{j \in J}
\phi(B_jAB_j) \quad \mbox{ for } \quad  A\in \cA.\]
Since the positive functionals span $\cB^*$, it follows that the
limit on the right hand side exists for every $\phi \in \cB^*$ and
every $A \in \cA$. Moreover, for the Jordan decomposition
$\phi = \phi_+ - \phi_-$, we obtain
\[ \tilde\phi = (\phi_+)\,\tilde{} - (\phi_-)\,\tilde{} \]
with $(\phi_\pm)\,\tilde{}$ positive and
$\|(\phi_\pm)\,\tilde{}\| = \|\phi_\pm\|.$
This implies in particular that $\tilde\phi$ is continuous with
\[ \|\tilde\phi\| \leq
\|(\phi_+)\,\tilde{}\| + \|(\phi_-)\,\tilde{}\|
= \|\phi_+\| + \|\phi_-\| = \|\phi\|.\]
As $\tilde\phi$ extends $\phi$, we also have
$\|\phi\| \leq \|\tilde\phi\|$ for trivial reasons, hence equality.
We also derive from the uniqueness of the Jordan decomposition
that $\tilde\phi_\pm = (\phi_\pm)\,\tilde{}$.
\end{prf}

\section{Regularity properties of unitary group representations} \mlabel{app:3}

\begin{defn} \mlabel{def:rep-concepts}
(a) For a continuous unitary representation $(U,\cH)$ of the compact group $G$,
we write $\cH = \bigoplus_{\pi \in \hat G} \cH_\pi$ for the isotypic decomposition,
i.e., the representation on the subspace $\cH_\pi$ is a multiple of the
irreducible representation $\pi$.
The
\[ \Spec(U) := \{ \pi \in \hat G \: \cH_\pi \not=\{0\} \} \]
is called the {\it spectrum of $U$}.
In other words, $\Spec(U)$ is the set of
equivalence classes of the irreducible unitary
representations  which are contained in the direct sum decomposition of $U$.
We say that $U$ is {\it of finite multiplicity} if, for every $V \in \hat G$, the
corresponding multiplicity $\dim \Hom_G(V,\cH)$ is finite.

(b) Let $G$ be a locally compact abelian group and $\hat G= \Hom(G,\T)$ be its
character group, endowed with the topology of uniform convergence on compact subsets of
$G$. For a continuous unitary representation $(U,\cH)$ of $G$,
we define its {\it spectrum} $\Spec(U) \subeq \hat G$ as the support of the corresponding
spectral measure $P$, i.e., as the set of all those characters $\chi \in \hat G$ with
the property that every open neighborhood $U_\chi$ of $\chi$ in $\hat G$ satisfies
$P(U_\chi)\not=0$.
\end{defn}

\subsection{Abelian groups}

\begin{lem} \mlabel{lem:c.2} A continuous unitary representation of
the locally compact abelian group $G$ is norm continuous if and only
if $\Spec(U) \subeq \hat G$ is compact.
\end{lem}

\begin{prf} Using the spectral measure $P$ of $U$ and the notation
$\hat g(\chi) := \chi(g)$ for $g \in G, \chi \in \hat G$, it follows that the
representation $U$ can be written as
$U_g = \hat U(\hat g)$, where
\[ \hat U \: C_0(\Spec(U)) \to \cB(\cH), \quad f \mapsto P(f) = \int_{\Spec(U)} f(\chi)\, dP(\chi) \]
is given by the spectral integral. Since $\hat U$ is an isometric embedding, it suffices
to show that, for a closed subset $\Sigma \subeq \hat G$, the map
\[ \beta \: G \to C_0(\Sigma), \quad  g \mapsto \beta(g)(\chi)= \hat g(\chi) = \chi(g) \]
is norm continuous if and only if $\Sigma$ is compact.
As
$\|\hat g\res_\Sigma - 1\|_\infty = \sup_{\chi \in \Sigma} |\chi(g) - 1|,$
the norm continuity of $\beta$ is equivalent to the equicontinuity of $\Sigma$.
According to \cite[Prop.~7.6]{HoMo98},  the equicontinuity
of the closed subset $\Sigma \subeq \hat G$ is equivalent to its compactness.
\end{prf}

\begin{lem} \mlabel{lem:c.5}
Let $(U,\cH)$ be a continuous unitary
representation of the locally compact abelian group $G$ and
$U_{C^*(G)} \: C^*(G) \to \cB(\cH)$ be the
associated integrated representation. Then
\[
U_{C^*(G)}(C^*(G))\subseteq\cK(\cH)
\]
if and only if the spectral measure $P$ of $U$ is a locally finite sum of point measures
with finite-dimensional ranges.

For $G = \R$ and $U_t = e^{itA}$, this condition is equivalent to the compactness of the
resolvent $(A + i \1)^{-1}$.
\end{lem}

\begin{prf} Suppose first that $U_{C^*(G)}(C^*(G)) \subeq \cK(\cH)$.
Recall that, for each selfadjoint compact operator, all eigenspaces
for eigenvalues $\lambda\not=0$ are finite-dimensional.
As $U_{C^*(G)}(C^*(G)) $ acts nondegenerately, the intersection of the zero-eigenspaces
of all selfadjoint elements of $U_{C^*(G)}(C^*(G)) $  is zero.
It follows easily that $\cH$ is an orthogonal sum of
$U(G)$-invariant finite-dimensional subspaces.
Since $G$ is abelian,
$\cH$ is a direct sum $\hat\oplus_{\chi \in \hat G} \cH_\chi$ of the $G$-eigenspaces.
Now $C^*(G) \cong C_0(\hat G)\ni f$ acts in this decomposition by multiplication
$U_{C^*(G)}(f)(v_\chi)= f(\chi)v_\chi$, $v_\chi\in\cH_\chi$, so the compactness of all these operators implies
that all projections $P(\{\lambda\})$ have finite-dimensional range and that,
for each compact subset $C \subeq \hat G$, at most finitely many spaces
$\cH_\chi$, $\chi \in C$, are non-zero. Therefore $P$ is a locally finite sum
of finite rank projections.

Suppose, conversely, that $P$ is a locally finite sum of finite rank projections.
Let $f \in C_0(\hat G) \cong C^*(G)$ and $\eps > 0$.
Then $E := \{ \chi \in \hat G \: |f(\chi)| \geq  \eps\}$ is compact, so that
there exist $\lambda_1, \ldots, \lambda_k \in \hat G$ with
$P(E) =P(\{\lambda_1\}) + \cdots + P(\{\lambda_k\})$.
Now $P(E)$ is a finite rank projection commuting with
$U_{C^*(G)}(f)$ and $\| U_{C^*(G)}(f)(\1-P(E))\| \leq \eps$.
Therefore $U_{C^*(G)}(f)$ is compact.

For $G = \R$, the image of $C^*(\R) \cong C_0(\R)$ is generated by the resolvent operator
$(A + i \1)^{-1}$. Therefore its compactness is equivalent to the compactness of the image of
$C^*(\R)$.
\end{prf}

\begin{rem} For a  continuous unitary representation
$(U,\cH)$  of $\R$, we have the following description of the
$C^*$-algebra of continuous vectors $\cB(\cH)_c \subeq \cB(\cH)$
for the conjugation representation
$\alpha_t(A) = U_t A U_t^*$.

Let $P$ denote the spectral measure of $U$.
For $n \in \Z$, let $P_n:= P([n,n+1))$ and  $\cH_n := P_n\cH$, thus
we  obtain an orthogonal decomposition $\cH = \hat\oplus_{n \in \Z}\cH_n$.
We use this decomposition
for a convenient factorization of $U$. We write
$U_t = U^1_t U^2_t$, where
$U^1_t v = e^{int}v$ for $v \in \cH_n$. Then
$U^1$ is a continuous unitary representation with
$U^1_{2\pi} = \1$, so that it actually defines a representation
of the circle group $\T \cong \R/2\pi \Z$. The representation
$U^2_t = U_t U^1_{-t}$ has spectral measure supported by $[0,1]$,
so that it is norm-continuous. Therefore the two representations
$U$ and $U^1$ define the same space $\cB(\cH)_c$ of continuous elements for the
conjugation action.

Let $\alpha^1:\R\to\Aut\cB(\cH)$ be the conjugation $\alpha_t^1(A) = U_t^1 A U^1_{-t}$, and denote by
\[ \cB(\cH)_n := \{ A \in \cB(\cH)\mid (\forall t \in \R)\
\alpha_t^1(A)  = e^{int} A\}\]
its $\T$-eigenspaces in $\cB(\cH)$. Then $\cB(\cH)_n\subset \cB(\cH)_c$, hence
$\cB(\cH)_c\cap\cB(\cH)_n=\cB(\cH)_n$.
The Peter--Weyl Theorem generalizes to continuous Banach-space actions of $G$
(cf.~\cite[Th.~2]{Sh55} and~\cite[Th.~3.51]{HoMo98}), hence an application of it
to $\alpha^1\restriction\cB(\cH)_c$ implies that
\begin{equation}
\label{BHc}
 \cB(\cH)_c = \oline{\Spann\Big(\mathop{\cup}_{n \in \Z} \cB(\cH)_n\Big)}.
 \end{equation}
Write $A = (A_{jk})_{j,k\in \Z}$ as a matrix with
$A_{jk} \in B(\cH_k, \cH_j)$, and keep in mind that the convergence
$A=\sum\limits_{j\in\Z}\sum\limits_{k\in\Z}A_{jk}=\sum\limits_{j\in\Z}\sum\limits_{k\in\Z}P_jAP_k$ is in general
w.r.t. the strong operator topology.
We have
\[ \alpha_t^1(A) = (e^{it(j-k)}A_{jk})_{j,k\in \Z},\]
so that
\[ A \in \cB(\cH)_n \quad \Longleftrightarrow \quad
(j-k \not=n \Rarrow A_{jk} =0).\]
For $A = (A_{jk})_{j,k\in \Z} \in \cB(\cH)$, let
$A_n := (A_{jk} \delta_{j-k,n})_{j,k\in \Z}$ and observe that $A_n$ defines a bounded
operator on $\cH$, hence an element of $\cB(\cH)_n$.
If $A = \sum\limits_{n \in \Z} A_n$ converges w.r.t. norm topology, it follows that
$A \in \cB(\cH)_c$, but this condition is not necessary for
$A \in \cB(\cH)_c$ because Fourier series of continuous functions need not converge in norm.
By (\ref{BHc}), the space $\cB(\cH)_c$ is the norm closure of the space
$\Spann\big(\mathop{\cup}\limits_{n \in \Z} \cB(\cH)_n\big)$
which consists of those operators
whose matrix has only finitely
many non-zero diagonals.
\end{rem}

\subsection{Compact groups}

\begin{prop} \mlabel{prop:c.8}
For a continuous unitary representation $(U,\cH)$ of a compact group $G$,
the image $U_{C^*(G)}(C^*(G))$ of the integrated representation consists of compact
operators if and only if $U$ is of finite multiplicity.
\end{prop}

\begin{prf}  By Proposition~3.4 in~\cite{Wil07},
we have $C^*(G)\cong\mathop{\bigoplus}^\infty\limits_{V\in\hat G}
 \cB(\cH_V)$ and by \cite[Thm.~38.3]{HR97},
the (finite-dimensional) summand $\cB(\cH_V)$ is the minimal ideal $\cI_V$ of $C^*(G)$
 corresponding to $V\in\hat G$ in the sense that $V_{C^*(G)}$ is faithful on $\cI_V$ and zero on each
 summand $\cI_{V'}$ if $V\not=V'\in\hat G$.
 Then $U_{C^*(G)}(\cI_V) = \cB(\cH_V) \otimes \1_{\cM_V}$, where
$\cM_V := \Hom_G(V,\cH)$ denotes the corresponding multiplicity space.
This space consists of compact operators if and only if the multiplicity of
$V$ in $U$ is finite. Since
$\Spann\{ \cI_V\mid V\in \hat G\}$ is dense in $C^*(G)$, this proves the
assertion.
\end{prf}

\begin{lem} \mlabel{lem:c.9} For a compact group $G$, a unitary representation
$U \: G \to \U(\cH)$ is norm continuous if and only if its spectrum is finite.
\end{lem}

\begin{prf} Since each irreducible representation of $G$ is finite-dimensional,
hence norm continuous, this property is also inherited by any representation
with finite spectrum.

Suppose, conversely, that $U$ is norm continuous. Then
its range $U(G) \subeq \U(\cH)$ is a subgroup which is compact in the norm
topology, hence a Lie group by Yosida's Theorem (\cite{Yo36}). We may therefore
assume that $G$ is a Lie group. Let $T \subeq G$ be a maximal torus.
Then the boundedness of $U\res_T$ implies that the spectrum of this representation
in the discrete character group $\hat T$ is bounded (Lemma~\ref{lem:c.2}), hence finite.
Therefore the  classification of the irreducible representations
of the identity component $G_0$ in terms of dominant weights in $\hat T$ implies that
the $G_0$-spectrum of $U$ is finite because $\Spec(U\res_T)$ is finite and it contains
the ``highest weights'' of the irreducible $G_0$-representations occurring in~$U$.

For any irreducible $G$-subrepresentation $V$ of $U$ there exists an irreducible
$G_0$-subrepresentation $V_0$ of~$V$. Then
\[ \{0\} \not = \Hom_{G_0}(V_0, V\res_{G_0}) \cong \Hom_G(\Ind_{G_0}^G(V_0), V) \]
implies that $V$ is contained in the induced representation $\Ind_{G_0}^G(V_0)$.
Since $G_0$ is of finite index in $G$, the induced representation decomposes into finitely
many irreducible $G$-representation. Therefore the finiteness of $\Spec(U_{G_0})$ implies
the finiteness of $\Spec(U)$.
\end{prf}

\section*{Acknowledgments}

Hendrik Grundling gratefully acknowledges support by the
Emerging Fields Project ``Quantum Geometry'' at the Universit\" at Erlangen--N\"urnberg.
He is also very grateful to Detlev Buchholz, who pointed out a mistake in an earlier version
of Example~9.3.

\end{document}